\documentclass[12pt]{amsart}
\usepackage[utf8]{inputenc}
\usepackage[T1]{fontenc}
\usepackage{float}
\usepackage{nicematrix}
\usepackage{amsmath,amssymb,amsthm}
\usepackage[hidelinks]{hyperref}
\usepackage{enumitem}
\usepackage{tikz}
\usetikzlibrary{arrows.meta}
\usetikzlibrary{calc}
\usetikzlibrary{backgrounds}
\usepackage[margin=2cm]{geometry}

\theoremstyle{plain}
\newtheorem{theorem}{Theorem}[section]
\newtheorem{proposition}[theorem]{Proposition}
\newtheorem{lemma}[theorem]{Lemma}
\newtheorem{corollary}[theorem]{Corollary}
\theoremstyle{definition}
\newtheorem{definition}[theorem]{Definition}
\newtheorem{example}[theorem]{Example}
\newtheorem{remark}[theorem]{Remark}
\newcommand{\dcon}{\partial_{\mathrm{con}}}
\newcommand{\kk}{\Bbbk}
\newcommand{\F}{ {F_\bullet} }
\newcommand{\Rbar}{\overline{R}}
\newcommand{\M}{\mathcal{M}}
\DeclareMathAlphabet{\mathpzc}{OT1}{pzc}{m}{it}
\newcommand{\lift}[1]{\mathpzc{#1}}

\newcommand{\A}{\mathcal{A}}
\newcommand{\B}{\mathcal{B}}
\newcommand{\Grep}{G^{\mathrm{rep}}}
\newcommand{\id}{\operatorname{id}}
\newcommand{\lcm}{\operatorname{lcm}}
\newcommand{\hdeg}{\operatorname{hdeg}}
\newcommand{\gr}{\operatorname{gr}}
\newcommand{\pd}{\operatorname{pd}}
\newcommand{\aaa}{\mathfrak{a}}
\newcommand{\N}{\mathbb{N}}
\newcommand{\Sh}{\mathrm{Sh}}
\newcommand{\eps}{\varepsilon}
\newcommand{\wt}{\operatorname{wt}}

\newcommand{\cbl}{}

\colorlet{sh0}{cyan!10}
\colorlet{sh1}{yellow!28}
\colorlet{sh2}{violet!10}
\newcommand{\shent}[2]{{\setlength{\fboxsep}{1.5pt}\colorbox{#1}{$\vphantom{\sigma_{10}}#2$}}}

\title[Pruning Free Resolutions over a Complete Intersection]{Pruning Free Resolutions over a Complete Intersection}
\author{Benjamin Betts}
\author{Janet Vassilev}
\date{}
\begin{document}
\begin{abstract}
Given a free resolution $\F$ over a polynomial ring $R$ and a system of higher homotopies on it for a regular sequence $a_1,\dots,a_c$, we can use the Eisenbud--Shamash construction to build a free resolution over $\bar{R}=R/(a_1,\dots,a_c)$. We provide a combinatorial formula for the new system of higher homotopies on a pruned free resolution $\F^\M$ that corresponds to a pruning of the Eisenbud--Shamash resolution.
\end{abstract}

\maketitle

\section*{Introduction}

Free resolutions over a complete intersection are typically infinite, and the standard way to build one, due to Shamash \cite{Shamash} and Eisenbud \cite{Eisenbud1980}, starts from a finite free resolution over a polynomial ring. One such free resolution is the Taylor resolution \cite{Taylor} of a monomial ideal over a polynomial ring. In \cite{Sobieska}, Sobieska gives a compact formula for a system of higher homotopies on the Taylor resolution, which gives a free resolution over a complete intersection using the construction of Eisenbud and Shamash.
Unfortunately, Taylor resolutions are rarely minimal. Over the polynomial ring, algebraic Morse theory is a standard tool for pruning a resolution down toward the minimal one, although sometimes there are obstructions to achieving minimality. There are two main algorithms which produce Morse matchings on the Taylor resolution of a monomial ideal: the pruning algorithm \cite{AlvarezFernandezGimenez} and the Barile--Macchia algorithm \cite{BarileMacchia,ChauKara}.
Several authors have exploited these algorithms to produce minimal resolutions for monomial ideals with five or fewer generators \cite{Faridi4gen,Chau5gen}, for edge ideals of certain graphs \cite{BarileMacchia,ChauKara,Chauetalpath,ChauHagenBarileMacchia}, for powers of monomial ideals which have minimal pruned resolutions \cite{CooperetalMorseres,Cooperetalsimplicialres,ChauHaMinCellpowersofgraphs} or for ideals where the Morse matching produces DG structure \cite{CameronChauMaitraTribone,GelleretalDGprune}.
We show that algebraic Morse theory can also be used to prune these infinite resolutions over complete intersections.

Algebraic Morse theory is due independently to Sk\"oldberg \cite{Skoldberg} and to J\"ollenbeck and Welker \cite{JollenbeckWelker}, building on Forman's discrete Morse theory \cite{Forman}, its reformulation through acyclic matchings by Chari \cite{Chari}, and its use on cellular resolutions by Batzies and Welker \cite{BatziesWelker}. The \emph{complex graph} $G(\F)$ has the basis elements of $\F$, the \emph{cells}, as vertices and an edge $\alpha \to \beta$ whenever the coefficient $d_{\alpha,\beta}$ of $\beta$ in $d(\alpha)$ is nonzero. A \emph{Morse matching} $\M$ is a set of edges with unit coefficient, no two sharing a cell, such that reversing them leaves no infinite directed path. The \emph{Morse complex} $\F^{\M}$ is free on the unmatched, \emph{critical}, cells, with differential given by sums over paths, and a contraction $(i^{\M}, p^{\M}, h^{\M})$ which makes it a deformation retract of $\F$ (Theorem~\ref{thm:amt}).

Let $R$ be a polynomial ring over a field, $\F$ a free resolution of $R/I$ over $R$ with differential $d$, and $a_1, \dots, a_c$ a regular sequence contained in $I$, with $\Rbar = R/(a_1, \dots, a_c)$. A \emph{system of higher homotopies} $\sigma$ is a family of maps $\sigma_{\boldsymbol{u}}$ on $\F$, for $\boldsymbol{u} \in \N^c$, raising homological degree by $2|\boldsymbol{u}| - 1$, with $\sigma_{\boldsymbol{0}} = d$ and $\sum_{\boldsymbol{v} + \boldsymbol{v'} = \boldsymbol{u}} \sigma_{\boldsymbol{v}} \sigma_{\boldsymbol{v'}}$ equal to $a_j \id$ if $\boldsymbol{u} = \boldsymbol{e}_j$ and $0$ otherwise. It is \emph{truncated} if $\sigma_{\boldsymbol{u}} = 0$ for $|\boldsymbol{u}| \ge 2$. The \emph{Eisenbud--Shamash resolution} $\Sh(\F, \sigma)$ is free over $\Rbar$ on the cells $y^{\boldsymbol{u}} \otimes \alpha \otimes\overline{1}$, for $\boldsymbol{u} \in \N^c$ and $\alpha$ in the basis of $\F$, with differential $\sum_{\boldsymbol{u}} t^{\boldsymbol{u}} \otimes \sigma_{\boldsymbol{u}} \otimes \id$, where the \emph{CI operator} $t_j$ sends $y^{\boldsymbol{u}}$ to $y^{\boldsymbol{u} - \boldsymbol{e}_j}$. {\cbl We will spend Sections~\ref{sec:shamash} to~\ref{sec:sobieska} reviewing the Eisenbud--Shamash construction at a level of detail that we have not found elsewhere in the literature.}

For a truncated system $\sigma = \{\sigma^{(j)}\}_{j=1}^c$, with $\sigma^{(j)} = \sigma_{\boldsymbol{e}_j}$, and a Morse matching $\M$ on $\F$, let $\lift{M}$ be one copy of $\M$ over each $\boldsymbol{u} \in \N^c$. The main result (Theorem~\ref{thm:commute}) is that, under the identity on the common basis,
\[
    \Sh(\F, \sigma)^{\lift{M}} = \Sh\bigl(\F^{\M}, \sigma^{\M}\bigr), \qquad \sigma^{\M}_{\boldsymbol{u}} = \sum_{\boldsymbol{e}_{j_1} + \dots + \boldsymbol{e}_{j_k} = \boldsymbol{u}} p^{\M} \sigma^{(j_1)} h^{\M} \sigma^{(j_2)} \cdots h^{\M} \sigma^{(j_k)} i^{\M},
\]
and that $\sigma^{\M}$ is again a system of higher homotopies, in general no longer truncated.

The proof applies the homological perturbation lemma, which goes back to Brown \cite{Brown1965} and Gugenheim \cite{Gugenheim1972} (see Crainic \cite{Crainic} for a modern account). Sk\"oldberg \cite{SkoldbergHPT} derives algebraic Morse theory from it for a filtered complex. We replace his filtration by the assumption that the relevant geometric series is finite on each element (Lemma~\ref{lem:perturbation}), and we apply the lemma twice with an additivity property of perturbations (Lemma~\ref{lem:additive}). We split the Shamash differential into the part along each copy of $\F$ and the connecting part $\dcon$ between copies, and perturb first by the relevant perturbation at each level and then by the connecting part. Lemma~\ref{lem:additive}, the Perturbation Additivity Lemma, is the main tool to show that performing two perturbations produces the same resolution as doing the single perturbation by the sum of the perturbations. Lemma~\ref{lem:additive} is a version of a result of Real \cite[Proposition~3.6]{Real2000} under different hypotheses on the second perturbation, and our main theorem puts it to significant use. 

The paper is structured as follows. We fix notation in Section~\ref{s:Notation} for easy reference. In Section~\ref{sec:amt}, we introduce algebraic Morse theory and commence our running example with the Taylor resolution in Example \ref{ex:runtaylor} and a pruned resolution in Example \ref{ex:runSec3}.
In Section~\ref{sec:shamash}, we recall the Eisenbud--Shamash resolution. In Section~\ref{sec:shdiff}, we continue with some insights on how to decompose the Shamash differentials and view them as block matrices. In Section \ref{sec:sobieska} we continue our running example with the Taylor resolution in Example \ref{ex:runtaylortruncated} and recall Sobieska's homotopies for the Taylor resolution of a complete intersection.
 In Section~\ref{sec:hyp}, we define lifted matchings to the Eisenbud--Shamash resolution for any matching on the Taylor resolution; we prove that these lifted matchings are Morse when the original matching is. In Section \ref{sec:examples} as preparation for our Main Theorem \ref{thm:commute}, we present some motivating examples. In Section~\ref{sec:commute}, we prove the main theorem, Theorem \ref{thm:commute}.  The paper concludes with Section \ref{sec:moreexamples}, where we introduce new examples in which the transferred higher homotopies may no longer be truncated even though those of the Taylor resolution were.

\subsection*{Acknowledgments}
We thank Luigi Ferraro for the DG insight behind Sobieska's homotopy formula, and Josep \`Alvarez Montaner, Margherita Barile and Olivia Strahan for insightful conversations about matching algorithms.

\section{Notation} \label{s:Notation}

Throughout, $A$ will be a commutative Noetherian ring, $R=\kk[x_1, \ldots, x_n]$ will be a polynomial ring over a field $\kk$ and $\overline{R}=R/\aaa$ where $\aaa=(a_1, \ldots, a_c)$ and $a_1, \ldots, a_c$ forms a regular sequence in $R$.

This paper pertains to complexes of finite free $A$-modules: $(F_\bullet, d_\bullet)$ where $F_q$ is free of rank $s_q$ for some $s_q \in \N$ with basis $\B_q$ and $d_q: F_q \rightarrow F_{q-1}$ can be represented by an $s_{q-1} \times s_q$ matrix with an appropriate choice of basis. We call elements $x \in F_q$, elements of \emph{homological degree} $q$ and we will write $\hdeg x=q$. The basis for the complex is $\B = \bigcup_q \B_q$ so that $d$ has matrix representation from $\B$.

For free complexes $(F_\bullet, d)$ and $(F'_\bullet, d')$, we will consider various collections of maps between them, which sometimes have adornments such as subscripts or superscripts. Because these adornments are important descriptors of these maps, we will use the convention that if $f :F_\bullet \rightarrow F'_\bullet$ is such a collection of maps,
possibly not a chain map,
possibly of nonzero external degree, $f|_q$ will refer to the map $f_\bullet$ which starts at the homological degree $q$ in $F_\bullet$. Sometimes, however, it will be clear from the representation of our elements in $F_q$ in terms of its basis elements which homological degree the map starts in; in this case we will drop $|_q$.

We use the standard fine multigrading for elements in a polynomial ring, $R=\kk [x_1, \ldots, x_n]$.
The \emph{homogeneous elements} are the terms $\lambda x_1^{b_1}\cdots x_n^{b_n}$, where $\lambda \in \kk$ and $(b_1, \ldots, b_n) \in \N^n$. If $r$ is a homogeneous element of $R$ of multidegree ${\boldsymbol{b}}=(b_1, \ldots, b_n)$, we will write $\gr(r)={\boldsymbol{b}}$. $\N^n$ is endowed with the componentwise partial order, so $\boldsymbol{b} \le \boldsymbol{b'}$ means $b_i \le b_i'$ for every $1 \le i \le n$, and $\boldsymbol{b} < \boldsymbol{b'}$ means $\boldsymbol{b} \le \boldsymbol{b'}$ with $\boldsymbol{b} \neq \boldsymbol{b'}$. For $\boldsymbol{u} \in \N^n$ we write $|\boldsymbol{u}| = u_1 + \dots + u_n$, and $\boldsymbol{e}_1, \dots, \boldsymbol{e}_n$ are the standard basis vectors of $\N^n$. Sometimes we will represent the multidegree of an element by the monomial $x^{\boldsymbol{b}}=x_1^{b_1}\cdots x_n^{b_n}$ instead of by $\boldsymbol{b}$. Then $\boldsymbol{b} \leq \boldsymbol{b}'$ if and only if $x^{\boldsymbol{b}} \mid x^{\boldsymbol{b}'}$.

We have compiled a list of notation we use in the paper and where we first reference it.

\subsection*{List of notation}
\begin{itemize}[leftmargin=*,itemsep=0pt,parsep=0pt]
    \item $(\F, d)$: free complex with basis
     $\B=\bigcup_n \B_n$ where $\B_n$ is the basis of $F_n$. {\color{blue} (Section \ref{s:Notation}, Definition \ref{def:based})}
    \item $f|_q$: $f$ restricted to $F_q$. {\color{blue} (Section \ref{s:Notation})}
    \item $|\boldsymbol{u}| = u_1 + \dots + u_c$ where $\boldsymbol{u}=(u_1, \ldots, u_c)$. {\color{blue} (Section \ref{s:Notation})}
    \item $\boldsymbol{e}_j$ the $j$th standard basis vector. {\color{blue} (Section \ref{s:Notation})}
    \item $d_{\alpha,\beta}$: the coefficient of $\beta$ in $d_q (\alpha)$ where $\alpha \in \B_q$ and $\beta \in \B_{q-1}$, the $(\beta,\alpha)$ entry of matrix $d_q$. {\color{blue} (Definition \ref{def:based})}
    \item $f_{\alpha,\beta}$: the coefficient of $\beta$ in $f|_q(\alpha)$ where $\alpha , \beta \in \B$. {\color{blue} (Definition \ref{def:based})}
    \item $G(\F)$: weighted directed graph vertices $\B$, edges $\alpha \to \beta$ of weight $d_{\alpha,\beta} \neq 0$. {\color{blue} (Definition \ref{def:matching})}
    \item $\M$: partial matching of unit edges of $G(\F)$. {\color{blue} (Definition \ref{def:matching})}
    \item $\Grep(\F, \M)$: $G(\F)$ with each $(\alpha, \beta) \in \M$ reversed, of weight $-(d_{\alpha,\beta})^{-1}$. {\color{blue} (Definition \ref{def:matching})}
    \item $\mu$: the entries $d_{\alpha,\beta}$ with $(\alpha, \beta) \in \M$; $\delta = d - \mu$. {\color{blue} (Section \ref{sec:amt}, just before Definition \ref{def:contraction})}
    \item $(i, p, h)$: contraction, $pi = \id$, $ip - \id = dh + hd$, $h^2 = hi = ph = 0$. {\color{blue} (Definition \ref{def:contraction})}
    \item $\Phi = \sum_{n \ge 0} (h\delta)^n = (\id - h\delta)^{-1}$. {\color{blue} (Section \ref{sec:amt}, just after Definition \ref{def:contraction})}
    \item $\Gamma(\alpha, \beta) = \sum_\pi \wt(\pi)$ over gradient paths $\pi$ from $\alpha$ to $\beta$. {\color{blue} (Section \ref{sec:amt}, just before Theorem \ref{thm:amt})}
    \item $d^{\M} = p \delta \Phi i$, $d^{\M}_{\alpha,\beta} = \Gamma(\alpha, \beta)$; $\F^{\M} = (K_\bullet, d^{\M})$. {\color{blue} (Theorem \ref{thm:amt})}
    \item $(i^{\M}, p^{\M}, h^
{\M}) = (\Phi i, \, p + p \delta \Phi h, \, \Phi h)$. {\color{blue} (Theorem \ref{thm:amt})}
    \item $t_j$: CI operator of degree $-2$, $t_j(y^{\boldsymbol{v}}) = y^{\boldsymbol{v} - \boldsymbol{e}_j}$; $t^{\boldsymbol{u}} = t_1^{u_1} \cdots t_c^{u_c}$. {\color{blue} (Definition \ref{def:dp})}
    \item $R\{y_j\}_{j=1}^{c} = \bigoplus_{\boldsymbol{u}} R y^{\boldsymbol{u}}$, $y^{\boldsymbol{u}} = y_1^{(u_1)} \cdots y_c^{(u_c)}$ of degree $2|\boldsymbol{u}|$. {\color{blue} (Definition \ref{def:dp})}
    \item $\sigma = \{\sigma_{\boldsymbol{u}}\}$: system of higher homotopies, $\sigma_{\boldsymbol{0}} = d$; $\sigma^{(j)} = \sigma_{\boldsymbol{e}_j}$. {\color{blue} (Definition \ref{def:higher})}
    \item $\Sh(\F, \sigma) = R\{y_j\}_{j=1}^{c} \otimes_R \F \otimes_R \Rbar$, $d_{\Sh} = \sum_{\boldsymbol{u}} t^{\boldsymbol{u}} \otimes \sigma_{\boldsymbol{u}} \otimes \id$ $= \overline{D}$, where $D = \sum_{\boldsymbol{u}} t^{\boldsymbol{u}} \otimes \sigma_{\boldsymbol{u}}$ is a map of $R\{y_j\}_{j=1}^{c} \otimes_R \F$. {\color{blue} (Definition \ref{def:shamash})}
    \item $\overline{F}_\bullet = \F \otimes_R \Rbar$ with basis $\bar\alpha = \alpha \otimes \overline{1}$, and $\bar g = g \otimes \id$ for an $R$-linear map $g$: a bar is reduction modulo $\aaa$. {\color{blue} (Definition \ref{def:shamash})}
    \item $y^{\boldsymbol{u}} \alpha$: the cell $y^{\boldsymbol{u}} \otimes \bar\alpha$ of $\Sh(\F, \sigma)$. In matrices, examples and figures we omit the tensor sign and the bar, write $\alpha$ for $1 \otimes \bar\alpha$, and write $y^{\boldsymbol{u}} F_\bullet$ for the copy $y^{\boldsymbol{u}} \otimes \overline{F}_\bullet$. {\color{blue} (Section \ref{sec:decomp}, Section \ref{sec:hyp})}
    \item $\partial_0 = \id \otimes d$, $\partial_j = t_j \otimes \sigma^{(j)}$, $\dcon = \partial_1 + \dots + \partial_c$, maps of $R\{y_j\}_{j=1}^{c} \otimes_R \F$, with $D = \partial_0 + \dcon$ and $d_{\Sh} = \bar\partial_0 + \bar\partial_{\mathrm{con}}$. {\color{blue} (Remark \ref{rem:piecessquare})}
    \item $\Theta = \sum_j t_j \otimes a_j \id$, equal to $D^2$. {\color{blue} (Remark \ref{rem:piecessquare})}
    \item $N_k = \binom{k+c-1}{c-1}$, the number of $\boldsymbol{u}$ with $|\boldsymbol{u}| = k$. {\color{blue} (Section \ref{sec:block})}
    \item $T^{(k)}_{\boldsymbol{v}}$: matrix of $t^{\boldsymbol{v}}$ from level $k$ to level $k - |\boldsymbol{v}|$; $S_k = \sum_j T^{(k)}_{\boldsymbol{e}_j} \otimes \sigma^{(j)}$. {\color{blue} (Section \ref{sec:block})}
    \item $\eps_S$: Taylor cell, $S \subseteq \{1, \dots, r\}$, $\hdeg \eps_S = |S|$, $m_S = \lcm(m_i : i \in S) = \gr(\eps_S)$. {\color{blue} (Section \ref{sec:sobieska})}
    
    \item $\lift{M} = \{(y^{\boldsymbol{u}} \otimes \bar\alpha, y^{\boldsymbol{u}} \otimes \bar\beta) : (\alpha, \beta) \in \M\}$. {\color{blue} (Definition \ref{def:lift})}
    \item $\sigma^{\M}_{\boldsymbol{0}} = d^{\M}$, $\sigma^{\M}_{\boldsymbol{u}} = \sum p^{\M} \sigma^{(j_1)} h^{\M} \cdots h^{\M} \sigma^{(j_k)} i^{\M}$ over $\boldsymbol{e}_{j_1} + \dots + \boldsymbol{e}_{j_k} = \boldsymbol{u}$. {\color{blue} (Definition \ref{def:transferred})}
    \item $\lift{f} = \id \otimes f$, a map of $R\{y_j\}_{j=1}^{c} \otimes_R \F$, so $\lift{f}(y^{\boldsymbol{u}} \otimes \alpha) = y^{\boldsymbol{u}} \otimes f(\alpha)$. On $\Sh(\F, \sigma)$ we use $\overline{\lift{f}}$. For the Greek maps $\delta$ and $\mu$ we write $\id \otimes \delta$ and $\id \otimes \mu$. {\color{blue} (Section \ref{sec:commute}, just before Proposition \ref{prop:closed})}
\end{itemize}

\section{Algebraic Morse theory}\label{sec:amt}

This section fixes the form of algebraic Morse theory we use. For monomial ideals, Batzies and Welker \cite{BatziesWelker} apply Forman's discrete Morse theory \cite{Forman} to a CW complex supporting a cellular resolution. An acyclic matching of cells with equal monomial labels collapses the complex to a smaller CW complex, which again supports a cellular resolution, so the reduction respects cellularity. Algebraic Morse theory \cite{JollenbeckWelker,Skoldberg} works with the free complex directly and needs no supporting cell complex. For monomial ideals there are algorithms that produce Morse matchings on the Taylor resolution: the pruning algorithm of \`Alvarez Montaner, Fern\'andez-Ramos and Gimenez \cite[Algorithm~3.1]{AlvarezFernandezGimenez}, and the Barile--Macchia algorithm, introduced in \cite{BarileMacchia} for edge ideals of forests and extended to all monomial ideals by Chau and Kara \cite[Algorithm~2.9]{ChauKara}. The results below apply to any of these matchings. Everything through Theorem~\ref{thm:amt} holds over an arbitrary commutative ring $A$, and we apply it over $R$ and over $\Rbar$.

\begin{definition}\label{def:based}
    Let $A$ be a commutative ring. A \emph{free complex} over $A$ is a chain complex $(\F, d)$ of free $A$-modules together with a chosen basis $\B_q$ of each $F_q$. A basis element $\alpha \in \B_q$ is a \emph{$q$-cell}, and $\B = \bigsqcup_q \B_q$.

    The \emph{entries} $d_{\alpha,\beta} \in A$ of $d$ are defined by $d_q(\alpha) = \sum_{\beta \in \B_{q-1}} d_{\alpha,\beta} \beta$ for $\alpha \in \B_q$, a finite sum, and any
    collection of $A$-linear maps $f \colon F_\bullet \to F_{\bullet+r}$ has entries $f_{\alpha,\beta}$ defined the same way. The \emph{complex graph} $G(\F)$ has vertex set $\B$ and an edge $\alpha \to \beta$ of \emph{weight} $d_{\alpha,\beta}$ whenever $d_{\alpha,\beta} \neq 0$. An edge is a \emph{unit edge} if its weight is a unit of $A$.
\end{definition}

\begin{remark}
     As seen from the definition, the basis is an important part of the data.
    Sk\"oldberg calls the same object a \emph{based complex} \cite[\S2]{Skoldberg}.
\end{remark}

\begin{definition}[\cite{Taylor}]
    For a monomial ideal $I = (m_1, \dots, m_r)$ the Taylor resolution is the free complex $\F$ whose cells $\eps_S$ are indexed by the subsets $S \subseteq \{1, \dots, r\}$, where $\hdeg(\eps_S) = |S|$ and $\gr(\eps_S) = m_S = \lcm(m_i : i \in S)$, and whose differential is
    \[
        d(\eps_S) = \sum_{i \in S} (-1)^{\#\{s \in S : s > i\}} \frac{m_S}{m_{S \setminus i}} \eps_{S \setminus i} .
    \]
    The sign ensures that the first nonzero entry of every column in every matrix of $d$ is positive. It counts the elements of $S$ above the one being deleted, so the term that deletes the largest index carries a plus. This is the sign convention of \cite[Def.~3.1]{Sobieska}, where the sign is written $(-1)^{k-p}$ at the $p$th element of $S = \{s_1 < \dots < s_k\}$, with $k = |S|$, counting the same elements.
\end{definition}

\begin{example}[The running example]\label{ex:runtaylor} Take $R = \kk[x, y, z, w]$, the ideal $I = (xy, yz, zw, xw)$ with generators $m_1, m_2, m_3, m_4$ in that order, and $\F$ the Taylor resolution of $R/I$, with cells $\eps_S$ in homological degree $|S|$, ordered lexicographically.  The differentials of $\F$ are shown in Figure~\ref{fig:runTaylor} beside $G(\F)$, each at the height of the edges it describes.

\begin{figure}[H]
    \centering
    \resizebox{0.68\textwidth}{!}{%
  \begin{tikzpicture}[>=Stealth, x=1cm, y=1cm,
        vtx/.style={inner sep=1.2pt, font=\scriptsize, draw, rounded corners=1pt, fill=white, minimum width=5mm, minimum height=3.5mm},
        uedge/.style={->, black, thick},
        dedge/.style={->, gray!60},
        wt/.style={font=\tiny, fill=white, inner sep=0.4pt, text=black},
        mat/.style={font=\scriptsize, anchor=west}]
        \renewcommand{\arraystretch}{0.75}
        \node[vtx, align=center] (e0) at (0.00,0.00) {$\eps_\varnothing$\\[-1.5pt]{\tiny $1$}};
        \node[vtx, align=center] (e1) at (-2.17,1.90) {$\eps_{1}$\\[-1.5pt]{\tiny $xy$}};
        \node[vtx, align=center] (e2) at (-0.72,1.90) {$\eps_{2}$\\[-1.5pt]{\tiny $yz$}};
        \node[vtx, align=center] (e3) at (0.72,1.90) {$\eps_{3}$\\[-1.5pt]{\tiny $zw$}};
        \node[vtx, align=center] (e4) at (2.17,1.90) {$\eps_{4}$\\[-1.5pt]{\tiny $xw$}};
        \node[vtx, align=center] (e12) at (-3.62,3.80) {$\eps_{12}$\\[-1.5pt]{\tiny $xyz$}};
        \node[vtx, align=center] (e13) at (-2.17,3.80) {$\eps_{13}$\\[-1.5pt]{\tiny $xyzw$}};
        \node[vtx, align=center] (e14) at (-0.72,3.80) {$\eps_{14}$\\[-1.5pt]{\tiny $xyw$}};
        \node[vtx, align=center] (e23) at (0.72,3.80) {$\eps_{23}$\\[-1.5pt]{\tiny $yzw$}};
        \node[vtx, align=center] (e24) at (2.17,3.80) {$\eps_{24}$\\[-1.5pt]{\tiny $xyzw$}};
        \node[vtx, align=center] (e34) at (3.62,3.80) {$\eps_{34}$\\[-1.5pt]{\tiny $xzw$}};
        \node[vtx, align=center] (e123) at (-2.17,5.70) {$\eps_{123}$\\[-1.5pt]{\tiny $xyzw$}};
        \node[vtx, align=center] (e124) at (-0.72,5.70) {$\eps_{124}$\\[-1.5pt]{\tiny $xyzw$}};
        \node[vtx, align=center] (e134) at (0.72,5.70) {$\eps_{134}$\\[-1.5pt]{\tiny $xyzw$}};
        \node[vtx, align=center] (e234) at (2.17,5.70) {$\eps_{234}$\\[-1.5pt]{\tiny $xyzw$}};
        \node[vtx, align=center] (e1234) at (0.00,7.60) {$\eps_{1234}$\\[-1.5pt]{\tiny $xyzw$}};
        \draw[dedge] (e1) -- node[wt, pos=0.46] {$xy$} (e0);
        \draw[dedge] (e2) -- node[wt, pos=0.62] {$yz$} (e0);
        \draw[dedge] (e3) -- node[wt, pos=0.38] {$zw$} (e0);
        \draw[dedge] (e4) -- node[wt, pos=0.38] {$xw$} (e0);
        \draw[dedge] (e12) -- node[wt, pos=0.46] {$-x$} (e2);
        \draw[dedge] (e12) -- node[wt, pos=0.38] {$z$} (e1);
        \draw[dedge] (e13) -- node[wt, pos=0.46] {$-xy$} (e3);
        \draw[dedge] (e13) -- node[wt, pos=0.7] {$zw$} (e1);
        \draw[dedge] (e14) -- node[wt, pos=0.46] {$-y$} (e4);
        \draw[dedge] (e14) -- node[wt, pos=0.54] {$w$} (e1);
        \draw[dedge] (e23) -- node[wt, pos=0.7] {$-y$} (e3);
        \draw[dedge] (e23) -- node[wt, pos=0.54] {$w$} (e2);
        \draw[dedge] (e24) -- node[wt, pos=0.46] {$-yz$} (e4);
        \draw[dedge] (e24) -- node[wt, pos=0.3] {$xw$} (e2);
        \draw[dedge] (e34) -- node[wt, pos=0.38] {$-z$} (e4);
        \draw[dedge] (e34) -- node[wt, pos=0.7] {$x$} (e3);
        \draw[dedge] (e123) -- node[wt, pos=0.3] {$x$} (e23);
        \draw[uedge] (e123) -- node[wt, pos=0.46] {$-1$} (e13);
        \draw[dedge] (e123) -- node[wt, pos=0.62] {$w$} (e12);
        \draw[uedge] (e124) -- node[wt, pos=0.46] {$1$} (e24);
        \draw[dedge] (e124) -- node[wt, pos=0.7] {$-z$} (e14);
        \draw[dedge] (e124) -- node[wt, pos=0.62] {$w$} (e12);
        \draw[dedge] (e134) -- node[wt, pos=0.7] {$y$} (e34);
        \draw[dedge] (e134) -- node[wt, pos=0.3] {$-z$} (e14);
        \draw[uedge] (e134) -- node[wt, pos=0.46] {$1$} (e13);
        \draw[dedge] (e234) -- node[wt, pos=0.46] {$y$} (e34);
        \draw[uedge] (e234) -- node[wt, pos=0.46] {$-1$} (e24);
        \draw[dedge] (e234) -- node[wt, pos=0.54] {$x$} (e23);
        \draw[uedge] (e1234) -- node[wt, pos=0.46] {$-1$} (e234);
        \draw[uedge] (e1234) -- node[wt, pos=0.38] {$1$} (e134);
        \draw[uedge] (e1234) -- node[wt, pos=0.62] {$-1$} (e124);
        \draw[uedge] (e1234) -- node[wt, pos=0.38] {$1$} (e123);   
        \node[mat] at (4.5,0.95) {$d_{1} = \begin{bmatrix} xy & yz & zw & xw \end{bmatrix}$};
        \node[mat] at (4.5,2.85) {$d_{2} = \begin{bmatrix} z & zw & w & 0 & 0 & 0 \\ -x & 0 & 0 & w & xw & 0 \\ 0 & -xy & 0 & -y & 0 & x \\ 0 & 0 & -y & 0 & -yz & -z \end{bmatrix}$};
        \node[mat] at (4.5,4.75) {$d_{3} = \begin{bmatrix} w & w & 0 & 0 \\ -1 & 0 & 1 & 0 \\ 0 & -z & -z & 0 \\ x & 0 & 0 & x \\ 0 & 1 & 0 & -1 \\ 0 & 0 & y & y \end{bmatrix}$};
        \node[mat] at (4.5,6.65) {$d_{4} = \begin{bmatrix} 1 \\ -1 \\ 1 \\ -1 \end{bmatrix}$};
    \end{tikzpicture}}
    \caption{The running example: the Taylor resolution and its differentials}
    \label{fig:runTaylor}
\end{figure}
\end{example}

\begin{definition}\label{def:matching}
    A \emph{partial matching} $\M$ on a free complex $\F$ over $A$ is a set of unit edges of $G(\F)$, no two of which share a vertex. We use the following words for $\M$.

    \begin{enumerate}[label=(\arabic*),leftmargin=*,itemsep=1pt]
        \item For $(\alpha, \beta) \in \M$ we call $\alpha$ the \emph{upper} and $\beta$ the \emph{lower} cell of the pair, the homological degree of $\beta$ being one less than that of $\alpha$. A cell lying on no edge of $\M$ is \emph{critical}.

        \item The \emph{replacement graph} $\Grep(\F, \M)$ has the vertices of $G(\F)$, the edges of $G(\F)$ that lie outside $\M$ with their weights unchanged, and for each $(\alpha, \beta) \in \M$ a reversed edge $\beta \to \alpha$ of weight $-(d_{\alpha,\beta})^{-1}$.

        \item A \emph{directed path} in a directed graph is a sequence of edges, finite or infinite, each starting at the vertex where the previous one ends, with repeated vertices allowed.

        \item The matching is \emph{acyclic} if $\Grep(\F, \M)$ contains no directed cycle, and \emph{locally finite} if every directed path in $\Grep(\F, \M)$ is finite.

        \item A \emph{Morse matching} is a locally finite partial matching.
    \end{enumerate}
\end{definition}

\begin{remark}\label{rem:acyclic}
    Algebraic Morse theory is usually stated for \emph{acyclic} matchings, those for which $\Grep(\F, \M)$ has no directed cycle \cite{Skoldberg,JollenbeckWelker}. A locally finite matching is acyclic, since a directed cycle traversed over and over is an infinite directed path. When $\B$ is finite the converse holds, since an infinite directed path in a finite graph contains a directed cycle. We state it for locally finite matchings instead, as this makes the perturbation lemma (Lemma~\ref{lem:perturbation}) work without additional structure. 
\end{remark}

Given $\M$, split $d = \mu + \delta$, where the \emph{matched part} $\mu$ keeps the entries $d_{\alpha,\beta}$ with $(\alpha, \beta) \in \M$ and the \emph{unmatched part} is $\delta = d - \mu$.

\begin{remark} \label{rem:blockmu}
    For each homological degree $q$, we can partition the basis $\B_q$ into three sets $K_q$ - the critical cells, $U_q$ - the upper cells, and $L_q$ - the lower cells. We can rewrite $d_q$ as a block matrix with respect to the reordered bases $\{K_q, U_q,L_q\}$ and $\{K_{q-1}, U_{q-1},L_{q-1}\}$ as illustrated below.
    $$\begin{array}{cc}
     & \begin{array}{ccc}
     K_q \ \  & U_q \ \ & L_q\\
     \end{array}\\
    \begin{array}{c}
         K_{q-1} \\
         U_{q-1}\\
         L_{q-1}\\
    \end{array}&\left[\begin{array}{c|c|c}
     d_{KK}&d_{UK}  & d_{LK} \\
    \hline
    d_{KU}  & d_{UU}  & d_{LU} \\
    \hline
    d_{KL} & d_{UL} & d_{LL} \\\end{array}\right]\end{array},$$
where $d_{UL}=\delta_{UL}+\mu_{UL}$ with $\delta_{UL}$ representing the edges between upper and lower cells that are not matched edges, and $\mu_{UL}$ being all those edges that are in the matching.
In terms of this ordering of the bases, $$\delta= \left[\begin{array}{c|c|c}
     d_{KK}&d_{UK}  & d_{LK} \\
    \hline
    d_{KU}  & d_{UU}  & d_{LU} \\
    \hline
    d_{KL} & \delta_{UL} & d_{LL} \\\end{array}\right] \text{ and } \mu=\left[\begin{array}{c|c|c}
     0& 0  & 0 \\
    \hline
    0  & 0 & 0 \\
    \hline
    0 & \mu_{UL} & 0 \\\end{array}\right],$$

    \noindent where $\mu_{UL}$ is an invertible matrix by construction.

\end{remark}

\begin{definition}\label{def:contraction}
    A \emph{contraction} $(i, p, h)$ of $(\F, d)$ onto $(\F', d')$ consists of chain maps $i \colon \F' \to \F$ and $p \colon \F \to \F'$ of homological degree $0$ and a collection of $A$-linear maps $h \colon F_\bullet \to F_{\bullet+1}$ of homological degree 1 with $pi = \id$, $ip - \id = dh + hd$ and $h^2 = hi = ph = 0$. In particular $i$ and $p$ are mutually inverse homotopy equivalences.
\end{definition}

A \emph{perturbation} of a complex $(\F, d)$ is an $A$-linear map $\delta \colon F_\bullet \to F_{\bullet-1}$ with $(d + \delta)^2 = 0$. For $A$-linear maps $h$ and $\delta$ of a module to itself, we say $\delta$ is \emph{small} with respect to $h$ if $h\delta$ is locally nilpotent, that is, for every element $x$ we have $(h\delta)^n x = 0$ for some $n$. Then
\[
    \Phi = \sum_{n \ge 0} (h\delta)^n
\]
is a finite sum on each element, and $\Phi = (\id - h\delta)^{-1}$.

We state the perturbation lemma as Sk\"oldberg does in \cite[Theorem~1]{SkoldbergHPT}, with one change. He assumes a filtration, which he uses only to show that $\delta$ is small with respect to $h$. Since smallness is immediately implied by local finiteness of the matching, we do not need to worry about the filtration.

\begin{lemma}[Perturbation lemma, cf.~{\cite[Theorem~1]{SkoldbergHPT}}]\label{lem:perturbation}
    Let $(i, p, h)$ be a contraction of a complex $(\F, d)$ of $A$-modules onto a complex $(\F', d')$ of $A$-modules and let $\delta$ be a small perturbation of $(\F, d)$. Then
    \[
        \bigl( \Phi i, \; p + p \delta \Phi h, \; \Phi h \bigr)
    \]
    is a contraction of $(\F, d + \delta)$ onto $(\F', d' + p \delta \Phi i)$.
\end{lemma}

The main theorem of algebraic Morse theory is due independently to Sk\"oldberg \cite{Skoldberg} and to J\"ollenbeck and Welker \cite[Theorem~2.2]{JollenbeckWelker}. Sk\"oldberg later showed that algebraic Morse theory follows from the perturbation lemma \cite[Theorem~2]{SkoldbergHPT}, and the formulas for the contraction below are the ones his proof produces.

A \emph{gradient path} $\pi$ from $\alpha$ to $\beta$ is a directed path $\alpha = v_0 \to v_1 \to \dots \to v_{2m+1} = \beta$ in $\Grep(\F, \M)$, $m \ge 0$, in which $v_{2i} \to v_{2i+1}$ is an edge of $G(\F)$ outside $\M$ and $v_{2i+1} \to v_{2i+2}$ is a reversed matched edge. Its weight $\wt(\pi)$ is the product of the weights of its edges. A gradient path lowers homological degree by one, and for cells $\alpha$ and $\beta$ with $\hdeg(\beta) = \hdeg(\alpha) - 1$ we write $\Gamma(\alpha, \beta) = \sum_\pi \wt(\pi)$, the sum over all gradient paths from $\alpha$ to $\beta$. For a Morse matching it has only finitely many terms, because every cell has finitely many outgoing edges and every directed path is finite.

\begin{theorem}[Algebraic Morse theory, {\cite{Skoldberg,JollenbeckWelker}}]\label{thm:amt}
   Let $\M$ be a Morse matching on a free complex $(\F, d)$ over a commutative ring $A$, with $d = \mu + \delta$. Let $K_\bullet$ be the graded free $A$-module on the critical cells, $i \colon K_\bullet \to \F$ the inclusion of the critical cells and $p \colon \F \to K_\bullet$ the projection killing the noncritical cells. The map $\mu$ restricts to an isomorphism from the span of the upper cells onto the span of the lower cells, and we let $h$ be minus its inverse on the lower cells and zero on the critical and upper cells.

   Then $(i, p, h)$ is a contraction of $(\F, \mu)$ onto $(K_\bullet, 0)$, the sum $\Phi = \sum_{n \ge 0} (h\delta)^n$ is finite on each element, and perturbing by $\delta$ gives a contraction of $(\F, d)$ onto the \emph{Morse complex} ${\F^{\M} = (K_\bullet, d^{\M})}$,
    \[
        \begin{tikzpicture}[>=Stealth, x=1cm, y=1cm, font=\small, baseline=(current bounding box.center)]
            \node (L1) at (0, 1.9) {$(\F, \mu)$};
            \node (R1) at (5.2, 1.9) {$(K_\bullet, 0)$};
            \node (L2) at (0, 0) {$(\F, d)$};
            \node (R2) at (5.2, 0) {$\F^{\M} = (K_\bullet, d^{\M})$};
            \foreach \r/\p/\i in {1/{p}/{i}, 2/{p^{\M}}/{i^{\M}}}{
            \draw[->] ([yshift=3.5pt]L\r.east) -- node[above, font=\footnotesize]{$\p$} ([yshift=3.5pt]R\r.west);
            \draw[<-] ([yshift=-3.5pt]L\r.east) -- node[below, font=\footnotesize]{$\i$} ([yshift=-3.5pt]R\r.west);}
            \foreach \r/\h in {1/{h}, 2/{h^{\M}}}{
            \draw[->] ([yshift=3pt]L\r.west) .. controls +(-0.9, 0.4) and +(-0.9, -0.4) .. node[left, font=\footnotesize]{$\h$} ([yshift=-3pt]L\r.west);}
            \draw[->, thick] (L1) -- node[right, font=\footnotesize]{$+\delta$} (L2);
        \end{tikzpicture}
    \]
    where
    \[
        d^{\M} = p \delta \Phi i, \qquad d^{\M}_{\alpha,\beta} = \Gamma(\alpha, \beta) \quad \text{for critical } \alpha, \beta, \qquad (i^{\M}, p^{\M}, h^{\M}) = \bigl( \Phi i, \; p + p \delta \Phi h, \; \Phi h \bigr) .
    \]
\end{theorem}

Theorem~\ref{thm:amt} is Lemma~\ref{lem:perturbation} applied to $(i, p, h)$ and the perturbation $\delta$, since $(\mu + \delta)^2 = d^2 = 0$, with local finiteness making $\Phi$ a finite sum on each cell. The contraction $(i, p, h)$ uses only that the matched entries are invertible. In the block form of Remark~\ref{rem:blockmu}, with columns $K_q, U_q, L_q$, rows $K_{q+1}, U_{q+1}, L_{q+1}$, and $\mu_{UL}$ the matched block of $d_{q+1}$, $$h|_q=\left[\begingroup
        \renewcommand{\arraystretch}{1.5}\begin{array}{c|c|c}
     0& 0  & 0 \\
    \hline
    0  & 0 & -\mu_{UL}^{-1} \\
    \hline
    0 & 0 & 0 \\\end{array}\endgroup\right].$$

Hence, when every matched weight is $\pm 1$, as for the Taylor resolution, $h|_q = -\mu_{q+1}^{T}$.

In the graph, $\delta$ follows an edge outside $\M$, and $h$ carries a lower cell up its reversed matched edge, with weight $-(d_{\alpha,\beta})^{-1}$, and kills every other cell. So $d^{\M} = p \delta \Phi i$ sums over gradient paths between critical cells. Likewise $i^{\M}$ and $p^{\M}$ sum over paths that alternate edges outside $\M$ and reversed matched edges, starting with an edge outside $\M$ for $i^{\M}$ and with a reversed edge for $p^{\M}$. So $i^{\M}(\kappa)$ involves only $\kappa$ and upper cells, and $p^{\M}$ agrees with $p$ except on the lower cells.

\begin{example}[The running example]\label{ex:runSec3}  We continue Example \ref{ex:runtaylor} with $R = \kk[x, y, z, w]$, $I = (m_1, m_2, m_3, m_4) = (xy, yz, zw, xw)$ and the matching: \[
    \M = \bigl\{ (\eps_{124}, \eps_{24}), (\eps_{123}, \eps_{13}), (\eps_{1234}, \eps_{234}) \bigr\}.
\]
In Figure~\ref{fig:runex} we write $\Grep(\F, \M)$ along with the maps $\mu$ and $\delta$. $\M$ produces no directed cycle, so $\M$ is a Morse matching. It is the matching produced by the pruning algorithm \cite[Algorithm~3.1]{AlvarezFernandezGimenez}.

\begin{figure}[htbp]
    \centering

    \begin{tikzpicture}[>=Stealth, x=1cm, y=1cm,
        vtx/.style={inner sep=1.2pt, font=\scriptsize, draw, rounded corners=1pt, fill=white, minimum width=5mm, minimum height=3.5mm},
        crit/.style={vtx, fill=gray!10},
        mu/.style={vtx, fill=red!12},
        md/.style={vtx, fill=orange!18},
        uedge/.style={->, black, thick},
        medge/.style={->, red!70!black, thick},
        dedge/.style={->, gray!60},
        wt/.style={font=\tiny, fill=white, inner sep=0.4pt, text=black},
        mat/.style={font=\scriptsize, anchor=west}]
        \node[crit] (e0) at (0.00,0.00) {$\eps_\varnothing$};
        \node[crit] (e1) at (-2.17,1.90) {$\eps_{1}$};
        \node[crit] (e2) at (-0.72,1.90) {$\eps_{2}$};
        \node[crit] (e3) at (0.72,1.90) {$\eps_{3}$};
        \node[crit] (e4) at (2.17,1.90) {$\eps_{4}$};
        \node[crit] (e12) at (-3.62,3.80) {$\eps_{12}$};
        \node[md] (e13) at (-2.17,3.80) {$\eps_{13}$};
        \node[crit] (e14) at (-0.72,3.80) {$\eps_{14}$};
        \node[crit] (e23) at (0.72,3.80) {$\eps_{23}$};
        \node[md] (e24) at (2.17,3.80) {$\eps_{24}$};
        \node[crit] (e34) at (3.62,3.80) {$\eps_{34}$};
        \node[mu] (e123) at (-2.17,5.70) {$\eps_{123}$};
        \node[mu] (e124) at (-0.72,5.70) {$\eps_{124}$};
        \node[crit] (e134) at (0.72,5.70) {$\eps_{134}$};
        \node[md] (e234) at (2.17,5.70) {$\eps_{234}$};
        \node[mu] (e1234) at (0.00,7.60) {$\eps_{1234}$};
        \draw[dedge] (e1) -- (e0);
        \draw[dedge] (e2) -- (e0);
        \draw[dedge] (e3) -- (e0);
        \draw[dedge] (e4) -- (e0);
        \draw[dedge] (e12) -- (e2);
        \draw[dedge] (e12) -- (e1);
        \draw[dedge] (e13) -- (e3);
        \draw[dedge] (e13) -- (e1);
        \draw[dedge] (e14) -- (e4);
        \draw[dedge] (e14) -- (e1);
        \draw[dedge] (e23) -- (e3);
        \draw[dedge] (e23) -- (e2);
        \draw[dedge] (e24) -- (e4);
        \draw[dedge] (e24) -- (e2);
        \draw[dedge] (e34) -- (e4);
        \draw[dedge] (e34) -- (e3);
        \draw[dedge] (e123) -- (e23);
        \draw[medge] (e13) -- node[wt, pos=0.35, text=red!70!black] {$1$} (e123);
        \draw[dedge] (e123) -- (e12);
        \draw[medge] (e24) -- node[wt, pos=0.35, text=red!70!black] {$-1$} (e124);
        \draw[dedge] (e124) -- (e14);
        \draw[dedge] (e124) -- (e12);
        \draw[dedge] (e134) -- (e34);
        \draw[dedge] (e134) -- (e14);
        \draw[uedge] (e134) -- node[wt, pos=0.3] {$1$} (e13);
        \draw[dedge] (e234) -- (e34);
        \draw[uedge] (e234) -- node[wt, pos=0.3] {$-1$} (e24);
        \draw[dedge] (e234) -- (e23);
        \draw[medge] (e234) -- node[wt, pos=0.35, text=red!70!black] {$1$} (e1234);
        \draw[uedge] (e1234) -- node[wt, pos=0.3] {$1$} (e134);
        \draw[uedge] (e1234) -- node[wt, pos=0.3] {$-1$} (e124);
        \draw[uedge] (e1234) -- node[wt, pos=0.3] {$1$} (e123);
        \node[mat] at (4.75,0.95) {$\mu_1 = 0, \quad \delta_1 = d_1$};
        \node[mat] at (4.75,2.70) {$\mu_2 = 0, \quad \delta_2 = d_2$};
        \node[mat] at (4.75,4.84) {$\mu_3 = \begin{bmatrix} 0&0&0&0\\ -1&0&0&0\\ 0&0&0&0\\ 0&0&0&0\\ 0&1&0&0\\ 0&0&0&0 \end{bmatrix}, \quad \delta_3 = \begin{bmatrix} w&w&0&0\\ 0&0&1&0\\ 0&-z&-z&0\\ x&0&0&x\\ 0&0&0&-1\\ 0&0&y&y \end{bmatrix}$};
        \node[mat] at (4.75,6.88) {$\mu_4 = \begin{bmatrix} 0\\ 0\\ 0\\ -1 \end{bmatrix}, \quad \delta_4 = \begin{bmatrix} 1\\ -1\\ 1\\ 0 \end{bmatrix}$};
        \path (-4.03,0) (12.65,0);
    \end{tikzpicture}

    \caption{The running example. Red, orange, gray: upper, lower, critical cells.}
    \label{fig:runex}
\end{figure}

No entry of $\delta$ joins an upper cell to a lower cell, so $h^{\M} = h = -\mu^{T}$, which is zero except at the lower cells,
\[
    h(\eps_{13}) = \eps_{123}, \qquad h(\eps_{24}) = -\eps_{124}, \qquad h(\eps_{234}) = \eps_{1234} .
\]
The map $p^{\M}$ differs from $p$ only on the lower cells, for example $p^{\M}(\eps_{13}) = w \eps_{12} + x \eps_{23}$ from the paths ${\eps_{13} \xrightarrow{1} \eps_{123} \xrightarrow{w} \eps_{12}}$ and ${\eps_{13} \xrightarrow{1} \eps_{123} \xrightarrow{x} \eps_{23}}$. In degrees $0$ and $1$ every cell is critical, $i$, $p$, $i^{\M}$, $p^{\M}$ are identities and $h = h^{\M} = 0$. In degrees $2$ and $3$, with the critical cells $\eps_{12}, \eps_{14}, \eps_{23}, \eps_{34}$ and $\eps_{134}$, the contraction $(i, p, h)$ is
\[
    i|_2 = \begin{bmatrix} 1 & 0 & 0 & 0\\ 0 & 0 & 0 & 0\\ 0 & 1 & 0 & 0\\ 0 & 0 & 1 & 0\\ 0 & 0 & 0 & 0\\ 0 & 0 & 0 & 1 \end{bmatrix}, \quad   
    i|_3 = \begin{bmatrix} 0\\ 0\\ 1\\ 0 \end{bmatrix}, \quad p|_2 = \begin{bmatrix} 1 & 0 & 0 & 0 & 0 & 0\\ 0 & 0 & 1 & 0 & 0 & 0\\ 0 & 0 & 0 & 1 & 0 & 0\\ 0 & 0 & 0 & 0 & 0 & 1 \end{bmatrix}, \quad
    p|_3 = \begin{bmatrix} 0 & 0 & 1 & 0 \end{bmatrix},
\]
\[
    h|_2 = \begin{bmatrix} 0 & 1 & 0 & 0 & 0 & 0\\ 0 & 0 & 0 & 0 & -1 & 0\\ 0 & 0 & 0 & 0 & 0 & 0\\ 0 & 0 & 0 & 0 & 0 & 0 \end{bmatrix}, \quad
    h|_3 = \begin{bmatrix} 0 & 0 & 0 & 1 \end{bmatrix},
\]
and the contraction $(i^{\M}, p^{\M}, h^{\M})$ is
\[
    i^{\M}|_2 = i|_2, \quad
    i^{\M}|_3 = \begin{bmatrix} 1\\ 0\\ 1\\ 0 \end{bmatrix}, \quad
    p^{\M}|_2 = \begin{bmatrix} 1 & w & 0 & 0 & -w & 0\\ 0 & 0 & 1 & 0 & z & 0\\ 0 & x & 0 & 1 & 0 & 0\\ 0 & 0 & 0 & 0 & 0 & 1 \end{bmatrix}, \quad
    p^{\M}|_3 = \begin{bmatrix} 0 & 0 & 1 & 1 \end{bmatrix}, \quad
    h^{\M} = h .
\]
The Morse complex $\F^{\M}$ is the minimal free resolution of $R/I$ with differentials
\[
    d^{\M}_1 = d_1, \qquad
    d^{\M}_2 = \begin{bmatrix} z & w & 0 & 0\\ -x & 0 & w & 0\\ 0 & 0 & -y & x\\ 0 & -y & 0 & -z \end{bmatrix}, \qquad
    d^{\M}_3 = \begin{bmatrix} w\\ -z\\ x\\ y \end{bmatrix}.
\]
\end{example}

\begin{remark}[Row and column operations]\label{rem:rowcol}
    Algebraic Morse theory can really be viewed as a combinatorial solution to the linear algebra problem of isolating the units in the matrices of a free resolution. The contraction $(i^{\M}, p^{\M}, h^{\M})$ of Theorem~\ref{thm:amt} gives a change of basis that splits $\F$ into $\F^{\M}$ and trivial complexes. Split $\B_q$ into the critical cells $K_q$, the upper matched cells $U_q$ and the lower matched cells $L_q$, and let $\iota_{U_q}$ be the inclusion of the span of $U_q$. The identity $i^{\M} p^{\M} - \id = d h^{\M} + h^{\M} d$ lets each $x$ in $F_q$ be expressed by
    \[
        x = i^{\M}\bigl(p^{\M}(x)\bigr) + \bigl(-h^{\M} d(x)\bigr) + d\bigl(-h^{\M}(x)\bigr),
    \]
    with $p^{\M}(x)$ in the span of $K_q$, $-h^{\M} d(x)$ in the span of $U_q$ and $-h^{\M}(x)$ in the span of $U_{q+1}$. So the new basis of $F_q$ is $i^{\M}(\kappa)$ for $\kappa \in K_q$, then $\alpha \in U_q$, then $d(\alpha')$ for $\alpha' \in U_{q+1}$, and the change of basis matrix and its inverse are
    \[
        V_q = \begin{bmatrix} i^{\M}|_q & \iota_{U_q} & d\, \iota_{U_{q+1}} \end{bmatrix}, \qquad
        V_q^{-1} = \begin{bmatrix} p^{\M}|_q\\ -h^{\M} d|_q\\ -h^{\M}|_q \end{bmatrix}.
    \]
    The columns of $i^{\M}|_q$ are the column operations on $d_q$, and the rows of $p^{\M}|_q$ are the row operations on $d_{q+1}$. In the new bases
    \[
        V_{q-1}^{-1}\, d_q\, V_q = \begin{bmatrix} p^{\M} d\, i^{\M} & 0 & 0\\ 0 & 0 & 0\\ 0 & \mathbf{1} & 0 \end{bmatrix} = \begin{bmatrix} d^{\M}_q & 0 & 0\\ 0 & 0 & 0\\ 0 & \mathbf{1} & 0 \end{bmatrix},
    \]
    so $\F \cong \F^{\M} \oplus \bigoplus_{(\alpha, \beta) \in \M} \bigl(0 \to A\alpha \to A\, d(\alpha) \to 0\bigr)$. For a single pair $(\alpha, \beta)$ of weight $w$, $i^{\M}(\gamma) = \gamma - w^{-1} d_{\gamma,\beta}\, \alpha$ subtracts $w^{-1} d_{\gamma,\beta}$ times column $\alpha$ from column $\gamma$, $p^{\M}(\beta) = -w^{-1} \sum_{\varepsilon} d_{\alpha,\varepsilon}\, \varepsilon$ subtracts $w^{-1} d_{\alpha,\varepsilon}$ times row $\beta$ from row $\varepsilon$, and $d^{\M}_{\gamma,\varepsilon} = d_{\gamma,\varepsilon} - d_{\gamma,\beta}\, w^{-1} d_{\alpha,\varepsilon}$ is the Schur complement. The spirit of J\"ollenbeck and Welker's proof of algebraic Morse theory is to apply this inductively. 
\end{remark}

\section{The Eisenbud--Shamash construction}\label{sec:shamash}
This section recalls the Eisenbud--Shamash construction. The construction, with the divided powers below, is Shamash's \cite{Shamash} for a single element $a_1$. Systems of higher homotopies for a regular sequence and the operators $t_j$ were introduced by Eisenbud in \cite{Eisenbud1980}. The operators appeared earlier on $\operatorname{Ext}$ in Gulliksen's work \cite{Gulliksen}, and are described combinatorially in Section \ref{sec:block}.

\begin{definition}\label{def:dp} \cite{Eisenbud1980}
    Let $R[t_j]_{j=1}^{c} = R[t_1, \dots, t_c]$ be the polynomial ring over $R$ with additional variables $t_j$ in degree $-2$. Let $R\{y_j\}_{j=1}^{c} = \bigoplus_{\boldsymbol{u} \in \N^c} R y^{\boldsymbol{u}}$ be its graded dual, the free $R$-module on the symbols $y^{\boldsymbol{u}} = y_1^{(u_1)} \cdots y_c^{(u_c)}$ dual to the monomials $t^{\boldsymbol{u}}$, for $\boldsymbol{u} \in \N^c$, with $y^{\boldsymbol{u}}$ in homological degree $2|\boldsymbol{u}|$.

    The dual of an $R[t_j]_{j=1}^{c}$-module is an $R[t_j]_{j=1}^{c}$-module by $(t_j \varphi)(g) = \varphi(t_j g)$, for $\varphi$ in the dual and $g$ in the module. In the dual basis $t_j(y^{\boldsymbol{v}}) = y^{\boldsymbol{v}-\boldsymbol{e}_j}$ when $v_j \ge 1$ and $t_j(y^{\boldsymbol{v}}) = 0$ when $v_j = 0$. Following \cite{EisenbudPeeva} we call $t_j$, acting this way, the \emph{CI operator} $t_j$.

    The CI operators commute. For $\boldsymbol{u} \in \N^c$ the composite $t^{\boldsymbol{u}} = t_1^{u_1} \cdots t_c^{u_c}$ satisfies $t^{\boldsymbol{u}}(y^{\boldsymbol{v}}) = y^{\boldsymbol{v}-\boldsymbol{u}}$ when $\boldsymbol{u} \le \boldsymbol{v}$ and $t^{\boldsymbol{u}}(y^{
    \boldsymbol{v}}) = 0$ otherwise. Then $t^{\boldsymbol{v}} t^{\boldsymbol{v'}} = t^{\boldsymbol{v} + \boldsymbol{v'}}$ for all $\boldsymbol{v}, \boldsymbol{v'} \in \N^c$.

    We call an index $j \in \{1, \dots, c\}$ a \emph{color}. We will later represent our base free resolution by a weighted directed graph. We will have $c$ homotopies, each corresponding to added arrows on the graph, and so we color the arrows to distinguish which homotopy they describe.
\end{definition}

The module $R\{y_j\}_{j=1}^{c}$ is the divided power algebra on $y_1, \dots, y_c$, but only its module structure and the $R[t_j]_{j=1}^{c}$-action above are used here. The CI operators are also called Eisenbud operators or cohomology operators \cite{AvramovSixLectures}. Every $y^{\boldsymbol{u}}$ has even homological degree, so tensor factors from $R\{y_j\}_{j=1}^{c}$ contribute no Koszul signs anywhere below.

\begin{example}\label{ex:buildup}
    We continue Example~\ref{ex:runtaylor} with the single element $a = xy = m_1$, so $\Rbar = R/(xy)$. Over $R$ we have $d(\eps_1) = xy \eps_\varnothing$. Over $\Rbar$ this is $0$, so $\eps_1$ is a cycle of $\F \otimes_R \Rbar$. It is not a boundary, because every entry of $d_2$ lies in $(x, y, z, w)$. So reducing $\F$ modulo $a$ introduces homology, and $\F \otimes_R \Rbar$ is not a resolution of $R/I$ over $\Rbar$.

    To kill this homology we add a new basis element $y_1 \eps_\varnothing$ in homological degree $2$ with boundary $\eps_1$. The same happens at every cell. So for every cell $\alpha$ we add a new cell $y_1 \alpha$, two homological degrees above $\alpha$, and we look for a map $\sigma$ on $\F$, raising homological degree by one, with $\sigma(\eps_\varnothing) = \eps_1$ and
    \[
        d(y_1 \alpha) = y_1 d(\alpha) + \sigma(\alpha).
    \]
    Then $d^2(y_1 \alpha) = \sigma(d\alpha) + d(\sigma\alpha)$, so the two copies $\F$ and $y_1 \F$, glued by $\sigma$, form a complex over $\Rbar$ as soon as $d\sigma + \sigma d = a \id$.

    Killing that homology introduces new homology, one copy up. Suppose first that $\sigma^2 = 0$. Then $\sigma(\eps_1) = \sigma^2(\eps_\varnothing) = 0$, so $d(y_1 \eps_1) = y_1 \, xy \, \eps_\varnothing$, which is $0$ over $\Rbar$, and again $y_1 \eps_1$ is not a boundary. We kill it with a new cell $y_1^{(2)} \eps_\varnothing$ with boundary $y_1 \eps_1$, and for every cell $\alpha$ we add a new cell $y_1^{(2)} \alpha$ with
    \[
        d(y_1^{(2)} \alpha) = y_1^{(2)} d(\alpha) + y_1 \sigma(\alpha).
    \]
    This introduces new homology again, which we kill with a fourth copy, and so on. The process never stops, and the result is one copy $y_1^{(k)} \F$ for every $k \ge 0$.

    If $\sigma^2 \neq 0$, the last formula no longer squares to zero over $\Rbar$, because $$d^2(y_1^{(2)} \alpha) = y_1 (d\sigma + \sigma d)(\alpha) + \sigma^2(\alpha).$$ We correct it by a further map $\tau$ on $\F$, raising homological degree by three, setting $$d(y_1^{(2)} \alpha) = y_1^{(2)} d(\alpha) + y_1 \sigma(\alpha) + \tau(\alpha),$$ and this squares to zero over $\Rbar$ when $d\tau + \tau d + \sigma^2 = 0$. Each further copy may need a further map.

    The construction above does not use a grading. In this example $a$ is a monomial, and the fine grading can be kept. If $y_1$ is given the multidegree $xy$ of $a$ and $\sigma$ raises multidegree by $a$, then $d$ is multigraded, and the copy $y_1^{(k)} \F$ is $\F$ shifted in multidegree by $a^k$.

    The next definition collects the maps $d, \sigma, \tau, \dots$ and the identities they must satisfy, for any number of elements $a_1, \dots, a_c$.
\end{example}

\begin{definition}\label{def:higher}
    Let $\F$ be a free resolution over $R$ and $a_1, \dots, a_c \in R$. A \emph{system of higher homotopies} for $a_1, \dots, a_c$ on $\F$ is a family $\{\sigma_{\boldsymbol{u}}\}_{\boldsymbol{u} \in \N^c}$ of $R$-linear maps $\sigma_{\boldsymbol{u}} \colon F_\bullet \to F_{\bullet + 2|\boldsymbol{u}| - 1}$ with $\sigma_{\boldsymbol{0}} = d$ and
    \begin{equation}\label{eq:higher}
        \sum_{\boldsymbol{v} + \boldsymbol{v'} = \boldsymbol{u}} \sigma_{\boldsymbol{v}} \sigma_{\boldsymbol{v'}} =
        \begin{cases}
            a_j \id & \text{if } \boldsymbol{u} = \boldsymbol{e}_j,\\
            0 & \text{otherwise.}
        \end{cases}
    \end{equation}

    A system of higher homotopies is \emph{truncated} if $\sigma_{\boldsymbol{u}} = 0$ for all $\boldsymbol{u}$ with $|\boldsymbol{u}| \ge 2$. Write $\sigma^{(j)} = \sigma_{\boldsymbol{e}_j}$. For a truncated system, \eqref{eq:higher} is equivalent to
    \begin{enumerate}
        \item[(a)] $d^2 = 0$ at $\boldsymbol{u} = 0$,
        \item[(b)] $d\sigma^{(j)} + \sigma^{(j)} d = a_j \id$ at $\boldsymbol{u} = \boldsymbol{e}_j$,
        \item[(c)] $(\sigma^{(j)})^2 = 0$ at $\boldsymbol{u} = 2\boldsymbol{e}_j$ and $\sigma^{(j)}\sigma^{(k)} + \sigma^{(k)}\sigma^{(j)} = 0$ at $\boldsymbol{u} = \boldsymbol{e}_j + \boldsymbol{e}_k$,
    \end{enumerate}
    while at $|\boldsymbol{u}| \ge 3$ every term has a factor $\sigma_{\boldsymbol{v}}$ with $|\boldsymbol{v}| \ge 2$ and so vanishes.
    We call $| \boldsymbol{u} |$ a \emph{level}, and refer to the largest $| \boldsymbol{u}|$ that appears in a system of higher homotopies as its level.
\end{definition}

\begin{remark}\label{rem:length}
    A system of higher homotopies on a complex of finite length must be of finite level, because $\sigma_{\boldsymbol{u}}$, which raises homological degree by $2|\boldsymbol{u}| - 1$, must eventually be a map into the $0$ module. For example, on any complex of length $\ell \le 2$ every system of higher homotopies is truncated. 
\end{remark}

\begin{remark}\label{rem:updown}
    In the complex graph of Definition~\ref{def:based}, draw an arrow $\alpha \to \beta$ of weight $\sigma^{(j)}_{\alpha,\beta}$, going up one homological degree, for every nonzero entry of $\sigma^{(j)}$. The identity $d\sigma^{(j)} + \sigma^{(j)} d = a_j \id$ then reads as follows. From a cell $\alpha$, go up along an arrow of $\sigma^{(j)}$ and then down along an edge of $d$, or down along an edge of $d$ and then up along an arrow of $\sigma^{(j)}$. The sum of the weights of these two step paths from $\alpha$ to a cell $\beta$ is $a_j$ if $\beta = \alpha$ and $0$ otherwise. Over $\Rbar$ every one of these sums is $0$.
\end{remark}

We define a free resolution of $R/I$ over $\Rbar = R/(a_1,\dots,a_c)$ by finding a free resolution $\F$ of $R/I$ over $R$, writing one copy $y^{\boldsymbol{u}} \F$ of it for each $\boldsymbol{u} \in \N^c$, that is, for each monomial $a_1^{u_1} \cdots a_c^{u_c}$ in the $a_j$, placed $2|\boldsymbol{u}|$ homological degrees above $\F$, and gluing them together carefully. The $t_j$ from the $R[t_j]_{j=1}^c$ behave as a shift of homological degree $-2$ on each of these free resolutions.

\begin{definition}\label{def:shamash}
    \cite[Construction 2.4]{Sobieska} For a system of higher homotopies $\sigma = \{\sigma_{\boldsymbol{u}}\}_{\boldsymbol{u} \in \N^c}$ for $a_1, \dots, a_c$ on $\F$, the \emph{Eisenbud--Shamash resolution} is
    \[
        \Sh(\F, \sigma) = R\{y_j\}_{j=1}^{c} \otimes_R \F \otimes_R \Rbar, \qquad d_{\Sh} = \sum_{\boldsymbol{u} \in \N^c} t^{\boldsymbol{u}} \otimes \sigma_{\boldsymbol{u}} \otimes \id .
    \]
\end{definition}

We write
\[
    D = \sum_{\boldsymbol{u} \in \N^c} t^{\boldsymbol{u}} \otimes \sigma_{\boldsymbol{u}} \colon R\{y_j\}_{j=1}^{c} \otimes_R \F \to R\{y_j\}_{j=1}^{c} \otimes_R \F ,
\]
so that $d_{\Sh} = \overline{D}$, in the following notation. A bar means reduction modulo $\aaa$, that is, $- \otimes_R \Rbar$. So $\overline{F}_\bullet = \F \otimes_R \Rbar$ is free over $\Rbar$ with basis $\bar\alpha = \alpha \otimes \overline{1}$ for $\alpha \in \B$, and an $R$-linear map $g$ has reduction $\bar g = g \otimes \id$. Then $\Sh(\F, \sigma) = R\{y_j\}_{j=1}^{c} \otimes_R \overline{F}_\bullet$, with basis the $y^{\boldsymbol{u}} \otimes \bar\alpha$, and $d_{\Sh} = \overline{D} = \sum_{\boldsymbol{u}} t^{\boldsymbol{u}} \otimes \bar\sigma_{\boldsymbol{u}}$. The map $D$ over $R$ describes the differential of a curved complex \cite[Chapter 3]{LodayVallette}, while $\bar{D}$ is genuinely a differential over $\bar{R}$.

\begin{remark}
    We will denote $D^2$ by $\Theta$, and this will later be an important shorthand for computation in our proofs in Section \ref{sec:commute}.

    We will provide a block-structure for the matrix of $D$ in Section \ref{sec:block}.
\end{remark}
Although we call this resolution the Eisenbud--Shamash resolution, we will denote the differential $d_{\Sh}$, the Shamash differential. It might not be apparent that it is a resolution. This fact will follow from the following theorem due to Shamash \cite{Shamash} for $c = 1$ and Eisenbud \cite[Theorem 7.2]{Eisenbud1980} in general. We have stated it as it appears in \cite[Proposition~4.1.4]{EisenbudPeeva} and \cite[Theorem~2.5]{Sobieska}.

\begin{theorem}[{\cite{Eisenbud1980,EisenbudPeeva,Shamash,Sobieska}}]\label{thm:shamash}
    Let $\F$ be a free resolution of $R/I$ over $R$, let $a_1, \dots, a_c$ be an $R$-regular sequence in $I$, and let $\sigma$ be a system of higher homotopies for $a_1, \dots, a_c$ on $\F$. Then $\Sh(\F, \sigma)$ is a free $\Rbar$-resolution of $R/I$.
\end{theorem}

\begin{remark}\label{rem:shnontrunc}
    The sum defining $d_{\Sh}$ is finite on each element, because $t^{\boldsymbol{u}}$ kills $y^{\boldsymbol{v}}$ unless $\boldsymbol{u} \le \boldsymbol{v}$. So $d_{\Sh}$ is well defined for every system of higher homotopies.
\end{remark}

\begin{remark}\label{rem:shtrunc}
    For a truncated system $\sigma = \{\sigma^{(j)}\}_{j=1}^c$ satisfying (a)-(c) in Definition \ref{def:higher} only $\boldsymbol{u} = \boldsymbol{0}$ and $\boldsymbol{u} = \boldsymbol{e}_j$ contribute, and
    \[
        d_{\Sh} = \id \otimes d \otimes \id + \sum_{j=1}^{c} t_j \otimes \sigma^{(j)} \otimes \id .
    \]
\end{remark}

\begin{remark}[{\cite[Remark~3.4]{Sobieska}, \cite{AvramovSixLectures}}] \label{rem:mincond}
    Let $\mathfrak{m} = (x_1, \dots, x_n)$, and suppose that $I$ and $a_1, \dots, a_c$ are homogeneous. When $\F$ is the minimal free resolution of $R/I$ and $a_j \in \mathfrak{m} I$ for every $1 \leq j \leq c$, then $\Sh(\F,\sigma)$ is minimal for any system of higher homotopies $\sigma$ on $\F$ for $(a_1,\dots,a_c)$.
\end{remark}

\begin{example}\label{ex:nonmonomial}
    When the regular sequence is homogeneous but not monomial, the resolution over $\Rbar$ is graded but no longer multigraded. In Example~\ref{ex:buildup} replace $a = xy$ by
    \[
        a = x^2y + xyz = (x + z)\, xy ,
    \]
    which is homogeneous of degree $3$ but not a monomial. If $\sigma$ is a homotopy for $xy$ on $\F$ with $\sigma^2 = 0$, then $(x + z)\sigma$ is one for $a$, and it sends $\eps_\varnothing$ to $(x + z)\eps_1$. Over $\Rbar = R/(x^2y + xyz)$ the new cell $y_1 \eps_\varnothing$ therefore has boundary $(x + z)\eps_1$, and
    \[
        d_{\Sh}|_2 = \left[\begin{array}{c|c} \overline{d_2} & \overline{(x+z)\sigma|_0} \end{array}\right] = \ \ \begin{NiceArray}{cccccc|cc}
        \eps_{12} & \eps_{13} & \eps_{14} & \eps_{23} & \eps_{24} & \eps_{34} & y_1\eps_{\varnothing} &\\
        z & zw & w & 0 & 0 & 0 & x + z & \eps_1\\
        -x & 0 & 0 & w & xw & 0 & 0 & \eps_2\\
        0 & -xy & 0 & -y & 0 & x & 0 & \eps_3\\
        0 & 0 & -y & 0 & -yz & -z & 0 & \eps_4
        \CodeAfter
        \SubMatrix[{2-1}{5-7}]
        \end{NiceArray}
    \]
    The entry $x + z$ is not a monomial, so no multidegree of $y_1 \eps_\varnothing$ makes this matrix multigraded. If $\sigma$ is homogeneous, then so is $(x + z)\sigma$, and the complex is still graded in the standard grading, with $y_1$ in degree $3$.
\end{example}

\section{The Shamash differential}\label{sec:shdiff}

\subsection{Decomposing the Shamash differential}\label{sec:decomp}

For a truncated system of higher homotopies, the differential of the Eisenbud--Shamash resolution has $c + 1$ pieces, one for the base and one for each color. This subsection extends Remark~\ref{rem:shtrunc}, which will help make sense of the combinatorics of the Eisenbud--Shamash resolution.

\begin{definition}
\label{def:delcon}
    For a truncated system $\sigma = \{\sigma^{(j)}\}_{j=1}^c$ satisfying (a)-(c) in Definition \ref{def:higher}, we split $D$ into the following maps.
    \begin{itemize}
        \item $\partial_0 = \id \otimes d$.
        \item $\partial_j = t_j \otimes \sigma^{(j)}$ for $1 \leq j \leq c$, the \emph{$j$th colored piece of the Shamash differential}. It lowers homological degree by $1$, connecting the $\boldsymbol{u}$ copy of the base resolution in the Eisenbud--Shamash resolution to the $\boldsymbol{u}-\boldsymbol{e}_j$ copy.
        \item $\dcon = \partial_1 + \dots + \partial_c$, the \emph{connecting part of the Shamash differential}.
    \end{itemize}
    Hence, $D = \partial_0 + \dcon$ and $d_{\Sh} = \bar\partial_0 + \bar\partial_{\mathrm{con}}$.
\end{definition}

\begin{remark}\label{rem:piecessquare}
    Since every $t_j$ has even degree, no signs occur. Thus $\partial_j \partial_k = t_j t_k \otimes \sigma^{(j)} \sigma^{(k)}$ and $\partial_0 \partial_j + \partial_j \partial_0 = t_j \otimes (d\sigma^{(j)} + \sigma^{(j)} d)$. The identities of a truncated system (Definition~\ref{def:higher}), therefore give on $R\{y_j\}_{j=1}^{c} \otimes_R \F$ before the tensor with $\Rbar$,
    \[
        \partial_0^2 = 0, \qquad \dcon^2 = \sum_{j} \partial_j^2 + \sum_{j < k} (\partial_j \partial_k + \partial_k \partial_j) = 0, \qquad \partial_0 \dcon + \dcon \partial_0 = \sum_j t_j \otimes a_j \id.
    \]
    Hence,
    \[
        \Theta = D^2 = \partial_0^2 + (\partial_0 \dcon + \dcon \partial_0) + \dcon^2 = \sum_{j=1}^{c} t_j \otimes a_j \id ,
    \]
    and $d_{\Sh}^2 = \overline{D^2} = 0$ since each $a_j$ is zero in $\Rbar$.
    By this identity the CI operators $t_j$ of $R\{y_j\}_{j=1}^c$ induce the maps $t_j \otimes \id$ which are the CI operators of $\Sh(\F, \sigma)$ for the lift $D$ of $d_{\Sh}$, in the sense of \cite[Construction~4.1.1]{EisenbudPeeva}.
\end{remark}

\begin{figure}[htbp]
    \centering
    \resizebox{\textwidth}{!}{%
    \begin{tikzpicture}[>=Stealth, x=2.1cm, y=1.25cm,
        fnode/.style={inner sep=1.5pt, font=\scriptsize},
        op/.style={font=\tiny, inner sep=1pt},
        diaglbl/.style={font=\tiny, inner sep=1pt, fill=white},
        rowlbl/.style={anchor=east, font=\tiny},
        mat/.style={font=\scriptsize, inner sep=1pt, anchor=center}]
        \node[rowlbl] at (-1.95, 2) {$y_1^{(0)}$:};
        \node[fnode] (A6) at (-1, 2) {$F_6$};
        \node[fnode] (A5) at (0, 2) {$F_5$};
        \node[fnode] (A4) at (1, 2) {$F_4$};
        \node[fnode] (A3) at (2, 2) {$F_3$};
        \node[fnode] (A2) at (3, 2) {$F_2$};
        \node[fnode] (A1) at (4, 2) {$F_1$};
        \node[fnode] (A0) at (5, 2) {$F_0$};
        \node[rowlbl] at (-1.95, 1) {$y_1^{(1)}$:};
        \node[fnode] (B4) at (-1, 1) {$y_1 F_4$};
        \node[fnode] (B3) at (0, 1) {$y_1 F_3$};
        \node[fnode] (B2) at (1, 1) {$y_1 F_2$};
        \node[fnode] (B1) at (2, 1) {$y_1 F_1$};
        \node[fnode] (B0) at (3, 1) {$y_1 F_0$};
        \node[rowlbl] at (-1.95, 0) {$y_1^{(2)}$:};
        \node[fnode] (C2) at (-1, 0) {$y_1^{(2)} F_2$};
        \node[fnode] (C1) at (0, 0) {$y_1^{(2)} F_1$};
        \node[fnode] (C0) at (1, 0) {$y_1^{(2)} F_0$};
        \draw[->] (A6) -- node[op,above]{$d_6$} (A5);
        \draw[->] (A5) -- node[op,above]{$d_5$} (A4);
        \draw[->] (A4) -- node[op,above]{$d_4$} (A3);
        \draw[->] (A3) -- node[op,above]{$d_3$} (A2);
        \draw[->] (A2) -- node[op,above]{$d_2$} (A1);
        \draw[->] (A1) -- node[op,above]{$d_1$} (A0);
        \draw[->] (B4) -- node[op,above]{$d_4$} (B3);
        \draw[->] (B3) -- node[op,above]{$d_3$} (B2);
        \draw[->] (B2) -- node[op,above]{$d_2$} (B1);
        \draw[->] (B1) -- node[op,above]{$d_1$} (B0);
        \draw[->] (C2) -- node[op,above]{$d_2$} (C1);
        \draw[->] (C1) -- node[op,above]{$d_1$} (C0);
        \draw[->] (B4) -- node[diaglbl, midway]{$\sigma|_4$} (A5);
        \draw[->] (B3) -- node[diaglbl, midway]{$\sigma|_3$} (A4);
        \draw[->] (B2) -- node[diaglbl, midway]{$\sigma|_2$} (A3);
        \draw[->] (B1) -- node[diaglbl, midway]{$\sigma|_1$} (A2);
        \draw[->] (B0) -- node[diaglbl, midway]{$\sigma|_0$} (A1);
        \draw[->] (C2) -- node[diaglbl, midway]{$\sigma|_2$} (B3);
        \draw[->] (C1) -- node[diaglbl, midway]{$\sigma|_1$} (B2);
        \draw[->] (C0) -- node[diaglbl, midway]{$\sigma|_0$} (B1);
        \node[font=\scriptsize, anchor=east] at (-1.55, 2) {$\cdots$};
        \node[font=\scriptsize, anchor=east] at (-1.55, 1) {$\cdots$};
        \node[font=\scriptsize, anchor=east] at (-1.55, 0) {$\cdots$};
        \begin{scope}[on background layer]
            \fill[sh0, rounded corners=4pt] (-2.35, 1.62) rectangle (5.35, 2.38);
            \fill[sh1, rounded corners=4pt] (-2.35, 0.62) rectangle (3.4, 1.38);
            \fill[sh2, rounded corners=4pt] (-2.35, -0.38) rectangle (1.45, 0.38);
        \end{scope}
        \node[mat] at (-0.32, -1.05) {$\left[\begingroup\renewcommand{\arraystretch}{1.3}\begin{array}{c|c|c} \shent{sh0}{\overline{d_5}} & \shent{sh1}{\overline{\sigma|_3}} & 0\\
        \hline
        0 & \shent{sh1}{\overline{d_3}} & \shent{sh2}{\overline{\sigma|_1}}\\
        \hline
        0 & 0 & \shent{sh2}{\overline{d_1}}\end{array}\endgroup\right],$};
        \node[mat] at (1.2, -1.05) {$\left[\begingroup\renewcommand{\arraystretch}{1.3}\begin{array}{c|c|c} \shent{sh0}{\overline{d_4}} & \shent{sh1}{\overline{\sigma|_2}} & 0 \\
        \hline
        0 & \shent{sh1}{\overline{d_2}} & \shent{sh2}{\overline{\sigma|_0}} \end{array}\endgroup\right],$};
        \node[mat] at (2.45, -1.05) {$\left[\begingroup\renewcommand{\arraystretch}{1.3}\begin{array}{c|c} \shent{sh0}{\overline{d_3}} & \shent{sh1}{\overline{\sigma|_1}} \\
        \hline
        0 & \shent{sh1}{\overline{d_1}} \end{array}\endgroup\right],$};
        \node[mat] at (3.45, -1.05) {$\left[\begin{array}{c|c} \shent{sh0}{\overline{d_2}} & \shent{sh1}{\overline{\sigma|_0}} \end{array}\right],$};
        \node[mat] at (4.1, -1.05) {$\begin{bmatrix} \shent{sh0}{\overline{d_1}} \end{bmatrix}$};
    \end{tikzpicture}}
    \caption{Shamash resolution of a truncated system, $c=1$.}
    \label{fig:Sham1}
\end{figure}

Figure~\ref{fig:Sham1} shows the Shamash resolution of a truncated system for $c=1$, with $\sigma = \sigma^{(1)}$, the copy $y_1^{(i)} F_q$ in row $i$ and column $n = 2i + q$, and $\bar\partial_0$ and $\bar\partial_1$ the horizontal and diagonal arrows. Below it are the matrices of $d_{\Sh}|_n$ for $n = 5, \dots, 1$, with columns the copies in column $n$, rows those in column $n - 1$, and $\sigma|_q$ the restriction of $\sigma$ to $F_q$. Each matrix is block bidiagonal, corresponding to having a truncated system of higher homotopies. The general block structure has more superdiagonal terms because arrows do not stay between adjacent rows. We make this precise in Proposition \ref{prop:blocks}.

\begin{figure}[htbp]
    \centering
    \resizebox{0.55\textwidth}{!}{%
    \begin{tikzpicture}[>=Stealth, x=3.4cm, y=1.75cm,
        snode/.style={inner sep=2pt, font=\scriptsize},
        op/.style={font=\scriptsize, inner sep=1pt},
        a1/.style={->, thick, shorten >=2pt, shorten <=2pt},
        a2/.style={->, thick, densely dashed, shorten >=2pt, shorten <=2pt}]
        \node[snode] (u00) at (0,0) {$\F$};
        \node[snode] (u10) at (1,0) {$y_1 \F$};
        \node[snode] (u20) at (2,0) {$y_1^{(2)} \F$};
        \node[snode] (u01) at (0,1) {$y_2 \F$};
        \node[snode] (u11) at (1,1) {$y_1 y_2 \F$};
        \node[snode] (u21) at (2,1) {$y_1^{(2)} y_2 \F$};
        \node[snode] (u02) at (0,2) {$y_2^{(2)} \F$};
        \node[snode] (u12) at (1,2) {$y_1 y_2^{(2)} \F$};
        \node[snode] (u22) at (2,2) {$y_1^{(2)} y_2^{(2)} \F$};
        \draw[a1] (u10) -- node[op,below]{$\bar\partial_1$} (u00);
        \foreach \a/\b in {u20/u10,u11/u01,u21/u11,u12/u02,u22/u12}{\draw[a1] (\a) -- (\b);}
        \draw[a2] (u01) -- node[op,left]{$\bar\partial_2$} (u00);
        \foreach \a/\b in {u02/u01,u11/u10,u12/u11,u21/u20,u22/u21}{\draw[a2] (\a) -- (\b);}
        \node[font=\scriptsize] at (2.62,0) {$\cdots$};
        \node[font=\scriptsize] at (2.62,1) {$\cdots$};
        \node[font=\scriptsize] at (2.62,2) {$\cdots$};
        \node[font=\scriptsize] at (0,2.55) {$\vdots$};
        \node[font=\scriptsize] at (1,2.55) {$\vdots$};
        \node[font=\scriptsize] at (2,2.55) {$\vdots$};
    \end{tikzpicture}}
    \caption{Shamash complex for $c=2$. }
    \label{fig:Sham2_g}
\end{figure}

Figure \ref{fig:Sham2_g} illustrates a truncated system for $c=2$. Solid arrows are $\bar\partial_1$, from the node at $\boldsymbol{u}$ to the node at $\boldsymbol{u} - \boldsymbol{e}_1$ of the slice behind, and dashed arrows are $\bar\partial_2$, to the node at $\boldsymbol{u} - \boldsymbol{e}_2$. The arrows of $\bar\partial_0$ go straight back from each node to the node at the same $\boldsymbol{u}$ and are not drawn. Every square anticommutes, because $\sigma^{(1)}\sigma^{(2)} + \sigma^{(2)}\sigma^{(1)} = 0$.

A general system adds an arrow from the copy over $\boldsymbol{u}$ to the copy over $\boldsymbol{u} - \boldsymbol{v}$ for every $\boldsymbol{v} \le \boldsymbol{u}$ with $|\boldsymbol{v}| \ge 2$, so the lattice has additional diagonal arrows.
The squares then no longer anticommute. Taking $\boldsymbol{u} = (1,1)$ in \eqref{eq:higher} of Definition \ref{def:higher} gives $\sigma^{(1)}\sigma^{(2)} + \sigma^{(2)}\sigma^{(1)} = -(d\sigma_{(1,1)} + \sigma_{(1,1)}d)$, so
\[
    \partial_1 \partial_2 + \partial_2 \partial_1 = -t_1 t_2 \otimes (d\sigma_{(1,1)} + \sigma_{(1,1)}d)
\]
instead of zero, and the diagonal arrow $t_1 t_2 \otimes \sigma_{(1,1)}$ accounts for the difference.

\subsection{Block form of the Shamash differential}\label{sec:block}

The Eisenbud--Shamash resolution takes a copy of a finite free resolution $\F$ for each lattice point $\boldsymbol{u} \in \N^c$, and groups them by level. Fix a level, $k$, and we should think of these lattice points as the points that are on the hyperplane $u_1 + \dots + u_c = k$, $k$ steps from the origin in the nonnegative orthant. The copy at each of these points is shifted up by $2k$ homological degrees, so consecutive hyperplanes are two homological degrees apart. In this subsection, we flatten this picture so that for each homological degree $n$ and each hyperplane we determine which homological degree of the free resolutions that live there contributes to the degree $n$ part of the Eisenbud--Shamash resolution. This allows us to make precise the block structure of the Shamash differential.

For $k \ge 0$ put $N_k = \#\{\boldsymbol{u} \in \N^c : |\boldsymbol{u}| = k\} = \binom{k+c-1}{c-1}$. Gathering summands by level $k = |\boldsymbol{u}|$,
\[
    (R\{y_j\}_{j=1}^{c} \otimes_R \F)_n = \bigoplus_{2k \le n} L_k(n), \qquad L_k(n) = \bigoplus_{|\boldsymbol{u}| = k} y^{\boldsymbol{u}} \otimes F_{n-2k} ,
\]
and $\Sh(\F, \sigma)_n = \bigoplus_{2k \le n} L_k(n) \otimes_R \Rbar$, with matrix entries read modulo $\aaa$. The copy at a point $\boldsymbol{u}$ with $|\boldsymbol{u}| = k$ contributes $y^{\boldsymbol{u}} \otimes F_{n-2k}$ to the degree $n$ part, as in Figure~\ref{fig:Sham2}.

\begin{figure}[htbp]
    \centering
    \resizebox{0.72\textwidth}{!}{%
    \begin{tikzpicture}[>=Stealth, x=1.15cm, y=1.15cm,
        snode/.style={inner sep=1.5pt, font=\scriptsize},
        znode/.style={inner sep=1.5pt, font=\scriptsize, text=gray},
        a1/.style={->, thick, shorten >=2pt, shorten <=2pt},
        a2/.style={->, thick, densely dashed, shorten >=2pt, shorten <=2pt}]
        \begin{scope}[shift={(0,0)}]
            \node[snode] (n00) at (0,0) {$F_{4}$};
            \node[snode] (n01) at (0,1) {$F_{2}$};
            \node[snode] (n02) at (0,2) {$F_{0}$};
            \node[snode] (n10) at (1,0) {$F_{2}$};
            \node[snode] (n11) at (1,1) {$F_{0}$};
            \node[snode] (n20) at (2,0) {$F_{0}$};
            \draw[a2] (n01) -- (n00);
            \draw[a2] (n02) -- (n01);
            \draw[a1] (n10) -- (n00);
            \draw[a1] (n11) -- (n01);
            \draw[a2] (n11) -- (n10);
            \draw[a1] (n20) -- (n10);
            \node[font=\scriptsize] at (1.0,-0.8) {$n=4$};
        \end{scope}
        \begin{scope}[shift={(3.9,0)}]
            \node[znode] (n00) at (0,0) {$0$};
            \node[snode] (n01) at (0,1) {$F_{3}$};
            \node[snode] (n02) at (0,2) {$F_{1}$};
            \node[snode] (n10) at (1,0) {$F_{3}$};
            \node[snode] (n11) at (1,1) {$F_{1}$};
            \node[snode] (n20) at (2,0) {$F_{1}$};
            \draw[a2] (n01) -- (n00);
            \draw[a2] (n02) -- (n01);
            \draw[a1] (n10) -- (n00);
            \draw[a1] (n11) -- (n01);
            \draw[a2] (n11) -- (n10);
            \draw[a1] (n20) -- (n10);
            \node[font=\scriptsize] at (1.0,-0.8) {$n=5$};
        \end{scope}
        \begin{scope}[shift={(7.8,0)}]
            \node[znode] (n00) at (0,0) {$0$};
            \node[snode] (n01) at (0,1) {$F_{4}$};
            \node[snode] (n02) at (0,2) {$F_{2}$};
            \node[snode] (n03) at (0,3) {$F_{0}$};
            \node[snode] (n10) at (1,0) {$F_{4}$};
            \node[snode] (n11) at (1,1) {$F_{2}$};
            \node[snode] (n12) at (1,2) {$F_{0}$};
            \node[snode] (n20) at (2,0) {$F_{2}$};
            \node[snode] (n21) at (2,1) {$F_{0}$};
            \node[snode] (n30) at (3,0) {$F_{0}$};
            \draw[a2] (n02) -- (n01);
            \draw[a2] (n03) -- (n02);
            \draw[a1] (n11) -- (n01);
            \draw[a2] (n11) -- (n10);
            \draw[a1] (n12) -- (n02);
            \draw[a2] (n12) -- (n11);
            \draw[a1] (n20) -- (n10);
            \draw[a1] (n21) -- (n11);
            \draw[a2] (n21) -- (n20);
            \draw[a1] (n30) -- (n20);
            \node[font=\scriptsize] at (1.5,-0.8) {$n=6$};
        \end{scope}
        \begin{scope}[shift={(0,-5.6)}]
            \node[znode] (n00) at (0,0) {$0$};
            \node[znode] (n01) at (0,1) {$0$};
            \node[snode] (n02) at (0,2) {$F_{3}$};
            \node[snode] (n03) at (0,3) {$F_{1}$};
            \node[znode] (n10) at (1,0) {$0$};
            \node[snode] (n11) at (1,1) {$F_{3}$};
            \node[snode] (n12) at (1,2) {$F_{1}$};
            \node[snode] (n20) at (2,0) {$F_{3}$};
            \node[snode] (n21) at (2,1) {$F_{1}$};
            \node[snode] (n30) at (3,0) {$F_{1}$};
            \draw[a2] (n02) -- (n01);
            \draw[a2] (n03) -- (n02);
            \draw[a1] (n11) -- (n01);
            \draw[a2] (n11) -- (n10);
            \draw[a1] (n12) -- (n02);
            \draw[a2] (n12) -- (n11);
            \draw[a1] (n20) -- (n10);
            \draw[a1] (n21) -- (n11);
            \draw[a2] (n21) -- (n20);
            \draw[a1] (n30) -- (n20);
            \node[font=\scriptsize] at (1.5,-0.8) {$n=7$};
        \end{scope}
        \begin{scope}[shift={(4.9,-5.6)}]
            \node[znode] (n00) at (0,0) {$0$};
            \node[znode] (n01) at (0,1) {$0$};
            \node[snode] (n02) at (0,2) {$F_{4}$};
            \node[snode] (n03) at (0,3) {$F_{2}$};
            \node[snode] (n04) at (0,4) {$F_{0}$};
            \node[znode] (n10) at (1,0) {$0$};
            \node[snode] (n11) at (1,1) {$F_{4}$};
            \node[snode] (n12) at (1,2) {$F_{2}$};
            \node[snode] (n13) at (1,3) {$F_{0}$};
            \node[snode] (n20) at (2,0) {$F_{4}$};
            \node[snode] (n21) at (2,1) {$F_{2}$};
            \node[snode] (n22) at (2,2) {$F_{0}$};
            \node[snode] (n30) at (3,0) {$F_{2}$};
            \node[snode] (n31) at (3,1) {$F_{0}$};
            \node[snode] (n40) at (4,0) {$F_{0}$};
            \draw[a2] (n03) -- (n02);
            \draw[a2] (n04) -- (n03);
            \draw[a1] (n12) -- (n02);
            \draw[a2] (n12) -- (n11);
            \draw[a1] (n13) -- (n03);
            \draw[a2] (n13) -- (n12);
            \draw[a1] (n21) -- (n11);
            \draw[a2] (n21) -- (n20);
            \draw[a1] (n22) -- (n12);
            \draw[a2] (n22) -- (n21);
            \draw[a1] (n30) -- (n20);
            \draw[a1] (n31) -- (n21);
            \draw[a2] (n31) -- (n30);
            \draw[a1] (n40) -- (n30);
            \node[font=\scriptsize] at (2.0,-0.8) {$n=8$};
        \end{scope}
        \begin{scope}[shift={(10.8,-5.6)}]
            \node[znode] (n00) at (0,0) {$0$};
            \node[znode] (n01) at (0,1) {$0$};
            \node[znode] (n02) at (0,2) {$0$};
            \node[snode] (n03) at (0,3) {$F_{3}$};
            \node[snode] (n04) at (0,4) {$F_{1}$};
            \node[znode] (n10) at (1,0) {$0$};
            \node[znode] (n11) at (1,1) {$0$};
            \node[snode] (n12) at (1,2) {$F_{3}$};
            \node[snode] (n13) at (1,3) {$F_{1}$};
            \node[znode] (n20) at (2,0) {$0$};
            \node[snode] (n21) at (2,1) {$F_{3}$};
            \node[snode] (n22) at (2,2) {$F_{1}$};
            \node[snode] (n30) at (3,0) {$F_{3}$};
            \node[snode] (n31) at (3,1) {$F_{1}$};
            \node[snode] (n40) at (4,0) {$F_{1}$};
            \draw[a2] (n03) -- (n02);
            \draw[a2] (n04) -- (n03);
            \draw[a1] (n12) -- (n02);
            \draw[a2] (n12) -- (n11);
            \draw[a1] (n13) -- (n03);
            \draw[a2] (n13) -- (n12);
            \draw[a1] (n21) -- (n11);
            \draw[a2] (n21) -- (n20);
            \draw[a1] (n22) -- (n12);
            \draw[a2] (n22) -- (n21);
            \draw[a1] (n30) -- (n20);
            \draw[a1] (n31) -- (n21);
            \draw[a2] (n31) -- (n30);
            \draw[a1] (n40) -- (n30);
            \node[font=\scriptsize] at (2.0,-0.8) {$n=9$};
        \end{scope}
    \end{tikzpicture}}
    \caption{Degree $n$ slices of the Shamash resolution of a truncated system, $c = 2$, $\F$ of length $4$.}
    \label{fig:Sham2}
\end{figure}

Let $\boldsymbol{v} \in \N^c$ and $k \ge |\boldsymbol{v}|$. The columns of $T^{(k)}_{\boldsymbol{v}}$ are indexed by the $y^{\boldsymbol{u}}$ with $|\boldsymbol{u}| = k$, in decreasing order of $u_1$, then of $u_2$, and so on, and its rows by the $y^{\boldsymbol{u'}}$ with $|\boldsymbol{u'}| = k - |\boldsymbol{v}|$, in the same order. Its entry in column $y^{\boldsymbol{u}}$ and row $y^{\boldsymbol{u'}}$ is $1$ if $\boldsymbol{u} = \boldsymbol{u'} + \boldsymbol{v}$ and $0$ otherwise, so $T^{(k)}_{\boldsymbol{v}}$ is the matrix of $t^{\boldsymbol{v}}$ from level $k$ to level $k - |\boldsymbol{v}|$, and $T^{(k)}_0 = \mathbf{1}_{N_k}$ is the identity. The columns are the points of the plane at $k$, and the rows are those of the plane at $k - |\boldsymbol{v}|$, and $t^{\boldsymbol{v}}$ translates by $-\boldsymbol{v}$, sending the points $\boldsymbol{u} \ge \boldsymbol{v}$ to $\boldsymbol{u} - \boldsymbol{v}$ and the others out of $\N^c$, to $0$ (Figure~\ref{fig:slide}).

\begin{figure}[htbp]
   \centering
    \resizebox{\textwidth}{!}{%
    \begin{tikzpicture}[>=Stealth, x=1.3cm, y=1.2cm,
        pt/.style={inner sep=1pt, font=\scriptsize, fill=white},
        lv/.style={gray!55, thin},
        ax/.style={->, gray!70, thin},
        tk/.style={font=\tiny, gray!90},
        mv/.style={->, thick, shorten >=2pt, shorten <=2pt},
        cut/.style={->, thick, densely dashed, gray, shorten >=2pt, shorten <=2pt}]
        \begin{scope}
            \draw[ax] (-1.5,0) -- (3.8,0) node[right,font=\scriptsize] {$u_1$};
            \draw[ax] (0,-1.5) -- (0,3.7) node[above,font=\scriptsize] {$u_2$};
            \foreach \i in {-1,1,2,3}{ \draw[gray!70] (\i,-0.06) -- (\i,0.06); \node[tk,below] at (\i,-0.08) {$\i$}; }
            \foreach \i in {-1,1,2,3}{ \draw[gray!70] (-0.06,\i) -- (0.06,\i); \node[tk,left] at (-0.08,\i) {$\i$}; }
            \foreach \k in {1,2,3}{ \draw[lv] (\k,0) -- (0,\k); }
            \foreach \u/\v in {0/0,1/0,0/1,2/0,1/1,0/2,3/0,2/1,1/2,0/3}{ \node[pt] at (\u,\v) {$y^{(\u,\v)}$}; }
            \foreach \k in {1,2,3}{ \node[font=\tiny,text=gray,fill=white,inner sep=0.5pt,rotate=-42.7] at (\k-0.5,0.5) {$|\boldsymbol{u}|=\k$}; }
            \node[font=\scriptsize] at (1.3,-2) {(a) $t^{(0,0)}$};
        \end{scope}
        \begin{scope}[shift={(6.3,0)}]
            \draw[ax] (-1.5,0) -- (3.8,0) node[right,font=\scriptsize] {$u_1$};
            \draw[ax] (0,-1.5) -- (0,3.7) node[above,font=\scriptsize] {$u_2$};
            \foreach \i in {-1,1,2,3}{ \draw[gray!70] (\i,-0.06) -- (\i,0.06); \node[tk,below] at (\i,-0.08) {$\i$}; }
            \foreach \i in {-1,1,2,3}{ \draw[gray!70] (-0.06,\i) -- (0.06,\i); \node[tk,left] at (-0.08,\i) {$\i$}; }
            \foreach \k in {2,3}{ \draw[lv] (\k,0) -- (0,\k); }
            \foreach \u/\v in {3/0,2/1,1/2,0/3,2/0,1/1,0/2}{ \node[pt] (n\u\v) at (\u,\v) {$y^{(\u,\v)}$}; }
            \draw[mv] (n30) -- (n20); \draw[mv] (n21) -- (n11); \draw[mv] (n12) -- (n02);
            \node[pt,text=gray] (z) at (-1,3) {$0$};
            \draw[cut] (n03) -- (z);
            \node[font=\scriptsize] at (1.3,-2) {(b) $t^{(1,0)}$};
        \end{scope}
        \begin{scope}[shift={(12.6,0)}]
            \draw[ax] (-1.5,0) -- (3.8,0) node[right,font=\scriptsize] {$u_1$};
            \draw[ax] (0,-1.5) -- (0,3.7) node[above,font=\scriptsize] {$u_2$};
            \foreach \i in {-1,1,2,3}{ \draw[gray!70] (\i,-0.06) -- (\i,0.06); \node[tk,below] at (\i,-0.08) {$\i$}; }
            \foreach \i in {-1,1,2,3}{ \draw[gray!70] (-0.06,\i) -- (0.06,\i); \node[tk,left] at (-0.08,\i) {$\i$}; }
            \foreach \k in {1,3}{ \draw[lv] (\k,0) -- (0,\k); }
            \foreach \u/\v in {3/0,2/1,1/2,0/3,1/0,0/1}{ \node[pt] (n\u\v) at (\u,\v) {$y^{(\u,\v)}$}; }
            \draw[mv] (n21) -- (n10); \draw[mv] (n12) -- (n01);
            \node[pt,text=gray] (za) at (2,-1) {$0$};
            \node[pt,text=gray] (zb) at (-1,2) {$0$};
            \draw[cut] (n30) -- (za); \draw[cut] (n03) -- (zb);
            \node[font=\scriptsize] at (1.3,-2) {(c) $t^{(1,1)}$};
        \end{scope}
    \end{tikzpicture}}
    \caption{Sliding a plane by $t^{\boldsymbol{v}}$, $c = 2$.}
    \label{fig:slide}
\end{figure}

For the plane $|\boldsymbol{u}| = 3$ of Figure~\ref{fig:slide}, with columns $y^{(3,0)}$, $y^{(2,1)}$, $y^{(1,2)}$, $y^{(0,3)}$,
\[
    T^{(3)}_{(1,0)} = \begin{bmatrix} 1 & 0 & 0 & 0\\ 0 & 1 & 0 & 0\\ 0 & 0 & 1 & 0\end{bmatrix}, \qquad
    T^{(3)}_{(1,1)} = \begin{bmatrix} 0 & 1 & 0 & 0\\ 0 & 0 & 1 & 0\end{bmatrix} ,
\]
the zero columns being the points carried out of $\N^2$. Since $t^{\boldsymbol{w}} t^{\boldsymbol{v}} = t^{\boldsymbol{v}+\boldsymbol{w}}$, these matrices compose as $T^{(k-|\boldsymbol{v}|)}_{\boldsymbol{w}} T^{(k)}_{\boldsymbol{v}} = T^{(k)}_{\boldsymbol{v}+\boldsymbol{w}}$. For an $R$-linear map $f$, $T \otimes f$ is the Kronecker product, the block matrix with block $T_{\boldsymbol{u},\boldsymbol{u'}} f$ in the column of $y^{\boldsymbol{u}}$ and the row of $y^{\boldsymbol{u'}}$, and $(T \otimes f)(T' \otimes g) = TT' \otimes fg$.

\textbf{We fix the order of the basis here for the rest of the paper}. The cells $y^{\boldsymbol{u}} \otimes \alpha$ of a given homological degree are ordered first by level $|\boldsymbol{u}|$, starting with level $0$, then by $\boldsymbol{u}$, in decreasing order of $u_1$, then of $u_2$, and so on, and then by $\alpha$, in the order of the basis $\B$. For $c = 2$, for example, level $2$ is listed as $y_1^{(2)}, y_1 y_2, y_2^{(2)}$. With this order the matrix of $D$ is block upper triangular in the level, as the next proposition shows. The examples of \cite{Sobieska} list the levels in the opposite order, highest level first.

\begin{proposition}\label{prop:blocks}
    Let $\F$ be a free complex over $R$ with a finite basis of each $F_q$, let $a_1, \dots, a_c \in R$, and let $\sigma = \{\sigma_{\boldsymbol{u}}\}_{\boldsymbol{u} \in \N^c}$ be a system of higher homotopies for $a_1, \dots, a_c$ on $\F$. Let ${D = \sum_{\boldsymbol{u}} t^{\boldsymbol{u}} \otimes \sigma_{\boldsymbol{u}}}$, acting on $R\{y_j\}_{j=1}^{c} \otimes_R \F$, so that $d_{\Sh} = \overline{D}$. Write its matrix from $(R\{y_j\}_{j=1}^{c} \otimes_R \F)_n$ to $(R\{y_j\}_{j=1}^{c} \otimes_R \F)_{n-1}$ in blocks, with rows indexed by the target level $k'$ and columns by the source level $k$. Then:
    \begin{enumerate}[label=(\arabic*),leftmargin=*,itemsep=1pt]
        \item The blocks are
        \[ \bigl(D|_n\bigr)_{k,k'} = \sum_{|\boldsymbol{v}| = k - k'} T^{(k)}_{\boldsymbol{v}} \otimes \sigma_{\boldsymbol{v}}|_{n-2k} \quad (k' \le k), \qquad \bigl(D|_n\bigr)_{k,k'} = 0 \quad (k' > k),
        \]
        so the matrix is block upper triangular and its $m$th superdiagonal carries exactly the $\sigma_{\boldsymbol{v}}$ with $|\boldsymbol{v}| = m$.

        \item $\sigma$ is truncated if and only if the matrix of $D|_n$ is block bidiagonal for every $n$, in which case
        \[
            \bigl(D|_n\bigr)_{k,k} = \mathbf{1}_{N_k} \otimes d_{n-2k}, \qquad \bigl(D|_n\bigr)_{k,k-1} = S_k|_{n-2k} \quad \text{for} \quad S_k = \sum_{j=1}^{c} T^{(k)}_{\boldsymbol{e}_j} \otimes \sigma^{(j)} .
        \]
    \end{enumerate}
    The matrix of the differential of $\Sh(\F, \sigma)$ is the same, with its entries read modulo $\aaa$.
\end{proposition}

\begin{proof}
    Part (1): For $|\boldsymbol{u}| = k$ and $\alpha \in F_{n-2k}$, the term $t^{\boldsymbol{v}} \otimes \sigma_{\boldsymbol{v}}$ sends $y^{\boldsymbol{u}} \otimes \alpha$ to $y^{\boldsymbol{u}-\boldsymbol{v}} \otimes \sigma_{\boldsymbol{v}}(\alpha) \in L_{k-|\boldsymbol{v}|}(n-1)$ if $\boldsymbol{v} \le \boldsymbol{u}$ and to $0$ otherwise. So its block from $y^{\boldsymbol{u}} \otimes F_{n-2k}$ to $y^{\boldsymbol{u'}} \otimes F_{n-1-2k'}$ is $\sigma_{\boldsymbol{v}}|_{n-2k}$ if $\boldsymbol{u} = \boldsymbol{u'} + \boldsymbol{v}$ and $0$ otherwise, which is the $(\boldsymbol{u}, \boldsymbol{u'})$ block of $T^{(k)}_{\boldsymbol{v}} \otimes \sigma_{\boldsymbol{v}}|_{n-2k}$.

    Part (2): If $\sigma_{\boldsymbol{v}} = 0$ for all $|\boldsymbol{v}| \ge 2$ then only $k' = k$ and $k' = k - 1$ survive, and $T^{(k)}_0 = \mathbf{1}_{N_k}$ gives the diagonal blocks.

    Conversely, if $\sigma_{\boldsymbol{v}}|_q \neq 0$ with $|\boldsymbol{v}| = m \ge 2$, take $n = 2m + q$. The block $(m, 0)$ is ${\sum_{|\boldsymbol{v'}| = m} T^{(m)}_{\boldsymbol{v'}} \otimes \sigma_{\boldsymbol{v'}}|_q}$, where $T^{(m)}_{\boldsymbol{v'}}$ is the $1 \times N_m$ matrix with its single $1$ in column $\boldsymbol{v'}$, so the column of $\boldsymbol{v}$ carries $\sigma_{\boldsymbol{v}}|_q \neq 0$.
\end{proof}

If $\F$ has length $\ell$, then for $n > \ell$ a level $k$ has cells in homological degree $n$ only when $n - 2k \le \ell$, so the lowest levels are missing from $D|_n$. When $n - \ell$ is odd the columns of $D|_n$ start one level above its rows. Then $d$ sits below the diagonal, and the matrix, written without its empty blocks, looks block lower triangular.

By this proposition, we can write every block of an arbitrary Shamash differential down directly. We demonstrate this for $c=2$. The $(i, i')$ entry of $T^{(k)}_{(v_1, v_2)}$ is $1$ exactly when $i' = i + v_2$, reading rows and columns from $0$, so $T^{(k)}_{(v_1, v_2)}$ is a single band at offset $v_2$ and, for $0 \le m \le k$, the whole $m$th superdiagonal block, from the plane at level $k$ to the plane at level $k - m$, is the banded Toeplitz array
\[
    \bigl(D|_n\bigr)_{k,k-m} =
    \left[\arraycolsep=4pt\begin{array}{c|c|c|c|c|c|c}
        \sigma_{(m,0)} & \sigma_{(m-1,1)} & \sigma_{(m-2,2)} & \cdots & \sigma_{(0,m)} & & \\
        \hline
        & \sigma_{(m,0)} & \sigma_{(m-1,1)} & \cdots & \sigma_{(1,m-1)} &\sigma_{(0,m)} & \\
        \hline
        & & \ddots & \ddots & \ddots & \ddots & \ddots
    \end{array}\right]
\]
of size $(k-m+1) \times (k+1)$, the blank entries being zero and every entry restricted to homological degree $n - 2k$.

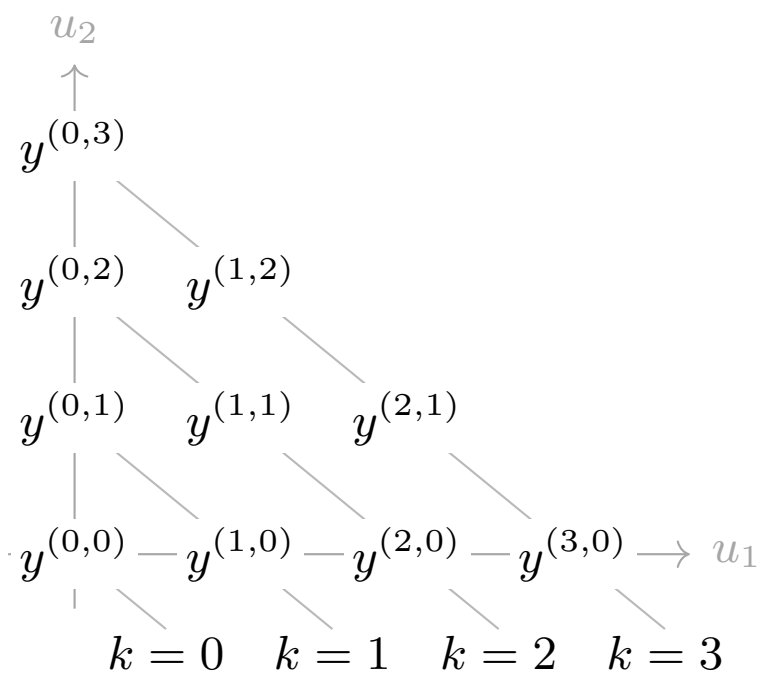
\begin{figure}[htbp]
    \centering
    \begin{tikzpicture}[x=1.1cm, y=0.9cm,
        pt/.style={inner sep=1.5pt, font=\scriptsize, fill=white}]
        \draw[->, gray!70, thin] (-0.4,0) -- (3.7,0) node[right,font=\scriptsize] {$u_1$};
        \draw[->, gray!70, thin] (0,-0.4) -- (0,3.6) node[above,font=\scriptsize] {$u_2$};
        \foreach \k in {0,1,2,3} {
            \draw[gray!55, thin] (0,\k) -- (\k+0.55,-0.55);
            \node[font=\scriptsize, anchor=north, inner sep=1pt] at (\k+0.55,-0.55) {$k = \k$};
        }
        \node[pt] at (0,0) {$y^{(0,0)}$};
        \node[pt] at (0,1) {$y^{(0,1)}$};
        \node[pt] at (0,2) {$y^{(0,2)}$};
        \node[pt] at (0,3) {$y^{(0,3)}$};
        \node[pt] at (1,0) {$y^{(1,0)}$};
        \node[pt] at (1,1) {$y^{(1,1)}$};
        \node[pt] at (1,2) {$y^{(1,2)}$};
        \node[pt] at (2,0) {$y^{(2,0)}$};
        \node[pt] at (2,1) {$y^{(2,1)}$};
        \node[pt] at (3,0) {$y^{(3,0)}$};
    \end{tikzpicture}

    \vspace{6pt}
    \resizebox{\textwidth}{!}{$\displaystyle
        \
        \def\arraystretch{2} \arraycolsep=4pt
        \begin{array}{c|cccc}
            & k = 0 & k = 1 & k = 2 & k = 3\\
            \hline
            k' = 0 & d_6 & \left[\def\arraystretch{1.3}\begin{array}{c|c}t_1\sigma_{(1,0)}|_4 & t_2\sigma_{(0,1)}|_4\end{array}\right] & \left[\def\arraystretch{1.3}\begin{array}{c|c|c}t_1^2\sigma_{(2,0)}|_2 & t_1t_2\sigma_{(1,1)}|_2 & t_2^2\sigma_{(0,2)}|_2\end{array}\right] & \left[\def\arraystretch{1.3}\begin{array}{c|c|c|c}t_1^3\sigma_{(3,0)}|_0 & t_1^2t_2\sigma_{(2,1)}|_0 & t_1t_2^2\sigma_{(1,2)}|_0 & t_2^3\sigma_{(0,3)}|_0\end{array}\right]\\[8pt]
            k' = 1 & 0 & \left[\def\arraystretch{1.3}\begin{array}{c|c}d_4 & 0\\
            \hline
            0 & d_4\end{array}\right] & \left[\def\arraystretch{1.3}\begin{array}{c|c|c}t_1\sigma_{(1,0)}|_2 & t_2\sigma_{(0,1)}|_2 & 0\\
            \hline
            0 & t_1\sigma_{(1,0)}|_2 & t_2\sigma_{(0,1)}|_2\end{array}\right] & \left[\def\arraystretch{1.3}\begin{array}{c|c|c|c}t_1^2\sigma_{(2,0)}|_0 & t_1t_2\sigma_{(1,1)}|_0 & t_2^2\sigma_{(0,2)}|_0 & 0\\
            \hline
            0 & t_1^2\sigma_{(2,0)}|_0 & t_1t_2\sigma_{(1,1)}|_0 & t_2^2\sigma_{(0,2)}|_0\end{array}\right]\\[12pt]
            k' = 2 & 0 & 0 & \left[\def\arraystretch{1.3}\begin{array}{c|c|c}d_2 & 0 & 0\\
            \hline
            0 & d_2 & 0\\
            \hline
            0 & 0 & d_2\end{array}\right] & \left[\def\arraystretch{1.3}\begin{array}{c|c|c|c}t_1\sigma_{(1,0)}|_0 & t_2\sigma_{(0,1)}|_0 & 0 & 0\\
            \hline
            0 & t_1\sigma_{(1,0)}|_0 & t_2\sigma_{(0,1)}|_0 & 0\\
            \hline
            0 & 0 & t_1\sigma_{(1,0)}|_0 & t_2\sigma_{(0,1)}|_0\end{array}\right]
        \end{array}
    $}
    \caption{The planes at levels $0, 1, 2, 3$ for $c = 2$, and the matrix of $D|_6$. The block in the column of $y^{\boldsymbol{u}}$ and the row of $y^{\boldsymbol{u'}}$ is the summand $t^{\boldsymbol{v}} \sigma_{\boldsymbol{v}}$ of $D$ with $\boldsymbol{v} = \boldsymbol{u} - \boldsymbol{u'}$.}
    \label{fig:Dsix}
\end{figure}

At $n = 6$, with source levels $k = 0, 1, 2, 3$ and target levels $k' = 0, 1, 2$, we illustrate the block block matrix, $D|_6$, as a table in Figure~\ref{fig:Dsix}. Although we normally omit the $t_i$, $\boldsymbol{v}$, and homological degree for these blocks, we include them in this case to demonstrate the function of the Kronecker product $T_{\boldsymbol{v}}^{(m)} \otimes \sigma_{\boldsymbol{v}}|_q$

For a truncated system, the block diagonal of Figure~\ref{fig:Dsix}, the blocks $\mathbf{1}_{N_k} \otimes d$ from each level to itself, is $\partial_0 = \id \otimes d$, and the first block superdiagonal, the blocks from level $k$ to level $k - 1$, is $\dcon = \partial_1 + \dots + \partial_c$. Each $\partial_j = t_j \otimes \sigma^{(j)}$ is the band $T^{(k)}_{\boldsymbol{e}_j} \otimes \sigma^{(j)}$ of these blocks. The blocks further above the diagonal carry the $\sigma_{\boldsymbol{v}}$ with $|\boldsymbol{v}| \ge 2$, which vanish for a truncated system. So only the diagonal blocks and the connecting blocks survive,
\[
    S_k = \sum_{j=1}^{c} T^{(k)}_{\boldsymbol{e}_j} \otimes \sigma^{(j)}, \qquad
    S_2 = \begin{bmatrix} \sigma^{(1)} & \sigma^{(2)} & 0\\ 0 & \sigma^{(1)} & \sigma^{(2)}\end{bmatrix} \quad (c = 2).
\]
\begin{remark}\label{rem:rank}
    Counting the copies on each hyperplane,
    \[
        \operatorname{rank} \Sh(\F, \sigma)_n = \sum_{2k \le n} N_k \operatorname{rank} F_{n-2k}.
    \]
\end{remark}

\begin{remark}[The square of $D$]\label{rem:hmf}
    Over $R$ the square of $D$ is not zero. With $\Theta = \sum_j t_j \otimes a_j \id$, the identity $D^2 = \partial_0 \dcon + \dcon \partial_0 = \Theta$ of Remark~\ref{rem:piecessquare} gives, for every $n$,
    \[
        D|_{n}\, D|_{n+1} = \Theta|_{n+1}, \qquad \bigl(\Theta|_{n+1}\bigr)_{k,k-1} = \sum_{j=1}^{c} a_j T^{(k)}_{\boldsymbol{e}_j} \otimes \id ,
    \]
    and every other block of $\Theta|_{n+1}$ is zero. So the product of two consecutive maps $D|_n$ is an explicit block matrix, with $a_j$ in the positions where $S_k$ has $\sigma^{(j)}$, one band for each color. This generalizes the identity $AB = BA = f \id$ of a matrix factorization to several elements $a_1, \dots, a_c$ at once.
\end{remark}

\section{Truncated systems for DG resolutions}\label{sec:sobieska}

Truncation is a restriction on the system. On any DG resolution, and in particular on the Taylor resolution of a monomial ideal, a truncated system always exists.

\begin{definition}\label{def:dg} \cite[Section 1.3]{AvramovSixLectures}
Let $A$ be a commutative ring. A \emph{DG algebra}, or differential graded algebra, over $A$ is a complex $(\F,d)$ of finitely generated free $A$-modules with unity $1 \in F_0$ and associative multiplication that respects the homological grading and satisfies the following:
    \begin{enumerate}
\item $d(ab)=d(a)b+(-1)^{\hdeg(a)}ad(b)$, for elements $a,b \in \F$, each in a single homological degree, 
\item $ab=(-1)^{\hdeg(a)\hdeg(b)}ba$, for $a,b \in \F$, and
\item $a^2=0$ for $\hdeg(a)$ odd.
\end{enumerate}
When $\F$ is a free resolution which is also a DG algebra, we call $\F$ a \emph{DG-resolution}.
\end{definition}

The Taylor resolution is a DG algebra \cite{Gemeda,AvramovSixLectures}, with multiplication
\[
    \eps_S \, \eps_T = (-1)^{\#\{(s,t) \in S \times T \, : \, s < t\}} \, \frac{m_S m_T}{m_{S \cup T}} \, \eps_{S \cup T} \quad \text{if } S \cap T = \varnothing, \qquad \eps_S \, \eps_T = 0 \quad \text{otherwise.}
\] The following proposition is a standard result, and we include a short proof for completeness.
\begin{proposition}\label{prop:dgtrunc}
    Let $\F$ be a DG-resolution of $R/I$, and let $a_1, \dots, a_c \in I$. Since $a_i \in I = d(F_1)$, we may choose $z_i \in F_1$ with $d z_i = a_i$. Let $\sigma^{(i)}$ be multiplication by $z_i$, and $\sigma_{\boldsymbol{u}} = 0$ for $|\boldsymbol{u}| \ge 2$. Then $\sigma$ is a truncated system of higher homotopies for $a_1, \dots, a_c$ on $\F$.
\end{proposition}

\begin{proof}
    Since $d(z_i x) = a_i x - z_i\,d(x)$, we have $d\sigma^{(i)} + \sigma^{(i)}d = a_i \id$. Since $z_i^2 = 0$ and $z_i z_k = -z_k z_i$, we have $(\sigma^{(i)})^2 = 0$ and $\sigma^{(i)}\sigma^{(k)} + \sigma^{(k)}\sigma^{(i)} = 0$. These are the identities (a)-(c) of Definition~\ref{def:higher}.
\end{proof}

\begin{remark}[Grading]\label{rem:multigrading}
    In Proposition~\ref{prop:dgtrunc}, suppose that $\F$ is graded, in the standard grading or in the fine multigrading, and that each $z_i$ is homogeneous. Then $\Sh(\F, \sigma)$ is graded, with $y_i$ in the degree of $z_i$, which is the degree of $a_i$. In the fine multigrading the $z_i$ can be chosen homogeneous exactly when every $a_i$ is a scalar times a monomial, and Example~\ref{ex:nonmonomial} shows what fails otherwise.
\end{remark}

\begin{example}[The running example]\label{ex:runtaylortruncated} Take $R = \kk[x, y, z, w]$, the ideal $I = (xy, yz, zw, xw)$ with generators $m_1, m_2, m_3, m_4$ in that order, and $\F$ the Taylor resolution of $R/I$.  
In Figure~\ref{fig:runTaylortruncated} the left figure represents $G(\F)$ for the Taylor resolution and the right figure is an illustration of a truncated system for $a = xy = m_1$ with 
\[\begin{aligned}
\sigma(\eps_\varnothing) &= \eps_1, & \sigma(\eps_2) &= -y \eps_{12}, & \sigma(\eps_3) &= -\eps_{13}, & \sigma(\eps_4) &= -x \eps_{14},\\
\sigma(\eps_{23}) &= y \eps_{123}, & \sigma(\eps_{24}) &= xy \eps_{124}, & \sigma(\eps_{34}) &= x \eps_{134}, & \sigma(\eps_{234}) &= -xy \eps_{1234},
\end{aligned}
    \]
    representing the dashed arrows. Note that our $\sigma$ has unit entries because we chose $a \notin \mathfrak{m} I$, so the condition from Remark \ref{rem:mincond} is not satisfied.

\begin{figure}[htbp]
    \centering
    \resizebox{\textwidth}{!}{%
  \begin{tikzpicture}[>=Stealth, x=1cm, y=1cm,
        vtx/.style={inner sep=1.2pt, font=\scriptsize, draw, rounded corners=1pt, fill=white, minimum width=5mm, minimum height=3.5mm},
        uedge/.style={->, gray!60},
        dedge/.style={->, gray!60},
        wt/.style={font=\tiny, fill=white, inner sep=0.4pt, text=black},
        mat/.style={font=\scriptsize, anchor=west}]
        \node[vtx, align=center] (e0) at (0.00,0.00) {$\eps_\varnothing$\\[-1.5pt]{\tiny $1$}};
        \node[vtx, align=center] (e1) at (-2.17,1.90) {$\eps_{1}$\\[-1.5pt]{\tiny $xy$}};
        \node[vtx, align=center] (e2) at (-0.72,1.90) {$\eps_{2}$\\[-1.5pt]{\tiny $yz$}};
        \node[vtx, align=center] (e3) at (0.72,1.90) {$\eps_{3}$\\[-1.5pt]{\tiny $zw$}};
        \node[vtx, align=center] (e4) at (2.17,1.90) {$\eps_{4}$\\[-1.5pt]{\tiny $xw$}};
        \node[vtx, align=center] (e12) at (-3.62,3.80) {$\eps_{12}$\\[-1.5pt]{\tiny $xyz$}};
        \node[vtx, align=center] (e13) at (-2.17,3.80) {$\eps_{13}$\\[-1.5pt]{\tiny $xyzw$}};
        \node[vtx, align=center] (e14) at (-0.72,3.80) {$\eps_{14}$\\[-1.5pt]{\tiny $xyw$}};
        \node[vtx, align=center] (e23) at (0.72,3.80) {$\eps_{23}$\\[-1.5pt]{\tiny $yzw$}};
        \node[vtx, align=center] (e24) at (2.17,3.80) {$\eps_{24}$\\[-1.5pt]{\tiny $xyzw$}};
        \node[vtx, align=center] (e34) at (3.62,3.80) {$\eps_{34}$\\[-1.5pt]{\tiny $xzw$}};
        \node[vtx, align=center] (e123) at (-2.17,5.70) {$\eps_{123}$\\[-1.5pt]{\tiny $xyzw$}};
        \node[vtx, align=center] (e124) at (-0.72,5.70) {$\eps_{124}$\\[-1.5pt]{\tiny $xyzw$}};
        \node[vtx, align=center] (e134) at (0.72,5.70) {$\eps_{134}$\\[-1.5pt]{\tiny $xyzw$}};
        \node[vtx, align=center] (e234) at (2.17,5.70) {$\eps_{234}$\\[-1.5pt]{\tiny $xyzw$}};
        \node[vtx, align=center] (e1234) at (0.00,7.60) {$\eps_{1234}$\\[-1.5pt]{\tiny $xyzw$}};
        \draw[dedge] (e1) -- node[wt, pos=0.46] {$xy$} (e0);
        \draw[dedge] (e2) -- node[wt, pos=0.62] {$yz$} (e0);
        \draw[dedge] (e3) -- node[wt, pos=0.38] {$zw$} (e0);
        \draw[dedge] (e4) -- node[wt, pos=0.38] {$xw$} (e0);
        \draw[dedge] (e12) -- node[wt, pos=0.46] {$-x$} (e2);
        \draw[dedge] (e12) -- node[wt, pos=0.38] {$z$} (e1);
        \draw[dedge] (e13) -- node[wt, pos=0.46] {$-xy$} (e3);
        \draw[dedge] (e13) -- node[wt, pos=0.7] {$zw$} (e1);
        \draw[dedge] (e14) -- node[wt, pos=0.46] {$-y$} (e4);
        \draw[dedge] (e14) -- node[wt, pos=0.54] {$w$} (e1);
        \draw[dedge] (e23) -- node[wt, pos=0.7] {$-y$} (e3);
        \draw[dedge] (e23) -- node[wt, pos=0.54] {$w$} (e2);
        \draw[dedge] (e24) -- node[wt, pos=0.46] {$-yz$} (e4);
        \draw[dedge] (e24) -- node[wt, pos=0.3] {$xw$} (e2);
        \draw[dedge] (e34) -- node[wt, pos=0.38] {$-z$} (e4);
        \draw[dedge] (e34) -- node[wt, pos=0.7] {$x$} (e3);
        \draw[dedge] (e123) -- node[wt, pos=0.3] {$x$} (e23);
        \draw[uedge] (e123) -- node[wt, pos=0.46] {$-1$} (e13);
        \draw[dedge] (e123) -- node[wt, pos=0.62] {$w$} (e12);
        \draw[uedge] (e124) -- node[wt, pos=0.46] {$1$} (e24);
        \draw[dedge] (e124) -- node[wt, pos=0.7] {$-z$} (e14);
        \draw[dedge] (e124) -- node[wt, pos=0.62] {$w$} (e12);
        \draw[dedge] (e134) -- node[wt, pos=0.7] {$y$} (e34);
        \draw[dedge] (e134) -- node[wt, pos=0.3] {$-z$} (e14);
        \draw[uedge] (e134) -- node[wt, pos=0.46] {$1$} (e13);
        \draw[dedge] (e234) -- node[wt, pos=0.46] {$y$} (e34);
        \draw[uedge] (e234) -- node[wt, pos=0.46] {$-1$} (e24);
        \draw[dedge] (e234) -- node[wt, pos=0.54] {$x$} (e23);
        \draw[uedge] (e1234) -- node[wt, pos=0.46] {$-1$} (e234);
        \draw[uedge] (e1234) -- node[wt, pos=0.38] {$1$} (e134);
        \draw[uedge] (e1234) -- node[wt, pos=0.62] {$-1$} (e124);
        \draw[uedge] (e1234) -- node[wt, pos=0.38] {$1$} (e123);
        
    \end{tikzpicture}
    \quad
  \begin{tikzpicture}[>=Stealth, x=1cm, y=1cm,
        vtx/.style={inner sep=1.2pt, font=\scriptsize, draw, rounded corners=1pt, fill=white, minimum width=5mm, minimum height=3.5mm},
        dedge/.style={->, gray!60},
        sedge/.style={->, green!55!black, thick, densely dashed},
        swt/.style={font=\tiny, fill=white, inner sep=0.4pt, text=green!40!black}]
        \node[vtx, align=center] (e0) at (0.00,0.00) {$\eps_\varnothing$\\[-1.5pt]{\tiny $1$}};
        \node[vtx, align=center] (e1) at (-2.17,1.90) {$\eps_{1}$\\[-1.5pt]{\tiny $xy$}};
        \node[vtx, align=center] (e2) at (-0.72,1.90) {$\eps_{2}$\\[-1.5pt]{\tiny $yz$}};
        \node[vtx, align=center] (e3) at (0.72,1.90) {$\eps_{3}$\\[-1.5pt]{\tiny $zw$}};
        \node[vtx, align=center] (e4) at (2.17,1.90) {$\eps_{4}$\\[-1.5pt]{\tiny $xw$}};
        \node[vtx, align=center] (e12) at (-3.62,3.80) {$\eps_{12}$\\[-1.5pt]{\tiny $xyz$}};
        \node[vtx, align=center] (e13) at (-2.17,3.80) {$\eps_{13}$\\[-1.5pt]{\tiny $xyzw$}};
        \node[vtx, align=center] (e14) at (-0.72,3.80) {$\eps_{14}$\\[-1.5pt]{\tiny $xyw$}};
        \node[vtx, align=center] (e23) at (0.72,3.80) {$\eps_{23}$\\[-1.5pt]{\tiny $yzw$}};
        \node[vtx, align=center] (e24) at (2.17,3.80) {$\eps_{24}$\\[-1.5pt]{\tiny $xyzw$}};
        \node[vtx, align=center] (e34) at (3.62,3.80) {$\eps_{34}$\\[-1.5pt]{\tiny $xzw$}};
        \node[vtx, align=center] (e123) at (-2.17,5.70) {$\eps_{123}$\\[-1.5pt]{\tiny $xyzw$}};
        \node[vtx, align=center] (e124) at (-0.72,5.70) {$\eps_{124}$\\[-1.5pt]{\tiny $xyzw$}};
        \node[vtx, align=center] (e134) at (0.72,5.70) {$\eps_{134}$\\[-1.5pt]{\tiny $xyzw$}};
        \node[vtx, align=center] (e234) at (2.17,5.70) {$\eps_{234}$\\[-1.5pt]{\tiny $xyzw$}};
        \node[vtx, align=center] (e1234) at (0.00,7.60) {$\eps_{1234}$\\[-1.5pt]{\tiny $xyzw$}};
        \draw[dedge] (e1) -- (e0);
        \draw[dedge] (e2) -- (e0);
        \draw[dedge] (e3) -- (e0);
        \draw[dedge] (e4) -- (e0);
        \draw[dedge] (e12) -- (e2);
        \draw[dedge] (e12) -- (e1);
        \draw[dedge] (e13) -- (e3);
        \draw[dedge] (e13) -- (e1);
        \draw[dedge] (e14) -- (e4);
        \draw[dedge] (e14) -- (e1);
        \draw[dedge] (e23) -- (e3);
        \draw[dedge] (e23) -- (e2);
        \draw[dedge] (e24) -- (e4);
        \draw[dedge] (e24) -- (e2);
        \draw[dedge] (e34) -- (e4);
        \draw[dedge] (e34) -- (e3);
        \draw[dedge] (e123) -- (e23);
        \draw[dedge] (e123) -- (e13);
        \draw[dedge] (e123) -- (e12);
        \draw[dedge] (e124) -- (e24);
        \draw[dedge] (e124) -- (e14);
        \draw[dedge] (e124) -- (e12);
        \draw[dedge] (e134) -- (e34);
        \draw[dedge] (e134) -- (e14);
        \draw[dedge] (e134) -- (e13);
        \draw[dedge] (e234) -- (e34);
        \draw[dedge] (e234) -- (e24);
        \draw[dedge] (e234) -- (e23);
        \draw[dedge] (e1234) -- (e234);
        \draw[dedge] (e1234) -- (e134);
        \draw[dedge] (e1234) -- (e124);
        \draw[dedge] (e1234) -- (e123);
        \draw[sedge] (e0) to[bend left=18] node[swt] {$1$} (e1);
        \draw[sedge] (e2) to[bend left=18] node[swt] {$-y$} (e12);
        \draw[sedge] (e3) to[bend left=18] node[swt] {$-1$} (e13);
        \draw[sedge] (e4) to[bend left=18] node[swt] {$-x$} (e14);
        \draw[sedge] (e23) to[bend left=18] node[swt] {$y$} (e123);
        \draw[sedge] (e24) to[bend left=18] node[swt] {$xy$} (e124);
        \draw[sedge] (e34) to[bend left=18] node[swt] {$x$} (e134);
        \draw[sedge] (e234) to[bend left=18] node[swt] {$-xy$} (e1234);
    \end{tikzpicture}}
    \caption{The running example: Taylor resolution and example truncated system}
    \label{fig:runTaylortruncated}
\end{figure}
\end{example}

On the Taylor resolution the homotopies can be written down, and the system they form is truncated. Sobieska gives the maps in \cite[Def.~3.1]{Sobieska} and the identities they satisfy in \cite[Thm.~3.2]{Sobieska}, for any regular sequence in $I$ and any choice of expressions for its members.

\begin{theorem}[{\cite[Thm.~3.2]{Sobieska}}]\label{thm:sobieska}
    Let $I = (m_1, \dots, m_r)$ be a monomial ideal, let $\F$ be its Taylor resolution, and let $a_1, \dots, a_c$ be an $R$-regular sequence in $I$, not necessarily homogeneous. Fix expressions $a_j = \sum_{t=1}^{r} f_{j,t} m_t$ with $f_{j,t} \in R$. Write $S = \{s_1 < \dots < s_k\}$ and let $p_{t,S}$ be the position of $t$ in $S \cup t$, so that $S \cup t = \{s_1, \dots, s_{p_{t,S}-1}, t, s_{p_{t,S}}, \dots, s_k\}$. Then
    \[
        \sigma^{(j)}(\eps_S) = \sum_{t \notin S} (-1)^{k - p_{t,S} - 1} f_{j,t} \frac{m_t m_S}{m_{S \cup t}} \eps_{S \cup t}, \qquad \sigma_{\boldsymbol{u}} = 0 \text{ for } |\boldsymbol{u}| \ge 2,
    \]
    is a system of higher homotopies for $a_1, \dots, a_c$ on $\F$. So every regular sequence in $I$ has a truncated system on the Taylor resolution.
\end{theorem}

Since $m_{S \cup t} = \lcm(m_t, m_S)$ and $p_{t,S} - 1 = \#\{s \in S : s < t\}$, the maps read
\[
    \sigma^{(j)}(\eps_S) = \sum_{t \notin S} (-1)^{\#\{s \in S : s > t\}} f_{j,t} \gcd(m_t, m_S) \eps_{S \cup t},
\]
with the sign convention of the Taylor differential. These are the maps of Proposition~\ref{prop:dgtrunc} for $z_j = \sum_{t=1}^{r} f_{j,t} \eps_t$.

\section{The Shamash complex as an input for Morse theory}\label{sec:hyp}

This section applies Theorem~\ref{thm:amt} to the Shamash resolution, lifts a matching on $\F$ to it, and defines the transferred maps $\sigma^{\M}$. Throughout, $\overline r$ is the image in $\Rbar$ of $r \in R$.

The Shamash resolution $\Sh(\F, \sigma)$ is free over $\Rbar$ on the cells $y^{\boldsymbol{u}} \otimes \bar\alpha$, for $\boldsymbol{u} \in \N^c$ and $\alpha \in \B$, the cell $y^{\boldsymbol{u}} \otimes \bar\alpha$ in homological degree $2|\boldsymbol{u}| + \hdeg(\alpha)$. On a cell
\[
    d_{\Sh}(y^{\boldsymbol{u}} \otimes \bar\alpha) = y^{\boldsymbol{u}} \otimes \bar d(\bar\alpha) + \sum_{j : u_j \ge 1} y^{\boldsymbol{u} - \boldsymbol{e}_j} \otimes \bar\sigma^{(j)}(\bar\alpha),
\]
with coefficients read modulo $\aaa$. In matrices, examples and figures we drop the bar and the tensor sign and write $y^{\boldsymbol{u}} \alpha$ for this cell, with $\alpha$ for $\boldsymbol{u} = \boldsymbol{0}$. So the complex graph of $\Sh(\F, \sigma)$ has two kinds of edges. A \emph{level edge} $y^{\boldsymbol{u}} \otimes \bar\alpha \to y^{\boldsymbol{u}} \otimes \bar\beta$, of weight $\overline{d_{\alpha,\beta}}$, is an edge of $\bar\partial_0$, and a \emph{color-$j$ connecting edge} $y^{\boldsymbol{u}} \otimes \bar\alpha \to y^{\boldsymbol{u} - \boldsymbol{e}_j} \otimes \bar\beta$, of weight $\overline{\sigma^{(j)}_{\alpha,\beta}}$, is an edge of $\bar\partial_j$. A \emph{connecting edge} is a color-$j$ connecting edge for some $j$. The level edges over a fixed $\boldsymbol{u}$ are copies of edges of $G(\F)$, and since $\kk \to \Rbar$ is injective they include a copy of every unit edge of $G(\F)$.

\begin{proposition}\label{thm:A}
    Let $\F$ be a free resolution of $R/I$ over $R$ with a finite basis of each $F_q$, let $a_1, \dots, a_c$ be an $R$-regular sequence in $I$, let $\sigma = \{\sigma^{(j)}\}_{j=1}^c$ be a truncated system of higher homotopies for $a_1, \dots, a_c$ on $\F$, and let $\A$ be a Morse matching on $\Sh(\F, \sigma)$. Then the Morse complex $\Sh(\F, \sigma)^{\A}$, free over $\Rbar$ on the cells critical for $\A$, is a free $\Rbar$-resolution of $R/I$, and Theorem~\ref{thm:amt} gives a contraction of $\Sh(\F, \sigma)$ onto it.
\end{proposition}

\begin{proof}
    By Theorem~\ref{thm:shamash}, $\Sh(\F, \sigma)$ is a free $\Rbar$-resolution of $R/I$, and Theorem~\ref{thm:amt} applies to it over $\Rbar$. A contraction is a homotopy equivalence, so $\Sh(\F, \sigma)^{\A}$ has the homology of $\Sh(\F, \sigma)$, which is $R/I$ in degree $0$.
\end{proof}

\subsection*{Lifting a matching}

Given a partial matching $\M$ on a free resolution $\F$ of $R/I$, we can also produce a matching on the Eisenbud--Shamash complex by lifting it to $\Sh(\F,\sigma)$.

\begin{definition}\label{def:lift}
    Let $\F$ be a free resolution of $R/I$ over $R$ with a finite basis of each $F_q$, let $a_1, \dots, a_c$ be an $R$-regular sequence in $I$, let $\sigma = \{\sigma^{(j)}\}_{j=1}^c$ be a truncated system of higher homotopies for $a_1, \dots, a_c$ on $\F$, and let $\M$ be a partial matching on $\F$. The \emph{lift} of $\M$ to $\Sh(\F, \sigma)$ is
    \[
        \lift{M} = \bigl\{ (y^{\boldsymbol{u}} \otimes \bar\alpha, y^{\boldsymbol{u}} \otimes \bar\beta) : (\alpha, \beta) \in \M, \boldsymbol{u} \in \N^c \bigr\},
    \]
    one copy of $\M$ over each $\boldsymbol{u}$. A cell $y^{\boldsymbol{u}} \otimes \bar\alpha$ is critical for $\lift{M}$ if and only if $\alpha$ is critical for $\M$, so the critical cells of $\lift{M}$ are the $y^{\boldsymbol{u}} \otimes \bar\kappa$ with $\kappa$ critical for $\M$.
\end{definition}

\begin{proposition}\label{prop:lift}
    In the setting of Definition~\ref{def:lift}, if $\M$ is a Morse matching on $\F$ then $\lift{M}$ is a Morse matching on $\Sh(\F, \sigma)$.
\end{proposition}

\begin{proof}
    Pairs of $\lift{M}$ over distinct $\boldsymbol{u}$ are disjoint, and over the same $\boldsymbol{u}$ they come from disjoint pairs of $\M$. Each lies over a matched edge of $G(\F)$, whose weight is a nonzero scalar and stays a unit in $\Rbar$. So $\lift{M}$ is a partial matching.

    Every edge of $\Grep(\Sh(\F, \sigma), \lift{M})$ is either a level or reversed edge, which preserves $\boldsymbol{u}$, or a connecting edge, which lowers $|\boldsymbol{u}|$ by one and is never reversed. Over a fixed $\boldsymbol{u}$, $y^{\boldsymbol{u}} \otimes \bar\alpha \mapsto \alpha$ carries edges injectively to edges of $\Grep(\F, \M)$, since the level edges over $\boldsymbol{u}$ are copies of edges of $G(\F)$. So a directed path from $y^{\boldsymbol{u}} \otimes \bar\alpha$ uses at most $|\boldsymbol{u}|$ connecting edges, and between them projects to directed paths of $\Grep(\F, \M)$, which are finite. Hence $\lift{M}$ is locally finite.
\end{proof}

\begin{definition}\label{def:transferred}
In the setting of Definition~\ref{def:lift}, let $\M$ be a Morse matching on $\F$, with Morse complex $\F^{\M} = (K_\bullet, d^{\M})$ and contraction $(i^{\M}, p^{\M}, h^{\M})$ of Theorem~\ref{thm:amt}. Set $\sigma^{\M}_0 = d^{\M}$ and, for $|\boldsymbol{u}| = k \ge 1$,
    \begin{equation}\label{eq:transfer}
        \sigma^{\M}_{\boldsymbol{u}}|_q = \sum_{\substack{(j_1, \dots, j_k) \\ \boldsymbol{e}_{j_1} + \dots + \boldsymbol{e}_{j_k} = \boldsymbol{u}}} p^{\M}|_{q+2k-1} \, \sigma^{(j_1)}|_{q+2k-2} \, h^{\M}|_{q+2k-3} \, \sigma^{(j_2)}|_{q+2k-4} \cdots h^{\M}|_{q+1} \, \sigma^{(j_k)}|_q \, i^{\M}|_q,
    \end{equation}
    for every $q$,
    that is,
    \[
        \sigma^{\M}_{\boldsymbol{u}} = \sum p^{\M} \sigma^{(j_1)} h^{\M} \sigma^{(j_2)} \cdots h^{\M} \sigma^{(j_k)} i^{\M}
    \]
    over the $\binom{k}{u_1, \dots, u_c}$ sequences in which the color $j$ occurs $u_j$ times.
    In particular $\sigma^{\M}_{\boldsymbol{e}_j} = p^{\M} \sigma^{(j)} i^{\M}$.
    \end{definition}

The formula \eqref{eq:transfer} can be read off a picture like Figure~\ref{fig:Gsigma1}, where the reversed matched edges of $\Grep(\F, \M)$ are red, its other edges are black (drawn in gray), and the arrows of the $\sigma^{(j)}$ are colored. Start at a critical cell $\kappa$. Follow a red-black path that stays in the homological degree of $\kappa$ and the one below it and ends in the degree of $\kappa$, and then a colored edge of color $j_k$. From the cell this reaches, follow a red-black path that starts up a red edge and then a new colored edge of color $j_{k-1}$, and continue in this way through all $k$ colored edges. After the last colored edge, follow a red-black path that eventually ends on a critical cell $\kappa'$. These three kinds of red-black paths are the terms of $i^{\M}$, $h^{\M}$ and $p^{\M}$.
The coefficient of a critical cell $\kappa'$ in $\sigma^{\M}_{\boldsymbol{u}}(\kappa)$ is the sum, over all such paths from $\kappa$ to $\kappa'$, of the product of the weights of the edges on the path.

  Each summand raises homological degree by $2k - 1$, the $k$ factors $\sigma^{(j_r)}$ and the $k - 1$ factors $h^{\M}$ each raising by one homological degree while $p^{\M}$ and $i^{\M}$ preserve the homological degree. So $\sigma^{\M}_{\boldsymbol{u}}$ raises homological degree by $2|\boldsymbol{u}| - 1$.

\section{Motivating the Main Theorem} \label{sec:examples}

The examples in this section continue Example~\ref{ex:runSec3} with Sobieska's homotopies from Section \ref{sec:sobieska} to motivate and illustrate the Main Theorem \ref{thm:commute} which we present in Section \ref{sec:commute}.

\begin{example}[The running example, $c = 1$]\label{ex:running}
    Take the running example, Example \ref{ex:runSec3}
    with $c = 1$, $a_1 = m_1 = xy$ and $\sigma$ Sobieska's homotopy (Theorem~\ref{thm:sobieska}).

    The nonzero values of $\sigma$ on cells are
    \[
        \begin{aligned}
            \sigma(\eps_\varnothing) &= \eps_1, & \sigma(\eps_2) &= -y \eps_{12}, & \sigma(\eps_3) &= -\eps_{13}, & \sigma(\eps_4) &= -x \eps_{14},\\
            \sigma(\eps_{23}) &= y \eps_{123}, & \sigma(\eps_{24}) &= xy \eps_{124}, & \sigma(\eps_{34}) &= x \eps_{134}, & \sigma(\eps_{234}) &= -xy \eps_{1234},
        \end{aligned}
    \]
    the dashed arrows of Figure~\ref{fig:Gsigma1}.

    \begin{figure}[htbp]
        \centering
        \resizebox{0.50\textwidth}{!}{%
        \begin{tikzpicture}[>=Stealth, x=1.32cm, y=1.5cm,
            vtx/.style={inner sep=1.6pt, font=\tiny, draw, rounded corners=1pt, fill=white, minimum width=6mm, minimum height=4mm},
            crit/.style={vtx, fill=gray!10},
            mu/.style={vtx, fill=red!12},
            md/.style={vtx, fill=orange!18},
            medge/.style={->, red!70!black, thick},
            sedge/.style={->, green!55!black, thick, densely dashed},
            dedge/.style={->, gray!50, thin}]
            \node[crit] (eempty) at (0, 0) {$\eps_\varnothing$};
            \node[crit] (e1) at (-1.5, 1) {$\eps_1$};
            \node[crit] (e2) at (1.5, 1) {$\eps_4$};
            \node[crit] (e3) at (-0.5, 1) {$\eps_2$};
            \node[crit] (e4) at (0.5, 1) {$\eps_3$};
            \node[crit] (e12) at (-0.5, 2) {$\eps_{14}$};
            \node[crit] (e13) at (-2.5, 2) {$\eps_{12}$};
            \node[md] (e14) at (-1.5, 2) {$\eps_{13}$};
            \node[md] (e23) at (1.5, 2) {$\eps_{24}$};
            \node[crit] (e24) at (2.5, 2) {$\eps_{34}$};
            \node[crit] (e34) at (0.5, 2) {$\eps_{23}$};
            \node[mu] (e123) at (-0.5, 3) {$\eps_{124}$};
            \node[crit] (e124) at (0.5, 3) {$\eps_{134}$};
            \node[mu] (e134) at (-1.5, 3) {$\eps_{123}$};
            \node[md] (e234) at (1.5, 3) {$\eps_{234}$};
            \node[mu] (e1234) at (-1, 4) {$\eps_{1234}$};
            \draw[dedge] (e1) -- (eempty);
            \draw[dedge] (e3) -- (eempty);
            \draw[dedge] (e4) -- (eempty);
            \draw[dedge] (e2) -- (eempty);
            \draw[dedge] (e13) -- (e3);
            \draw[dedge] (e13) -- (e1);
            \draw[dedge] (e14) -- (e4);
            \draw[dedge] (e14) -- (e1);
            \draw[dedge] (e12) -- (e2);
            \draw[dedge] (e12) -- (e1);
            \draw[dedge] (e34) -- (e4);
            \draw[dedge] (e34) -- (e3);
            \draw[dedge] (e23) -- (e2);
            \draw[dedge] (e23) -- (e3);
            \draw[dedge] (e24) -- (e2);
            \draw[dedge] (e24) -- (e4);
            \draw[dedge] (e134) -- (e34);
            \draw[medge] (e14) -- (e134);
            \draw[dedge] (e134) -- (e13);
            \draw[medge] (e23) -- (e123);
            \draw[dedge] (e123) -- (e12);
            \draw[dedge] (e123) -- (e13);
            \draw[dedge] (e124) -- (e24);
            \draw[dedge] (e124) -- (e12);
            \draw[dedge] (e124) -- (e14);
            \draw[dedge] (e234) -- (e24);
            \draw[dedge] (e234) -- (e23);
            \draw[dedge] (e234) -- (e34);
            \draw[medge] (e234) -- (e1234);
            \draw[dedge] (e1234) -- (e124);
            \draw[dedge] (e1234) -- (e123);
            \draw[dedge] (e1234) -- (e134);
            \draw[sedge] (eempty) to[bend left=18] (e1);
            \draw[sedge] (e3) to[bend left=18] (e13);
            \draw[sedge] (e4) to[bend left=18] (e14);
            \draw[sedge] (e2) to[bend left=18] (e12);
            \draw[sedge] (e34) to[bend left=18] (e134);
            \draw[sedge] (e23) to[bend left=18] (e123);
            \draw[sedge] (e24) to[bend left=18] (e124);
            \draw[sedge] (e234) to[bend left=18] (e1234);
        \end{tikzpicture}}
        \caption{$G(\F)$ with the arrows of $\sigma$, $c = 1$.}
        \label{fig:Gsigma1}
\end{figure}

    In the bases
    of Section~\ref{sec:amt},
    \[
        \sigma|_0 = \begin{bmatrix} 1\\ 0\\ 0\\ 0 \end{bmatrix}, \quad
        \sigma|_1 = \begin{bmatrix} 0 & -y & 0 & 0\\ 0 & 0 & -1 & 0\\ 0 & 0 & 0 & -x\\ 0 & 0 & 0 & 0\\ 0 & 0 & 0 & 0\\ 0 & 0 & 0 & 0 \end{bmatrix},
    \]
    \[
        \sigma|_2 = \begin{bmatrix} 0 & 0 & 0 & y & 0 & 0\\ 0 & 0 & 0 & 0 & xy & 0\\ 0 & 0 & 0 & 0 & 0 & x\\ 0 & 0 & 0 & 0 & 0 & 0 \end{bmatrix}, \quad
        \sigma|_3 = \begin{bmatrix} 0 & 0 & 0 & -xy \end{bmatrix},
    \]
    with $\sigma|_4 = 0$.

Over $R/(xy)$, $\Sh(\F, \sigma)$ has ranks $1, 4, 7, 8, 8, \dots$. With cells ordered by level, $d_{\Sh}$ has the block form of Figure~\ref{fig:Sham1}, with first matrices
\[d_{\Sh}|_1 = \overline{d_1} = \ \ \begin{NiceArray}{ccccc}
     \eps_1 & \eps_2 & \eps_3 & \eps_4 &\\
      0 & yz & zw & xw & \eps_{\varnothing}\\
      &&&&\\
      \CodeAfter
  \SubMatrix[{2-1}{2-4}]
\end{NiceArray},\]

    \[
        d_{\Sh}|_2 = \left[\begin{array}{c|c} \overline{d_2} & \overline{\sigma|_0} \end{array}\right] =\ \ \begin{NiceArray}{cccccc|cc}
        \eps_{12}  & \eps_{13}  & \eps_{14}  & \eps_{23}  & \eps_{24}  & \eps_{34}  & y_1\eps_{\varnothing} &\\
        z & zw & w & 0 & 0 & 0 & 1 & \ \eps_1 \\
        -x & 0 & 0 & w & xw & 0 & 0 & \ \eps_2\\
        0 & 0 & 0 & -y & 0 & x & 0 &\ \eps_3\\
        0 & 0 & -y & 0 & -yz & -z & 0 & \ \eps_4,
        \CodeAfter
  \SubMatrix[{2-1}{5-7}]
\end{NiceArray}
    \]

    \[
        d_{\Sh}|_3 = \left[\begingroup
    \renewcommand{\arraystretch}{1.5}\begin{array}{c|c} \overline{d_3} & \overline{\sigma|_1}\\
            \hline
            0 & \overline{d_1} \end{array}\endgroup\right] = \ \ \begin{NiceArray}{cccc|ccccc}
        \eps_{123}  & \eps_{124}  & \eps_{134}  & \eps_{234}  & y_1\eps_{1}  & y_1\eps_{2}  & y_1\eps_3 & y_1\eps_4 &\\
        w & w & 0 & 0 & 0 & -y & 0 & 0 & \eps_{12}\\
        -1 & 0 & 1 & 0 & 0 & 0 & -1 & 0 & \eps_{13} \\
        0 & -z & -z & 0 & 0 & 0 & 0 & -x & \eps_{14}\\
        x & 0 & 0 & x & 0 & 0 & 0 & 0 & \eps_{23}\\
        0 & 1 & 0 & -1 & 0 & 0 & 0 & 0 & \eps_{24} \\
        0 & 0 & y & y & 0 & 0 & 0 & 0 & \eps_{34}\\
        \hline
        0 & 0 & 0 & 0 & 0 & yz & zw & xw & y_1\eps_{\varnothing},
        \CodeAfter
  \SubMatrix[{2-1}{8-8}]
\end{NiceArray}
    \]

\[
        d_{\Sh}|_4 = \left[\begingroup
        \renewcommand{\arraystretch}{1.5}\begin{array}{c|c|c} \overline{d_4} & \overline{\sigma|_2}& 0\\
            \hline
            0 & \overline{d_2} & \overline{\sigma|_0}\end{array}\endgroup\right]=\ \ \begin{NiceArray}{c|cccccc|cc}
        \eps_{1234}  & y_1\eps_{12}  & y_1\eps_{13}  & y_1\eps_{14}  & y_1\eps_{23}  & y_1\eps_{24}  & y_1\eps_{34} & y_1^{(2)}\eps_{\varnothing} &\\
         1 & 0 & 0 & 0 & y & 0 & 0 & 0 & \eps_{123}\\
            -1 & 0 & 0 & 0 & 0 & 0 & 0 & 0 & \eps_{124}\\
            1 & 0 & 0 & 0 & 0 & 0 & x & 0 & \eps_{134}\\
            -1 & 0 & 0 & 0 & 0 & 0 & 0 & 0 & \eps_{234}\\
            \hline
            0& z & zw & w & 0 & 0 & 0 & 1 & y_1\eps_1 \\
            0 & -x & 0 & 0 & w & xw & 0 & 0 & y_1\eps_2 \\
            0 & 0 & 0 & 0 & -y & 0 & x & 0 & y_1\eps_3 \\
            0 & 0 & 0 & -y & 0 & -yz & -z & 0 & y_1\eps_4,
        \CodeAfter
  \SubMatrix[{2-1}{9-8}]
\end{NiceArray}
\]

\[
        d_{\Sh}|_5 = \left[\begingroup
        \renewcommand{\arraystretch}{1.5}\begin{array}{c|c} \overline{\sigma|_3} & 0\\
            \hline
             \overline{d_3} & \overline{\sigma|_1}\\
             \hline
             0 & \overline{d_1}\end{array}\endgroup\right]=\ \ \begin{NiceArray}{cccc|ccccc}
        y_1\eps_{123}  & y_1\eps_{124}  & y_1\eps_{134}  & y_1\eps_{234}  & y_1^{(2)}\eps_{1}  & y_1^{(2)}\eps_{2}  & y_1^{(2)}\eps_3 & y_1^{(2)}\eps_4 &\\
        0 & 0 & 0 & 0 & 0 & 0 & 0 & 0 & \eps_{1234}\\
        \hline
        w & w & 0 & 0 & 0 & -y & 0 & 0 & y_1\eps_{12}\\
        -1 & 0 & 1 & 0 & 0 & 0 & -1 & 0 & y_1\eps_{13} \\
        0 & -z & -z & 0 & 0 & 0 & 0 & -x & y_1\eps_{14}\\
        x & 0 & 0 & x & 0 & 0 & 0 & 0 & y_1\eps_{23}\\
        0 & 1 & 0 & -1 & 0 & 0 & 0 & 0 & y_1\eps_{24} \\
        0 & 0 & y & y & 0 & 0 & 0 & 0 & y_1\eps_{34}\\
        \hline
        0 & 0 & 0 & 0 & 0 & yz & zw & xw & y_1^{(2)}\eps_{\varnothing}.
        \CodeAfter
  \SubMatrix[{2-1}{9-8}]
\end{NiceArray}
\]

For $q \geq 2$,
$d_{\Sh}|_{2q}= d_{\Sh}|_{4}$ and $d_{\Sh}|_{2q+1} = d_{\Sh}|_{5}$. The bases for the columns and rows for
$d_{\Sh}|_{2q}$ are \begin{center}
columns: $y_1^{(q-2)}\eps_{1234}$, $y_1^{(q-1)}\eps_{12}$, $y_1^{(q-1)}\eps_{13}$, $y_1^{(q-1)}\eps_{14}$, $y_1^{(q-1)} \eps_{23}$, $y_1^{(q-1)} \eps_{24}$, $y_1^{(q-1)} \eps_{34}$, $y_1^{(q)} \eps_\varnothing$;\end{center}
\begin{center}rows: $y_1^{(q-2)}\eps_{123}$, $y_1^{(q-2)}\eps_{124}$, $y_1^{(q-2)}\eps_{134}$, $y_1^{(q-2)}\eps_{234}$, $y_1^{(q-1)}\eps_{1}$, $y_1^{(q-1)}\eps_{2}$, $y_1^{(q-1)}\eps_{3}$, $y_1^{(q-1)}\eps_{4}$.\end{center}
The bases for the columns and rows for
$d_{\Sh}|_{2q+1}$ are \begin{center}columns: $y_1^{(q-1)}\eps_{123}$, $y_1^{(q-1)}\eps_{124}$, $y_1^{(q-1)}\eps_{134}$, $y_1^{(q-1)}\eps_{234}$, $y_1^{(q)} \eps_1$, $y_1^{(q)} \eps_2$, $y_1^{(q)} \eps_3$, $y_1^{(q)} \eps_4$;\end{center}
\begin{center} rows: $y_1^{(q-2)}\eps_{1234}$, $y_1^{(q-1)}\eps_{12}$, $y_1^{(q-1)}\eps_{13}$, $y_1^{(q-1)}\eps_{14}$, $y_1^{(q-1)}\eps_{23}$, $y_1^{(q-1)}\eps_{24}$, $y_1^{(q-1)}\eps_{34}$, $y_1^{(q)} \eps_\varnothing$.
\end{center}
The edges in our matching $\M$ are $(\eps_{123}, \eps_{13})$, $(\eps_{124}, \eps_{24})$ and $(\eps_{1234}, \eps_{234})$. By the definition of $\lift{M}$, the edges of $\lift{M}$ are $$ \{(y_1^{(q)}\eps_{123}, y_1^{(q)}\eps_{13}), (y_1^{(q)}\eps_{124}, y_1^{(q)}\eps_{24}), (y_1^{(q)}\eps_{1234}, y_1^{(q)}\eps_{234})\}_{q=0}^{\infty}.$$

There are also unit connecting edges ${y_1^{(q+1)} \eps_\varnothing \xrightarrow{1} y_1^{(q)}\eps_1}$ and ${y_1^{(q+1)} \eps_3 \xrightarrow{-1} y_1^{(q)}\eps_{13}}$ for all $q \geq 0$. These are shown in Figure~\ref{fig:Grepc_16}, where dashed red arrows point toward partners in degree $7$. No lifted matching contains these edges.

To illustrate Theorem \ref{thm:commute}, note that the transferred system is given by the formula $\sigma^{\M} = p^{\M} \sigma i^{\M}$ of Definition \ref{def:transferred}, and the terms $\sigma^{\M}_k$ with $k \ge 2$ vanish since $p^{\M} \sigma (h^{\M} \sigma)^k i^{\M}=0$ for $k \geq 1$.
    Only $\sigma(\eps_3) = -\eps_{13}$ lands on a lower cell, so $\sigma^{\M} = p \sigma i$ except at $\eps_3$, where the paths ${\eps_{13} \xrightarrow{1} \eps_{123} \xrightarrow{w} \eps_{12}}$ and ${\eps_{13} \xrightarrow{1} \eps_{123} \xrightarrow{x} \eps_{23}}$ give
    \[
        \sigma^{\M}(\eps_3) = -p^{\M}(\eps_{13}) = -w \eps_{12} - x \eps_{23} ,
    \]
    so
    \[
        \sigma^{\M}|_0 = \begin{bmatrix} 1\\ 0\\ 0\\ 0 \end{bmatrix}, \qquad
        \sigma^{\M}|_1 = \begin{bmatrix} 0 & -y & -w & 0\\ 0 & 0 & 0 & -x\\ 0 & 0 & -x & 0\\ 0 & 0 & 0 & 0 \end{bmatrix}, \qquad
        \sigma^{\M}|_2 = \begin{bmatrix} 0 & 0 & 0 & x \end{bmatrix}.
    \]

    Note that the path $\eps_3 \to \eps_{13} \to \eps_{123} \to \eps_{23} \to \eps_{123}$ does not contribute to $\sigma^\M$ because a path described by $h^\M$ must start and end with a red arrow.
    
    \begin{figure}[htbp]
        \centering
        \resizebox{\textwidth}{!}{%
        \begin{tikzpicture}[>=Stealth, x=1.8cm, y=1.5cm,
            vtx/.style={inner sep=1.6pt, font=\tiny, draw, rounded corners=1pt, fill=white, minimum width=6mm, minimum height=4mm},
            crit/.style={vtx, fill=gray!10},
            mu/.style={vtx, fill=red!12},
            md/.style={vtx, fill=orange!18},
            medge/.style={->, red!70!black, thick},
            sedge/.style={->, green!55!black, thick},
            dedge/.style={->, gray!50, thin}]
            \fill[cyan!10, rounded corners=4pt] (-0.35,-0.40) rectangle (0.35,0.40);
            \fill[cyan!10] (-0.30,0.35) rectangle (0.30,0.65);
            \fill[cyan!10, rounded corners=4pt] (-1.85,0.60) rectangle (1.85,1.40);
            \fill[cyan!10] (-1.80,1.35) rectangle (1.80,1.65);
            \fill[cyan!10, rounded corners=4pt] (-2.85,1.60) rectangle (2.85,2.40);
            \fill[cyan!10] (-1.80,2.35) rectangle (1.80,2.65);
            \fill[cyan!10, rounded corners=4pt] (-1.85,2.60) rectangle (1.85,3.40);
            \fill[cyan!10] (-1.31,3.35) rectangle (-0.69,3.65);
            \fill[cyan!10, rounded corners=4pt] (-1.36,3.60) rectangle (-0.64,4.40);
            \fill[yellow!28, rounded corners=4pt] (3.65,1.60) rectangle (4.35,2.40);
            \fill[yellow!28] (3.70,2.35) rectangle (4.30,2.65);
            \fill[yellow!28, rounded corners=4pt] (2.65,2.60) rectangle (6.35,3.40);
            \fill[yellow!28] (2.70,3.35) rectangle (6.30,3.65);
            \fill[yellow!28, rounded corners=4pt] (0.64,3.60) rectangle (6.36,4.40);
            \fill[yellow!28] (0.69,4.35) rectangle (4.34,4.65);
            \fill[yellow!28, rounded corners=4pt] (0.61,4.60) rectangle (4.39,5.40);
            \fill[yellow!28] (1.14,5.35) rectangle (1.86,5.65);
            \fill[yellow!28, rounded corners=4pt] (1.09,5.60) rectangle (1.91,6.40);
            \fill[violet!10, rounded corners=4pt] (7.14,3.60) rectangle (7.86,4.40);
            \fill[violet!10] (7.19,4.35) rectangle (7.81,4.65);
            \fill[violet!10, rounded corners=4pt] (5.14,4.60) rectangle (8.86,5.40);
            \fill[violet!10] (5.19,5.35) rectangle (8.81,5.65);
            \fill[violet!10, rounded corners=4pt] (3.11,5.60) rectangle (8.89,6.40);
            \fill[lime!22, rounded corners=4pt] (9.64,5.60) rectangle (10.36,6.40);
            \node[font=\scriptsize, text=cyan!45!black, fill=white, fill opacity=0.85, text opacity=1, rounded corners=1pt, anchor=north, inner sep=1.2pt] at (0.00,-0.40) {$F_\bullet$};
            \node[font=\scriptsize, text=yellow!40!black, fill=white, fill opacity=0.85, text opacity=1, rounded corners=1pt, anchor=north, inner sep=1.2pt] at (4.00,1.60) {$y_1 F_\bullet$};
            \node[font=\scriptsize, text=violet!60!black, fill=white, fill opacity=0.85, text opacity=1, rounded corners=1pt, anchor=north, inner sep=1.2pt] at (7.50,3.60) {$y_1^{(2)} F_\bullet$};
            \node[font=\scriptsize, text=lime!40!black, fill=white, fill opacity=0.85, text opacity=1, rounded corners=1pt, anchor=north, inner sep=1.2pt] at (10.00,5.60) {$y_1^{(3)} F_\bullet$};
            \node[crit] (eempty) at (0, 0) {$\eps_\varnothing$};
            \node[crit] (e1) at (-1.5, 1) {$\eps_1$};
            \node[crit] (e2) at (1.5, 1) {$\eps_4$};
            \node[crit] (e3) at (-0.5, 1) {$\eps_2$};
            \node[crit] (e4) at (0.5, 1) {$\eps_3$};
            \node[crit] (e12) at (-0.5, 2) {$\eps_{14}$};
            \node[crit] (e13) at (-2.5, 2) {$\eps_{12}$};
            \node[md] (e14) at (-1.5, 2) {$\eps_{13}$};
            \node[md] (e23) at (1.5, 2) {$\eps_{24}$};
            \node[crit] (e24) at (2.5, 2) {$\eps_{34}$};
            \node[crit] (e34) at (0.5, 2) {$\eps_{23}$};
            \node[crit] (yempty) at (4, 2) {$y_1\eps_\varnothing$};
            \node[mu] (e123) at (-0.5, 3) {$\eps_{124}$};
            \node[crit] (e124) at (0.5, 3) {$\eps_{134}$};
            \node[mu] (e134) at (-1.5, 3) {$\eps_{123}$};
            \node[md] (e234) at (1.5, 3) {$\eps_{234}$};
            \node[crit] (y1) at (3, 3) {$y_1\eps_1$};
            \node[crit] (y2) at (6, 3) {$y_1\eps_4$};
            \node[crit] (y3) at (4, 3) {$y_1\eps_2$};
            \node[crit] (y4) at (5, 3) {$y_1\eps_3$};
            \node[mu] (e1234) at (-1, 4) {$\eps_{1234}$};
            \node[crit] (y12) at (3, 4) {$y_1\eps_{14}$};
            \node[crit] (y13) at (1, 4) {$y_1\eps_{12}$};
            \node[md] (y14) at (2, 4) {$y_1\eps_{13}$};
            \node[md] (y23) at (5, 4) {$y_1\eps_{24}$};
            \node[crit] (y24) at (6, 4) {$y_1\eps_{34}$};
            \node[crit] (y34) at (4, 4) {$y_1\eps_{23}$};
            \node[crit] (yyempty) at (7.5, 4) {$y_1^{(2)}\eps_\varnothing$};
            \node[mu] (y123) at (2, 5) {$y_1\eps_{124}$};
            \node[crit] (y124) at (3, 5) {$y_1\eps_{134}$};
            \node[mu] (y134) at (1, 5) {$y_1\eps_{123}$};
            \node[md] (y234) at (4, 5) {$y_1\eps_{234}$};
            \node[crit] (yy1) at (5.5, 5) {$y_1^{(2)}\eps_1$};
            \node[crit] (yy2) at (8.5, 5) {$y_1^{(2)}\eps_4$};
            \node[crit] (yy3) at (6.5, 5) {$y_1^{(2)}\eps_2$};
            \node[crit] (yy4) at (7.5, 5) {$y_1^{(2)}\eps_3$};
            \node[mu] (y1234) at (1.5, 6) {$y_1\eps_{1234}$};
            \node[crit] (yy12) at (5.5, 6) {$y_1^{(2)}\eps_{14}$};
            \node[crit] (yy13) at (3.5, 6) {$y_1^{(2)}\eps_{12}$};
            \node[md] (yy14) at (4.5, 6) {$y_1^{(2)}\eps_{13}$};
            \node[md] (yy23) at (7.5, 6) {$y_1^{(2)}\eps_{24}$};
            \node[crit] (yy24) at (8.5, 6) {$y_1^{(2)}\eps_{34}$};
            \node[crit] (yy34) at (6.5, 6) {$y_1^{(2)}\eps_{23}$};
            \node[crit] (yyyempty) at (10, 6) {$y_1^{(3)}\eps_\varnothing$};
            \draw[dedge] (e2) -- (eempty); \draw[dedge] (e3) -- (eempty); \draw[dedge] (e4) -- (eempty);
            \draw[dedge] (e12) -- (e1); \draw[dedge] (e12) -- (e2);
            \draw[dedge] (e13) -- (e1); \draw[dedge] (e13) -- (e3);
            \draw[dedge] (e14) -- (e1);
            \draw[dedge] (e23) -- (e2); \draw[dedge] (e23) -- (e3);
            \draw[dedge] (e24) -- (e2); \draw[dedge] (e24) -- (e4);
            \draw[dedge] (e34) -- (e3); \draw[dedge] (e34) -- (e4);
            \draw[dedge] (e123) -- (e12); \draw[dedge] (e123) -- (e13);
            \draw[dedge] (e124) -- (e12); \draw[dedge] (e124) -- (e14); \draw[dedge] (e124) -- (e24);
            \draw[dedge] (e134) -- (e13); \draw[dedge] (e134) -- (e34);
            \draw[dedge] (e234) -- (e23); \draw[dedge] (e234) -- (e24); \draw[dedge] (e234) -- (e34);
            \draw[dedge] (e1234) -- (e123); \draw[dedge] (e1234) -- (e124); \draw[dedge] (e1234) -- (e134);
            \draw[dedge] (y2) -- (yempty); \draw[dedge] (y3) -- (yempty); \draw[dedge] (y4) -- (yempty);
            \draw[dedge] (y12) -- (y1); \draw[dedge] (y12) -- (y2);
            \draw[dedge] (y13) -- (y1); \draw[dedge] (y13) -- (y3);
            \draw[dedge] (y14) -- (y1);
            \draw[dedge] (y23) -- (y2); \draw[dedge] (y23) -- (y3);
            \draw[dedge] (y24) -- (y2); \draw[dedge] (y24) -- (y4);
            \draw[dedge] (y34) -- (y3); \draw[dedge] (y34) -- (y4);
            \draw[dedge] (y123) -- (y12); \draw[dedge] (y123) -- (y13);
            \draw[dedge] (y124) -- (y12); \draw[dedge] (y124) -- (y14); \draw[dedge] (y124) -- (y24);
            \draw[dedge] (y134) -- (y13); \draw[dedge] (y134) -- (y34);
            \draw[dedge] (y234) -- (y23); \draw[dedge] (y234) -- (y24); \draw[dedge] (y234) -- (y34);
            \draw[dedge] (y1234) -- (y123); \draw[dedge] (y1234) -- (y124); \draw[dedge] (y1234) -- (y134);
            \draw[dedge] (yy2) -- (yyempty); \draw[dedge] (yy3) -- (yyempty); \draw[dedge] (yy4) -- (yyempty);
            \draw[dedge] (yy12) -- (yy1); \draw[dedge] (yy12) -- (yy2);
            \draw[dedge] (yy13) -- (yy1); \draw[dedge] (yy13) -- (yy3);
            \draw[dedge] (yy14) -- (yy1);
            \draw[dedge] (yy23) -- (yy2); \draw[dedge] (yy23) -- (yy3);
            \draw[dedge] (yy24) -- (yy2); \draw[dedge] (yy24) -- (yy4);
            \draw[dedge] (yy34) -- (yy3); \draw[dedge] (yy34) -- (yy4);
            \draw[sedge] (yempty) -- (e1);
            \draw[sedge] (y2) -- (e12); \draw[sedge] (y3) -- (e13); \draw[sedge] (y4) -- (e14);
            \draw[sedge] (y24) -- (e124); \draw[sedge] (y34) -- (e134);

            \draw[sedge] (yyempty) -- (y1);
            \draw[sedge] (yy2) -- (y12); \draw[sedge] (yy3) -- (y13); \draw[sedge] (yy4) -- (y14);
            \draw[sedge] (yy24) -- (y124); \draw[sedge] (yy34) -- (y134);
            \draw[sedge] (yyyempty) -- (yy1);
            \draw[medge] (e23) -- (e123);
            \draw[medge] (e14) -- (e134);
            \draw[medge] (e234) -- (e1234);
            \draw[medge] (y23) -- (y123);
            \draw[medge] (y14) -- (y134);
            \draw[medge] (y234) -- (y1234);
            \draw[medge, dashed] (yy14) -- ++(-0.35, 0.72) node[above left=-1pt, font=\tiny]{$y_1^{(2)}\eps_{123}$};
            \draw[medge, dashed] (yy23) -- ++(0.35, 0.72) node[above right=-1pt, font=\tiny]{$y_1^{(2)}\eps_{124}$};
        \end{tikzpicture}}
        \caption{$\Grep\bigl(\Sh(\F, \sigma), \lift{M}\bigr)$ through degree $6$, $c = 1$. Green: $\bar\partial_1$.}
        \label{fig:Grepc_16}
\end{figure}

The reduced complex $\Sh(\F, \sigma)^{\lift{M}}$ over $\Rbar = R/(xy)$ has ranks $1, 4, 5, 5, 5, \dots$ and matrices $$d_{\Sh}^{\lift{M}}|_1 = \begin{bmatrix}\overline{ d^{\M}_1} \end{bmatrix} = \ \ \begin{NiceArray}{ccccc}
        \eps_1 & \eps_2 & \eps_3 & \eps_4 &\\
        0 & yz & zw & xw & \eps_{\varnothing}\\
        &&&&
        \CodeAfter
  \SubMatrix[{2-1}{2-4}]
    \end{NiceArray},$$
    \[
        d_{\Sh}^{\lift{M}}|_2 = \left[\begin{array}{c|c} \overline{d^{\M}_2} & \overline{\sigma^{\M}|_0} \end{array}\right] =\ \  \begin{NiceArray}{cccc|cc}
           \eps_{12} &  \eps_{14} & \eps_{23} & \eps_{34} & y_1\eps_{\varnothing}&\\
            z & w & 0 & 0 & 1 & \eps_1 \\
            -x & 0 & w & 0 & 0 & \eps_2\\
            0 & 0 & -y & x & 0 & \eps_3\\
            0 & -y & 0 & -z & 0 & \eps_4,
        \CodeAfter
  \SubMatrix[{2-1}{5-5}]
  \end{NiceArray}
    \]

    \[
        d_{\Sh}^{\lift{M}}|_3 = \left[\begingroup
        \renewcommand{\arraystretch}{1.5}\begin{array}{c|c} \overline{d^{\M}_3} & \overline{\sigma^{\M}|_1}\\
            \hline
            0 & \overline{d^{\M}_1} \end{array}\endgroup\right] = \ \ \begin{NiceArray}{c|ccccc}
                \eps_{134} & y_1\eps_{1} & y_1\eps_2 & y_1\eps_3 & y_1\eps_4 &\\
                w &0 & -y & -w & 0 & \eps_{12}\\
                -z &  0 & 0 & 0 & -x & \eps_{14}\\
                x & 0 & 0 & -x & 0 & \eps_{23}\\
                y & 0 & 0 & 0 & 0 & \eps_{34}\\
                \hline
                0 & 0 & yz & zw & xw & y_1\eps_{\varnothing},
                \CodeAfter
  \SubMatrix[{2-1}{6-5}]
            \end{NiceArray}
    \]

    \[
        d_{\Sh}^{\lift{M}}|_4 =\left[\begingroup
        \renewcommand{\arraystretch}{1.5}\begin{array}{c|c} \overline{\sigma^{\M}|_2} & 0\\
            \hline
            \overline{d^{\M}_2} & \overline{\sigma^{\M}|_0} \end{array}\endgroup\right]
        = \ \ \begin{NiceArray}{cccc|cc}
            y_1\eps_{12} & y_1\eps_{14} & y_1\eps_{23} & y_1 \eps_{34} & y_1^{(2)}\eps_{\varnothing} & \\
            0 & 0 & 0 & x & 0 & \eps_{134}\\
                \hline
                z & w & 0 & 0 & 1 & y_1\eps_1\\
                -x & 0 & w & 0 & 0 & y_1\eps_2 \\
                0 & 0 & -y & x & 0 & y_1\eps_3\\
                0 & -y & 0 & -z & 0 & y_1\eps_4.\\
                \CodeAfter
  \SubMatrix[{2-1}{6-5}]
        \end{NiceArray}
    \]

   For all $q \geq 2$, $d_{\Sh}^{\lift{M}}|_{2q-1} = d_{\Sh}^{\lift{M}}|_3$ and $d_{\Sh}^{\lift{M}}|_{2q} = d_{\Sh}^{\lift{M}}|_4$.
    The bases for $d_{\Sh}^{\lift{M}}|_{2q-1}$ are
    \begin{center}
        columns
    $y_1^{(q-2)}\eps_{134}$, $y_1^{(q-1)} \eps_1$, $y_1^{(q-1)} \eps_2$, $y_1^{(q-1)} \eps_3$, $y_1^{(q-1)} \eps_4$;
    \end{center}

    \begin{center}
        rows $y_1^{(q-2)}\eps_{12}$, $y_1^{(q-2)}\eps_{14}$, $y_1^{(q-2)}\eps_{23}$, $y_1^{(q-2)}\eps_{34}$, $y_1^{(q-1)} \eps_\varnothing$.
    \end{center}
    The bases for $d_{\Sh}^{\lift{M}}|_{2q}$ are
    \begin{center}
        columns
     $y_1^{(q-1)} \eps_{12}$, $y_1^{(q-1)} \eps_{14}$, $y_1^{(q-1)} \eps_{23}$, $y_1^{(q-1)} \eps_{34}$, $y_1^{(q)} \eps_\varnothing$;
    \end{center}
    \begin{center}
        rows $y_1^{(q-2)}\eps_{134}$, $y_1^{(q-1)} \eps_1$, $y_1^{(q-1)} \eps_2$, $y_1^{(q-1)} \eps_3$, $y_1^{(q-1)} \eps_4$.
    \end{center}

    The only unit entries are those from $y_1^{(q)} \eps_\varnothing$ to $y_1^{(q-1)} \eps_1$ in each even degree $2q$, the edge $(y_1 \eps_\varnothing, \eps_1)$ of the diagram repeated on every level. Canceling those pairs by a second matching leaves no unit entry, so
    \[
        \beta^{\Rbar}_n(R/I) = 1, 3, 4, 4, 4, \dots \qquad \text{for } \Rbar = R/(xy).
    \]
\end{example}

\begin{example}[The running example, $c = 2$]\label{ex:runningtwo}
    Let our regular sequence consist of $a_1=m_1=xy$ (from Example \ref{ex:running}) and $a_2 = m_3 = zw$ and $\sigma^{(2)}$ be Sobieska's homotopy for $a_2$, keeping $\sigma^{(1)} = \sigma$ of Example~\ref{ex:running}.

    The nonzero values of $\sigma^{(2)}$ on cells are
    \[
        \begin{aligned}
            \sigma^{(2)}(\eps_\varnothing) &= \eps_3, & \sigma^{(2)}(\eps_1) &= \eps_{13}, & \sigma^{(2)}(\eps_2) &= z \eps_{23}, & \sigma^{(2)}(\eps_4) &= -w \eps_{34},\\
            \sigma^{(2)}(\eps_{12}) &= z \eps_{123}, & \sigma^{(2)}(\eps_{14}) &= -w \eps_{134}, & \sigma^{(2)}(\eps_{24}) &= -zw \eps_{234}, & \sigma^{(2)}(\eps_{124}) &= -zw \eps_{1234}.
        \end{aligned}
    \]
    In matrices,
    \[
        \sigma^{(2)}|_0 = \begin{bmatrix} 0\\ 0\\ 1\\ 0 \end{bmatrix}, \quad
        \sigma^{(2)}|_1 = \begin{bmatrix} 0 & 0 & 0 & 0\\ 1 & 0 & 0 & 0\\ 0 & 0 & 0 & 0\\ 0 & z & 0 & 0\\ 0 & 0 & 0 & 0\\ 0 & 0 & 0 & -w \end{bmatrix},
    \]

    \[
        \sigma^{(2)}|_2 = \begin{bmatrix} z & 0 & 0 & 0 & 0 & 0\\ 0 & 0 & 0 & 0 & 0 & 0\\ 0 & 0 & -w & 0 & 0 & 0\\ 0 & 0 & 0 & 0 & -zw & 0 \end{bmatrix}, \quad
        \sigma^{(2)}|_3 = \begin{bmatrix} 0 & -zw & 0 & 0 \end{bmatrix}.
    \]

We illustrate the complex graph $G(\F)$
including the arrows induced by $\sigma^{(1)}$ and $\sigma^{(2)}$ in Figure~\ref{fig:Gsigma2}.

\begin{figure}[htbp]
        \centering
        \resizebox{0.50\textwidth}{!}{%
        \begin{tikzpicture}[>=Stealth, x=1.32cm, y=1.5cm,
            vtx/.style={inner sep=1.6pt, font=\tiny, draw, rounded corners=1pt, fill=white, minimum width=6mm, minimum height=4mm},
            crit/.style={vtx, fill=gray!10},
            mu/.style={vtx, fill=red!12},
            md/.style={vtx, fill=orange!18},
            medge/.style={->, red!70!black, thick},
            sedge/.style={->, green!55!black, thick, densely dashed},
            tedge/.style={->, blue!65!black, thick, densely dashed},
            dedge/.style={->, gray!50, thin}]
            \node[crit] (eempty) at (0, 0) {$\eps_\varnothing$};
            \node[crit] (e1) at (-1.5, 1) {$\eps_1$};
            \node[crit] (e2) at (1.5, 1) {$\eps_4$};
            \node[crit] (e3) at (-0.5, 1) {$\eps_2$};
            \node[crit] (e4) at (0.5, 1) {$\eps_3$};
            \node[crit] (e12) at (-0.5, 2) {$\eps_{14}$};
            \node[crit] (e13) at (-2.5, 2) {$\eps_{12}$};
            \node[md] (e14) at (-1.5, 2) {$\eps_{13}$};
            \node[md] (e23) at (1.5, 2) {$\eps_{24}$};
            \node[crit] (e24) at (2.5, 2) {$\eps_{34}$};
            \node[crit] (e34) at (0.5, 2) {$\eps_{23}$};
            \node[mu] (e123) at (-0.5, 3) {$\eps_{124}$};
            \node[crit] (e124) at (0.5, 3) {$\eps_{134}$};
            \node[mu] (e134) at (-1.5, 3) {$\eps_{123}$};
            \node[md] (e234) at (1.5, 3) {$\eps_{234}$};
            \node[mu] (e1234) at (-1, 4) {$\eps_{1234}$};
            \draw[dedge] (e1) -- (eempty);
            \draw[dedge] (e3) -- (eempty);
            \draw[dedge] (e4) -- (eempty);
            \draw[dedge] (e2) -- (eempty);
            \draw[dedge] (e13) -- (e3);
            \draw[dedge] (e13) -- (e1);
            \draw[dedge] (e14) -- (e4);
            \draw[dedge] (e14) -- (e1);
            \draw[dedge] (e12) -- (e2);
            \draw[dedge] (e12) -- (e1);
            \draw[dedge] (e34) -- (e4);
            \draw[dedge] (e34) -- (e3);
            \draw[dedge] (e23) -- (e2);
            \draw[dedge] (e23) -- (e3);
            \draw[dedge] (e24) -- (e2);
            \draw[dedge] (e24) -- (e4);
            \draw[dedge] (e134) -- (e34);
            \draw[medge] (e14) -- (e134);
            \draw[dedge] (e134) -- (e13);
            \draw[medge] (e23) -- (e123);
            \draw[dedge] (e123) -- (e12);
            \draw[dedge] (e123) -- (e13);
            \draw[dedge] (e124) -- (e24);
            \draw[dedge] (e124) -- (e12);
            \draw[dedge] (e124) -- (e14);
            \draw[dedge] (e234) -- (e24);
            \draw[dedge] (e234) -- (e23);
            \draw[dedge] (e234) -- (e34);
            \draw[medge] (e234) -- (e1234);
            \draw[dedge] (e1234) -- (e124);
            \draw[dedge] (e1234) -- (e123);
            \draw[dedge] (e1234) -- (e134);
            \draw[sedge] (eempty) to[bend left=18] (e1);
            \draw[sedge] (e3) to[bend left=18] (e13);
            \draw[sedge] (e4) to[bend left=18] (e14);
            \draw[sedge] (e2) to[bend left=18] (e12);
            \draw[sedge] (e34) to[bend left=18] (e134);
            \draw[sedge] (e23) to[bend left=18] (e123);
            \draw[sedge] (e24) to[bend left=18] (e124);
            \draw[sedge] (e234) to[bend left=18] (e1234);
            \draw[tedge] (eempty) to[bend right=18] (e4);
            \draw[tedge] (e1) to[bend right=18] (e14);
            \draw[tedge] (e3) to[bend right=18] (e34);
            \draw[tedge] (e2) to[bend right=18] (e24);
            \draw[tedge] (e13) to[bend right=18] (e134);
            \draw[tedge] (e12) to[bend right=18] (e124);
            \draw[tedge] (e23) to[bend right=18] (e234);
            \draw[tedge] (e123) to[bend right=18] (e1234);
        \end{tikzpicture}}
        \caption{$G(\F)$ with the arrows of $\sigma^{(1)}$ and $\sigma^{(2)}$, $c = 2$.}
        \label{fig:Gsigma2}
\end{figure}

    Over $R/(xy, zw)$, $\Sh(\F, \sigma)$ has ranks $1, 4, 8, 12, 16, \dots$, and $d_{\Sh}$ is block bidiagonal with connecting blocks $S_k = \sum_j T^{(k)}_{\boldsymbol{e}_j} \otimes \sigma^{(j)}$ (Proposition~\ref{prop:blocks}). With cells ordered by level and within a level by $\boldsymbol{u} = (k, 0), \dots, (0, k)$, its matrices in degrees $2$ and $3$ are built from $S_1|_q = \begin{bmatrix} \sigma^{(1)}|_q & \sigma^{(2)}|_q \end{bmatrix}$:

        \[
        d_{\Sh}|_2 = \left[\begin{array}{c|c} \overline{d_2} & \overline{S_1|_0} \end{array}\right] = \ \ \begin{NiceArray}{cccccc|ccc}
          \eps_{12} & \eps_{13} & \eps_{14} & \eps_{23} & \eps_{24} & \eps_{34} & y_1\eps_{\varnothing} & y_2 \eps_{\varnothing}\\
           z & 0 & w & 0 & 0 & 0 & 1 & 0 & \eps_1\\ -x & 0 & 0 & w & xw & 0 & 0 & 0 & \eps_2\\ 0 & 0 & 0 & -y & 0 & x & 0 & 1 & \eps_3\\ 0 & 0 & -y & 0 & -yz & -z & 0 & 0 & \eps_4
          \CodeAfter
  \SubMatrix[{2-1}{5-8}]
        \end{NiceArray}
    \]

\begin{figure}[p]
        \centering
        \definecolor{sh00}{HTML}{D2F1F7}\definecolor{tx00}{HTML}{1F6F80}\definecolor{sh10}{HTML}{FBF0B8}\definecolor{tx10}{HTML}{7A6400}\definecolor{sh01}{HTML}{CFEDE4}\definecolor{tx01}{HTML}{1E6B57}\definecolor{sh20}{HTML}{EADCF3}\definecolor{tx20}{HTML}{5E3A80}\definecolor{sh11}{HTML}{F8DDE9}\definecolor{tx11}{HTML}{8A2E5A}\definecolor{sh02}{HTML}{DCE5F7}\definecolor{tx02}{HTML}{2F4D8A}\definecolor{sh30}{HTML}{E4F2C4}\definecolor{tx30}{HTML}{4F6B12}\definecolor{sh21}{HTML}{EEE0CC}\definecolor{tx21}{HTML}{6E4A1E}\definecolor{sh12}{HTML}{E2E2C0}\definecolor{tx12}{HTML}{5C5C1A}\definecolor{sh03}{HTML}{D9D9E0}\definecolor{tx03}{HTML}{45455A}
        \rotatebox{90}{\resizebox{0.88\textheight}{!}{%
        \begin{tikzpicture}[>=Stealth, x=1.00cm, y=3.00cm,
            vtx/.style={circle, inner sep=0pt, minimum size=2.4mm, draw, fill=white},
            crit/.style={vtx, fill=gray!35},
            mu/.style={vtx, fill=red!55},
            md/.style={vtx, fill=orange!60},
            medge/.style={->, red!70!black, thick},
            sedge/.style={->, green!55!black, thick},
            tedge/.style={->, blue!65!black, thick},
            dedge/.style={->, gray!50, thin}]
            \fill[sh00, rounded corners=4pt] (-0.42,-0.40) rectangle (0.42,0.40);
            \fill[sh00] (-0.37,0.35) rectangle (0.37,0.65);
            \fill[sh00, rounded corners=4pt] (-1.92,0.60) rectangle (1.92,1.40);
            \fill[sh00] (-1.87,1.35) rectangle (1.87,1.65);
            \fill[sh00, rounded corners=4pt] (-2.92,1.60) rectangle (2.92,2.40);
            \fill[sh00] (-1.87,2.35) rectangle (1.87,2.65);
            \fill[sh00, rounded corners=4pt] (-1.92,2.60) rectangle (1.92,3.40);
            \fill[sh00] (-1.18,3.35) rectangle (-0.42,3.65);
            \fill[sh00, rounded corners=4pt] (-1.23,3.60) rectangle (-0.37,4.40);
            \fill[sh01, rounded corners=4pt] (3.68,1.60) rectangle (4.52,2.40);
            \fill[sh01] (3.73,2.35) rectangle (4.47,2.65);
            \fill[sh01, rounded corners=4pt] (2.68,2.60) rectangle (6.52,3.40);
            \fill[sh01] (2.73,3.35) rectangle (6.47,3.65);
            \fill[sh01, rounded corners=4pt] (1.97,3.60) rectangle (7.83,4.40);
            \fill[sh01] (2.68,4.35) rectangle (6.52,4.65);
            \fill[sh01, rounded corners=4pt] (2.63,4.60) rectangle (6.57,5.40);
            \fill[sh01] (3.64,5.35) rectangle (4.56,5.65);
            \fill[sh01, rounded corners=4pt] (3.59,5.60) rectangle (4.61,6.40);
            \fill[sh02, rounded corners=4pt] (8.57,3.60) rectangle (9.43,4.40);
            \fill[sh02] (8.62,4.35) rectangle (9.38,4.65);
            \fill[sh02, rounded corners=4pt] (7.27,4.60) rectangle (11.13,5.40);
            \fill[sh02] (7.32,5.35) rectangle (11.08,5.65);
            \fill[sh02, rounded corners=4pt] (6.83,5.60) rectangle (12.77,6.40);
            \fill[sh03, rounded corners=4pt] (13.47,5.60) rectangle (14.33,6.40);
            \fill[sh10, rounded corners=4pt] (-4.52,1.60) rectangle (-3.68,2.40);
            \fill[sh10] (-4.47,2.35) rectangle (-3.73,2.65);
            \fill[sh10, rounded corners=4pt] (-6.52,2.60) rectangle (-2.68,3.40);
            \fill[sh10] (-6.47,3.35) rectangle (-2.73,3.65);
            \fill[sh10, rounded corners=4pt] (-7.83,3.60) rectangle (-1.97,4.40);
            \fill[sh10] (-6.52,4.35) rectangle (-2.68,4.65);
            \fill[sh10, rounded corners=4pt] (-6.57,4.60) rectangle (-2.63,5.40);
            \fill[sh10] (-6.16,5.35) rectangle (-5.24,5.65);
            \fill[sh10, rounded corners=4pt] (-6.21,5.60) rectangle (-5.19,6.40);
            \fill[sh11, rounded corners=4pt] (0.33,3.60) rectangle (1.27,4.40);
            \fill[sh11] (0.38,4.35) rectangle (1.22,4.65);
            \fill[sh11, rounded corners=4pt] (-1.97,4.60) rectangle (1.97,5.40);
            \fill[sh11] (-1.92,5.35) rectangle (1.92,5.65);
            \fill[sh11, rounded corners=4pt] (-3.01,5.60) rectangle (3.01,6.40);
            \fill[sh12, rounded corners=4pt] (5.19,5.60) rectangle (6.21,6.40);
            \fill[sh20, rounded corners=4pt] (-9.43,3.60) rectangle (-8.57,4.40);
            \fill[sh20] (-9.38,4.35) rectangle (-8.62,4.65);
            \fill[sh20, rounded corners=4pt] (-11.13,4.60) rectangle (-7.27,5.40);
            \fill[sh20] (-11.08,5.35) rectangle (-7.32,5.65);
            \fill[sh20, rounded corners=4pt] (-12.77,5.60) rectangle (-6.83,6.40);
            \fill[sh21, rounded corners=4pt] (-4.61,5.60) rectangle (-3.59,6.40);
            \fill[sh30, rounded corners=4pt] (-14.33,5.60) rectangle (-13.47,6.40);
            \node[font=\Large, text=tx00, fill=white, fill opacity=0.85, text opacity=1, rounded corners=1pt, anchor=north, inner sep=1.2pt] at (0.00,-0.40) {$F_\bullet$};
            \node[font=\Large, text=tx01, fill=white, fill opacity=0.85, text opacity=1, rounded corners=1pt, anchor=north, inner sep=1.2pt] at (4.10,1.60) {$y_2 F_\bullet$};
            \node[font=\Large, text=tx02, fill=white, fill opacity=0.85, text opacity=1, rounded corners=1pt, anchor=south, inner sep=1.2pt] at (9.80,6.45) {$y_2^{(2)} F_\bullet$};
            \node[font=\Large, text=tx03, fill=white, fill opacity=0.85, text opacity=1, rounded corners=1pt, anchor=south, inner sep=1.2pt] at (13.90,6.45) {$y_2^{(3)} F_\bullet$};
            \node[font=\Large, text=tx10, fill=white, fill opacity=0.85, text opacity=1, rounded corners=1pt, anchor=north, inner sep=1.2pt] at (-4.10,1.60) {$y_1 F_\bullet$};
            \node[font=\Large, text=tx11, fill=white, fill opacity=0.85, text opacity=1, rounded corners=1pt, anchor=south, inner sep=1.2pt] at (0.00,6.45) {$y_1y_2 F_\bullet$};
            \node[font=\Large, text=tx12, fill=white, fill opacity=0.85, text opacity=1, rounded corners=1pt, anchor=south, inner sep=1.2pt] at (5.70,6.45) {$y_1y_2^{(2)} F_\bullet$};
            \node[font=\Large, text=tx20, fill=white, fill opacity=0.85, text opacity=1, rounded corners=1pt, anchor=south, inner sep=1.2pt] at (-9.80,6.45) {$y_1^{(2)} F_\bullet$};
            \node[font=\Large, text=tx21, fill=white, fill opacity=0.85, text opacity=1, rounded corners=1pt, anchor=south, inner sep=1.2pt] at (-4.10,6.45) {$y_1^{(2)}y_2 F_\bullet$};
            \node[font=\Large, text=tx30, fill=white, fill opacity=0.85, text opacity=1, rounded corners=1pt, anchor=south, inner sep=1.2pt] at (-13.90,6.45) {$y_1^{(3)} F_\bullet$};
            \node[crit] (u00s0) at (0.00, 0) {};
            \node[crit] (u00s1) at (-1.50, 1) {};
            \node[crit] (u00s2) at (1.50, 1) {};
            \node[crit] (u00s3) at (-0.50, 1) {};
            \node[crit] (u00s4) at (0.50, 1) {};
            \node[crit] (u10s0) at (-4.10, 2) {};
            \node[crit] (u00s12) at (-0.50, 2) {};
            \node[crit] (u00s13) at (-2.50, 2) {};
            \node[md] (u00s14) at (-1.50, 2) {};
            \node[md] (u00s23) at (1.50, 2) {};
            \node[crit] (u00s24) at (2.50, 2) {};
            \node[crit] (u00s34) at (0.50, 2) {};
            \node[crit] (u01s0) at (4.10, 2) {};
            \node[crit] (u10s1) at (-6.10, 3) {};
            \node[crit] (u10s2) at (-3.10, 3) {};
            \node[crit] (u10s3) at (-5.10, 3) {};
            \node[crit] (u10s4) at (-4.10, 3) {};
            \node[mu] (u00s123) at (-0.50, 3) {};
            \node[crit] (u00s124) at (0.50, 3) {};
            \node[mu] (u00s134) at (-1.50, 3) {};
            \node[md] (u00s234) at (1.50, 3) {};
            \node[crit] (u01s1) at (3.10, 3) {};
            \node[crit] (u01s2) at (6.10, 3) {};
            \node[crit] (u01s3) at (4.10, 3) {};
            \node[crit] (u01s4) at (5.10, 3) {};
            \node[crit] (u20s0) at (-9.00, 4) {};
            \node[crit] (u10s12) at (-5.40, 4) {};
            \node[crit] (u10s13) at (-7.40, 4) {};
            \node[md] (u10s14) at (-6.40, 4) {};
            \node[md] (u10s23) at (-3.40, 4) {};
            \node[crit] (u10s24) at (-2.40, 4) {};
            \node[crit] (u10s34) at (-4.40, 4) {};
            \node[mu] (u00s1234) at (-0.80, 4) {};
            \node[crit] (u11s0) at (0.80, 4) {};
            \node[crit] (u01s12) at (4.40, 4) {};
            \node[crit] (u01s13) at (2.40, 4) {};
            \node[md] (u01s14) at (3.40, 4) {};
            \node[md] (u01s23) at (6.40, 4) {};
            \node[crit] (u01s24) at (7.40, 4) {};
            \node[crit] (u01s34) at (5.40, 4) {};
            \node[crit] (u02s0) at (9.00, 4) {};
            \node[crit] (u20s1) at (-10.70, 5) {};
            \node[crit] (u20s2) at (-7.70, 5) {};
            \node[crit] (u20s3) at (-9.70, 5) {};
            \node[crit] (u20s4) at (-8.70, 5) {};
            \node[mu] (u10s123) at (-5.10, 5) {};
            \node[crit] (u10s124) at (-4.10, 5) {};
            \node[mu] (u10s134) at (-6.10, 5) {};
            \node[md] (u10s234) at (-3.10, 5) {};
            \node[crit] (u11s1) at (-1.50, 5) {};
            \node[crit] (u11s2) at (1.50, 5) {};
            \node[crit] (u11s3) at (-0.50, 5) {};
            \node[crit] (u11s4) at (0.50, 5) {};
            \node[mu] (u01s123) at (4.10, 5) {};
            \node[crit] (u01s124) at (5.10, 5) {};
            \node[mu] (u01s134) at (3.10, 5) {};
            \node[md] (u01s234) at (6.10, 5) {};
            \node[crit] (u02s1) at (7.70, 5) {};
            \node[crit] (u02s2) at (10.70, 5) {};
            \node[crit] (u02s3) at (8.70, 5) {};
            \node[crit] (u02s4) at (9.70, 5) {};
            \node[crit] (u30s0) at (-13.90, 6) {};
            \node[crit] (u20s12) at (-10.30, 6) {};
            \node[crit] (u20s13) at (-12.30, 6) {};
            \node[md] (u20s14) at (-11.30, 6) {};
            \node[md] (u20s23) at (-8.30, 6) {};
            \node[crit] (u20s24) at (-7.30, 6) {};
            \node[crit] (u20s34) at (-9.30, 6) {};
            \node[mu] (u10s1234) at (-5.70, 6) {};
            \node[crit] (u21s0) at (-4.10, 6) {};
            \node[crit] (u11s12) at (-0.50, 6) {};
            \node[crit] (u11s13) at (-2.50, 6) {};
            \node[md] (u11s14) at (-1.50, 6) {};
            \node[md] (u11s23) at (1.50, 6) {};
            \node[crit] (u11s24) at (2.50, 6) {};
            \node[crit] (u11s34) at (0.50, 6) {};
            \node[mu] (u01s1234) at (4.10, 6) {};
            \node[crit] (u12s0) at (5.70, 6) {};
            \node[crit] (u02s12) at (9.30, 6) {};
            \node[crit] (u02s13) at (7.30, 6) {};
            \node[md] (u02s14) at (8.30, 6) {};
            \node[md] (u02s23) at (11.30, 6) {};
            \node[crit] (u02s24) at (12.30, 6) {};
            \node[crit] (u02s34) at (10.30, 6) {};
            \node[crit] (u03s0) at (13.90, 6) {};

            \draw[dedge] (u00s2) -- (u00s0);
            \draw[dedge] (u00s3) -- (u00s0);

            \draw[dedge] (u00s12) -- (u00s2);
            \draw[dedge] (u00s12) -- (u00s1);
            \draw[dedge] (u00s13) -- (u00s3);
            \draw[dedge] (u00s13) -- (u00s1);

            \draw[dedge] (u00s23) -- (u00s3);
            \draw[dedge] (u00s23) -- (u00s2);
            \draw[dedge] (u00s24) -- (u00s4);
            \draw[dedge] (u00s24) -- (u00s2);
            \draw[dedge] (u00s34) -- (u00s4);
            \draw[dedge] (u00s34) -- (u00s3);

            \draw[dedge] (u10s2) -- (u10s0);
            \draw[dedge] (u10s3) -- (u10s0);

            \draw[dedge] (u00s123) -- (u00s13);
            \draw[dedge] (u00s123) -- (u00s12);
            \draw[dedge] (u00s124) -- (u00s24);
            \draw[dedge] (u00s124) -- (u00s14);
            \draw[dedge] (u00s124) -- (u00s12);
            \draw[dedge] (u00s134) -- (u00s34);

            \draw[dedge] (u00s134) -- (u00s13);
            \draw[dedge] (u00s234) -- (u00s34);
            \draw[dedge] (u00s234) -- (u00s24);
            \draw[dedge] (u00s234) -- (u00s23);

            \draw[dedge] (u01s2) -- (u01s0);
            \draw[dedge] (u01s3) -- (u01s0);

            \draw[dedge] (u10s12) -- (u10s2);
            \draw[dedge] (u10s12) -- (u10s1);
            \draw[dedge] (u10s13) -- (u10s3);
            \draw[dedge] (u10s13) -- (u10s1);

            \draw[dedge] (u10s23) -- (u10s3);
            \draw[dedge] (u10s23) -- (u10s2);
            \draw[dedge] (u10s24) -- (u10s4);
            \draw[dedge] (u10s24) -- (u10s2);
            \draw[dedge] (u10s34) -- (u10s4);
            \draw[dedge] (u10s34) -- (u10s3);

            \draw[dedge] (u00s1234) -- (u00s134);
            \draw[dedge] (u00s1234) -- (u00s124);
            \draw[dedge] (u00s1234) -- (u00s123);
            \draw[dedge] (u01s12) -- (u01s2);
            \draw[dedge] (u01s12) -- (u01s1);
            \draw[dedge] (u01s13) -- (u01s3);
            \draw[dedge] (u01s13) -- (u01s1);

            \draw[dedge] (u01s23) -- (u01s3);
            \draw[dedge] (u01s23) -- (u01s2);
            \draw[dedge] (u01s24) -- (u01s4);
            \draw[dedge] (u01s24) -- (u01s2);
            \draw[dedge] (u01s34) -- (u01s4);
            \draw[dedge] (u01s34) -- (u01s3);

            \draw[dedge] (u20s2) -- (u20s0);
            \draw[dedge] (u20s3) -- (u20s0);

            \draw[dedge] (u10s123) -- (u10s13);
            \draw[dedge] (u10s123) -- (u10s12);
            \draw[dedge] (u10s124) -- (u10s24);
            \draw[dedge] (u10s124) -- (u10s14);
            \draw[dedge] (u10s124) -- (u10s12);
            \draw[dedge] (u10s134) -- (u10s34);

            \draw[dedge] (u10s134) -- (u10s13);
            \draw[dedge] (u10s234) -- (u10s34);
            \draw[dedge] (u10s234) -- (u10s24);
            \draw[dedge] (u10s234) -- (u10s23);

            \draw[dedge] (u11s2) -- (u11s0);
            \draw[dedge] (u11s3) -- (u11s0);

            \draw[dedge] (u01s123) -- (u01s13);
            \draw[dedge] (u01s123) -- (u01s12);
            \draw[dedge] (u01s124) -- (u01s24);
            \draw[dedge] (u01s124) -- (u01s14);
            \draw[dedge] (u01s124) -- (u01s12);
            \draw[dedge] (u01s134) -- (u01s34);

            \draw[dedge] (u01s134) -- (u01s13);
            \draw[dedge] (u01s234) -- (u01s34);
            \draw[dedge] (u01s234) -- (u01s24);
            \draw[dedge] (u01s234) -- (u01s23);

            \draw[dedge] (u02s2) -- (u02s0);
            \draw[dedge] (u02s3) -- (u02s0);

            \draw[dedge] (u20s12) -- (u20s2);
            \draw[dedge] (u20s12) -- (u20s1);
            \draw[dedge] (u20s13) -- (u20s3);
            \draw[dedge] (u20s13) -- (u20s1);

            \draw[dedge] (u20s23) -- (u20s3);
            \draw[dedge] (u20s23) -- (u20s2);
            \draw[dedge] (u20s24) -- (u20s4);
            \draw[dedge] (u20s24) -- (u20s2);
            \draw[dedge] (u20s34) -- (u20s4);
            \draw[dedge] (u20s34) -- (u20s3);

            \draw[dedge] (u10s1234) -- (u10s134);
            \draw[dedge] (u10s1234) -- (u10s124);
            \draw[dedge] (u10s1234) -- (u10s123);
            \draw[dedge] (u11s12) -- (u11s2);
            \draw[dedge] (u11s12) -- (u11s1);
            \draw[dedge] (u11s13) -- (u11s3);
            \draw[dedge] (u11s13) -- (u11s1);

            \draw[dedge] (u11s23) -- (u11s3);
            \draw[dedge] (u11s23) -- (u11s2);
            \draw[dedge] (u11s24) -- (u11s4);
            \draw[dedge] (u11s24) -- (u11s2);
            \draw[dedge] (u11s34) -- (u11s4);
            \draw[dedge] (u11s34) -- (u11s3);

            \draw[dedge] (u01s1234) -- (u01s134);
            \draw[dedge] (u01s1234) -- (u01s124);
            \draw[dedge] (u01s1234) -- (u01s123);
            \draw[dedge] (u02s12) -- (u02s2);
            \draw[dedge] (u02s12) -- (u02s1);
            \draw[dedge] (u02s13) -- (u02s3);
            \draw[dedge] (u02s13) -- (u02s1);

            \draw[dedge] (u02s23) -- (u02s3);
            \draw[dedge] (u02s23) -- (u02s2);
            \draw[dedge] (u02s24) -- (u02s4);
            \draw[dedge] (u02s24) -- (u02s2);
            \draw[dedge] (u02s34) -- (u02s4);
            \draw[dedge] (u02s34) -- (u02s3);
            \draw[sedge] (u10s0) -- (u00s1);
            \draw[tedge] (u01s0) -- (u00s4);
            \draw[sedge] (u10s2) -- (u00s12);
            \draw[sedge] (u10s3) -- (u00s13);
            \draw[sedge] (u10s4) -- (u00s14);
            \draw[tedge] (u01s1) -- (u00s14);
            \draw[tedge] (u01s2) -- (u00s24);
            \draw[tedge] (u01s3) -- (u00s34);
            \draw[sedge] (u20s0) -- (u10s1);

            \draw[sedge] (u10s24) -- (u00s124);
            \draw[sedge] (u10s34) -- (u00s134);
            \draw[sedge] (u11s0) -- (u01s1);
            \draw[tedge] (u11s0) -- (u10s4);
            \draw[tedge] (u01s12) -- (u00s124);
            \draw[tedge] (u01s13) -- (u00s134);

            \draw[tedge] (u02s0) -- (u01s4);
            \draw[sedge] (u20s2) -- (u10s12);
            \draw[sedge] (u20s3) -- (u10s13);
            \draw[sedge] (u20s4) -- (u10s14);

            \draw[tedge] (u11s1) -- (u10s14);
            \draw[sedge] (u11s2) -- (u01s12);
            \draw[tedge] (u11s2) -- (u10s24);
            \draw[sedge] (u11s3) -- (u01s13);
            \draw[tedge] (u11s3) -- (u10s34);
            \draw[sedge] (u11s4) -- (u01s14);

            \draw[tedge] (u02s1) -- (u01s14);
            \draw[tedge] (u02s2) -- (u01s24);
            \draw[tedge] (u02s3) -- (u01s34);
            \draw[sedge] (u30s0) -- (u20s1);

            \draw[sedge] (u20s24) -- (u10s124);
            \draw[sedge] (u20s34) -- (u10s134);
            \draw[sedge] (u21s0) -- (u11s1);
            \draw[tedge] (u21s0) -- (u20s4);
            \draw[tedge] (u11s12) -- (u10s124);
            \draw[tedge] (u11s13) -- (u10s134);

            \draw[sedge] (u11s24) -- (u01s124);
            \draw[sedge] (u11s34) -- (u01s134);
            \draw[sedge] (u12s0) -- (u02s1);
            \draw[tedge] (u12s0) -- (u11s4);
            \draw[tedge] (u02s12) -- (u01s124);
            \draw[tedge] (u02s13) -- (u01s134);

            \draw[tedge] (u03s0) -- (u02s4);
            \draw[medge] (u00s14) -- (u00s134);
            \draw[medge] (u00s23) -- (u00s123);
            \draw[medge] (u00s234) -- (u00s1234);
            \draw[medge] (u10s14) -- (u10s134);
            \draw[medge] (u10s23) -- (u10s123);
            \draw[medge] (u01s14) -- (u01s134);
            \draw[medge] (u01s23) -- (u01s123);
            \draw[medge] (u10s234) -- (u10s1234);
            \draw[medge] (u01s234) -- (u01s1234);
            \draw[medge, dashed] (u20s14) -- ++(-0.35, 0.25);
            \draw[medge, dashed] (u20s23) -- ++(0.35, 0.25);
            \draw[medge, dashed] (u11s14) -- ++(-0.35, 0.25);
            \draw[medge, dashed] (u11s23) -- ++(0.35, 0.25);
            \draw[medge, dashed] (u02s14) -- ++(-0.35, 0.25);
            \draw[medge, dashed] (u02s23) -- ++(0.35, 0.25);
            \node[anchor=east, font=\small, gray!80!black] at (-15.30, 0) {$n = 0$};
            \node[anchor=east, font=\small, gray!80!black] at (-15.30, 1) {$n = 1$};
            \node[anchor=east, font=\small, gray!80!black] at (-15.30, 2) {$n = 2$};
            \node[anchor=east, font=\small, gray!80!black] at (-15.30, 3) {$n = 3$};
            \node[anchor=east, font=\small, gray!80!black] at (-15.30, 4) {$n = 4$};
            \node[anchor=east, font=\small, gray!80!black] at (-15.30, 5) {$n = 5$};
            \node[anchor=east, font=\small, gray!80!black] at (-15.30, 6) {$n = 6$};
        \end{tikzpicture}
        }}
        \caption{$\Grep\bigl(\Sh(\F, \sigma), \lift{M}\bigr)$ through degree $6$, $c = 2$. Green: $\bar\partial_1$. Blue: $\bar\partial_2$.}\label{fig:ctwo}
\end{figure}

    \[
        d_{\Sh}|_3 = \left[\begingroup
        \renewcommand{\arraystretch}{1.5}\begin{array}{c|c} \overline{d_3} & \overline{S_1|_1}\\
            \hline
            0 & \overline{\mathbf{1}_2 \otimes d_1} \end{array}\endgroup\right]
        \]
which equals
\begin{center}\small $\begin{NiceArray}{cccc|cccc|ccccc}
            \eps_{123} & \eps_{124} & \eps_{134} & \eps_{234} & y_1\eps_1 & y_1\eps_2 & y_1\eps_3 & y_1\eps_4 & y_2\eps_1 & y_2\eps_2 & y_2\eps_3 & y_2\eps_4 &\\
            w & w & 0 & 0 & 0 & -y & 0 & 0 & 0 & 0 & 0 & 0 & \eps_{12}\\
            -1 & 0 & 1 & 0 & 0 & 0 & -1 & 0 & 1 & 0 & 0 & 0 & \eps_{13} \\
            0 & -z & -z & 0 & 0 & 0 & 0 & -x & 0 & 0 & 0 & 0 & \eps_{14} \\
            x & 0 & 0 & x & 0 & 0 & 0 & 0 & 0 & z & 0 & 0 & \eps_{23} \\
            0 & 1 & 0 & -1 & 0 & 0 & 0 & 0 & 0 & 0 & 0 & 0 & \eps_{24} \\
            0 & 0 & y & y & 0 & 0 & 0 & 0 & 0 & 0 & 0 & -w & \eps_{34}\\
            \hline
            0 & 0 & 0 & 0 & 0 & yz & 0 & xw & 0 & 0 & 0 & 0 & y_1\eps_{\varnothing} \\
            0 & 0 & 0 & 0 & 0 & 0 & 0 & 0 & 0 & yz & 0 & xw & y_2\eps_{\varnothing}
            \CodeAfter
  \SubMatrix[{2-1}{9-12}]
        \end{NiceArray}$\end{center}

Degree 4 is the first degree where an $S_2|_q=\left[\begingroup\renewcommand{\arraystretch}{1.5}\begin{array}{c|c|c}
    \sigma^{(1)}|_q & \sigma^{(2)}|_q & 0\\
    \hline
    0 & \sigma^{(1)}|_q & \sigma^{(2)}|_q \\
    \end{array}\endgroup\right]$ appears. We include below, but will not include the matrices in higher homological degree.

     \[
        d_{\Sh}|_4 = \left[\begingroup\renewcommand{\arraystretch}{1.5}\begin{array}{c|c|c} \overline{d_4} & \overline{S_1|_2} & 0\\
            \hline
            0 & \overline{\mathbf{1}_2 \otimes d_2}  & \overline{S_2|_0}\end{array}\endgroup\right] \]
         which equals
\begin{center}\footnotesize\setlength{\arraycolsep}{2.5pt}$\begin{NiceArray}{c|cccccc|cccccc|c|c|cc}
            \eps_{1234} & y_1\eps_{12} & y_1\eps_{13} & y_1\eps_{14} & y_1\eps_{23} & y_1\eps_{24} & y_1\eps_{34} &
            y_2\eps_{12} & y_2\eps_{13} & y_2\eps_{14} & y_2\eps_{23} & y_2\eps_{24} & y_2\eps_{34} & y_1^{(2)}\eps_{\varnothing} & y_1y_2\eps_{\varnothing} & y_2^{(2)}\eps_{\varnothing} &\\
            1 & 0 & 0 & 0 & y & 0 & 0 &  z & 0 & 0 & 0 & 0 & 0 & 0 & 0 & 0 & \eps_{123}\\
         -1 & 0 & 0 & 0 & 0 & 0 & 0 & 0 & 0 & 0 & 0 & 0 & 0 & 0 & 0 & 0 & \eps_{124}\\
         1 & 0 & 0 & 0 & 0 & 0 & x & 0 & 0 & -w & 0 & 0 & 0 & 0 & 0 & 0 & \eps_{134}\\
         -1 & 0 & 0 & 0 & 0 & 0 & 0 & 0 & 0 & 0 & 0 & 0 & 0 & 0 & 0 & 0 & \eps_{234}\\
         \hline
         0 & z & 0 & w & 0 & 0 & 0 & 0 & 0 & 0 & 0 & 0 & 0 & 1 & 0 & 0 & y_1\eps_1 \\
         0 & -x & 0 & 0 & w & xw & 0 & 0 & 0 & 0 & 0 & 0 & 0 & 0 & 0 & 0 & y_1 \eps_2 \\
         0 & 0 & 0 & 0 & -y & 0 & x & 0 & 0 & 0 & 0 & 0 & 0 & 0 & 1 & 0 & y_1\eps_3 \\
         0 & 0 & 0 & -y & 0 & -yz & -z & 0 & 0 & 0 & 0 & 0 & 0 & 0 & 0 & 0 & y_1 \eps_4\\
         \hline
         0 & 0 & 0 & 0 & 0 & 0 & 0 &z & 0 & w & 0 & 0 & 0 & 0& 1 & 0 & y_2\eps_1\\
         0 & 0 & 0 & 0 & 0 & 0 & 0 & -x & 0 & 0 & w & xw & 0 & 0 & 0 & 0 & y_2\eps_2\\
         0 & 0 & 0 & 0 & 0 & 0 & 0 & 0 & 0 & 0 & -y & 0 & x & 0 & 0 & 1 & y_2\eps_3\\
         0 & 0 & 0 & 0 & 0 & 0 & 0 & 0 & 0 & -y & 0 & -yz & -z &0 & 0 & 0 & y_2\eps_4
         \CodeAfter
  \SubMatrix[{2-1}{13-16}]
         \end{NiceArray}$\end{center}

As in Example \ref{ex:running}, the edges of $\lift{M}$ are $$ \{(y_1^{(p)}y_2^{(q)}\eps_{123}, y_1^{(p)}y_2^{(q)}\eps_{13}), (y_1^{(p)}y_2^{(q)}\eps_{124}, y_1^{(p)}y_2^{(q)}\eps_{24}), (y_1^{(p)}y_2^{(q)}\eps_{1234}, y_1^{(p)}y_2^{(q)}\eps_{234})\}_{p,q \geq 0}.$$
There are also unit connecting edges ${y_1^{(p)}y_2^{(q)} \eps_\varnothing \xrightarrow{1} y_1^{(p-1)}y_2^{(q)}\eps_1}$ and ${y_1^{(p)}y_2^{(q)} \eps_3 \xrightarrow{-1} y_1^{(p-1)}y_2^{(q)}\eps_{13}}$ for $p \ge 1$, and ${y_1^{(p)}y_2^{(q)} \eps_\varnothing \xrightarrow{1} y_1^{(p)}y_2^{(q-1)}\eps_3}$ and ${y_1^{(p)}y_2^{(q)} \eps_1 \xrightarrow{1} y_1^{(p)}y_2^{(q-1)}\eps_{13}}$ for $q \ge 1$, shown in Figure~\ref{fig:ctwo}. Along a connecting edge the multidegree does not change (Remark~\ref{rem:multigrading}). For example $y_1 y_2 \eps_\varnothing$, $y_2 \eps_1$, $y_1 \eps_3$ and $\eps_{13}$ all have multidegree $xyzw$, and the two paths ${y_1 y_2 \eps_\varnothing \xrightarrow{1} y_2 \eps_1 \xrightarrow{1} \eps_{13}}$ and ${y_1 y_2 \eps_\varnothing \xrightarrow{1} y_1 \eps_3 \xrightarrow{-1} \eps_{13}}$ stay in that multidegree.

As in Example \ref{ex:running}, the transferred system is truncated. Its map $\sigma^{\M}_{\boldsymbol{e}_1}$ is the $\sigma^{\M}$ of Example~\ref{ex:running}, and $\sigma^{\M}_{\boldsymbol{e}_2}$ is given by
    \[\sigma^{\M}_{\boldsymbol{e}_2}(\eps_1) = p^{\M}(\eps_{13}) = w \eps_{12} + x \eps_{23} ,
    \]
    so
    \[ \sigma^{\M}_{\boldsymbol{e}_2}|_0 = \begin{bmatrix} 0\\ 0\\ 1\\ 0 \end{bmatrix}, \qquad
        \sigma^{\M}_{\boldsymbol{e}_2}|_1 = \begin{bmatrix} w & 0 & 0 & 0\\ 0 & 0 & 0 & 0\\ x & z & 0 & 0\\ 0 & 0 & 0 & -w \end{bmatrix}, \qquad
        \sigma^{\M}_{\boldsymbol{e}_2}|_2 = \begin{bmatrix} 0 & -w & 0 & 0 \end{bmatrix}.
    \]

    The reduced complex over $\Rbar = R/(xy, zw)$ has ranks $1, 4, 6, 9, 11, 14, \dots$. Its connecting block is $S_k = \sum_{j=1}^{2} T^{(k)}_{\boldsymbol{e}_j} \otimes \sigma^{\M}_{\boldsymbol{e}_j}$, the block of Proposition~\ref{prop:blocks} for the base $\F^{\M}$, so that
    \[
        S_1|_q = \left[\begin{array}{c|c}
            \sigma^{\M}_{\boldsymbol{e}_1}|_q & \sigma^{\M}_{\boldsymbol{e}_2}|_q
        \end{array}\right],
        \qquad
        S_2|_q = \left[\begingroup
        \renewcommand{\arraystretch}{1.5}\begin{array}{c|c|c}
            \sigma^{\M}_{\boldsymbol{e}_1}|_q & \sigma^{\M}_{\boldsymbol{e}_2}|_q & 0\\
            \hline
            0 & \sigma^{\M}_{\boldsymbol{e}_1}|_q & \sigma^{\M}_{\boldsymbol{e}_2}|_q\\
        \end{array}\endgroup\right],
    \]
    and the matrices are $d_{\Sh}^{\lift{M}}|_1 = \begin{bmatrix} \overline{d^{\M}_1} \end{bmatrix} = \ \ \begin{NiceArray}{ccccc}
    \eps_1 & \eps_2 & \eps_3 & \eps_4 & \\
    0 & yz & 0 & xw & \eps_{\varnothing}\\
    &&&&
    \CodeAfter
  \SubMatrix[{2-1}{2-4}]\end{NiceArray}$,
    \[
        d_{\Sh}^{\lift{M}}|_2 = \left[\begin{array}{@{}c|c@{}}
            \overline{d^{\M}_2} & \overline{S_1|_0} \\
        \end{array}\right]=\ \
        \begin{NiceArray}{cccc|ccc}
        \eps_{12} & \eps_{14} & \eps_{23} & \eps_{34} & y_1\eps_{\varnothing} & y_2\eps_{\varnothing} &\\
            z & w & 0 & 0 & 1 & 0 & \eps_1\\
                -x & 0 & w & 0& 0 & 0 & \eps_2\\
                0 & 0 & -y & x & 0 & 1 & \eps_3\\
                0 & -y & 0 & -z& 0 & 0 &\eps_4,\\
                \CodeAfter
  \SubMatrix[{2-1}{5-6}]
        \end{NiceArray}
    \]

    \[
        d_{\Sh}^{\lift{M}}|_3 = \left[\begingroup
        \renewcommand{\arraystretch}{1.5}\begin{array}{c|c}
            \overline{d^{\M}_3} & \overline{S_1|_1} \\
            \hline
            0 & \overline{\mathbf{1}_2 \otimes d^{\M}_1} \\
        \end{array}\endgroup\right]=
        \ \ \begin{NiceArray}{c|cccc|ccccc}
        \eps_{134} & y_1\eps_1 & y_1\eps_2 & y_1\eps_3 & y_1\eps_4 &  y_2\eps_1 & y_2\eps_2 & y_2\eps_3 & y_2\eps_4 & \\
                w & 0 & -y & -w & 0 & w &  0 & 0 & 0 & \eps_{12}\\
                -z & 0 & 0 & 0 & -x &  0 &  0 & 0 & 0 & \eps_{14}\\
                x & 0 & 0 & -x & 0 & x &  z & 0 & 0 & \eps_{23}\\
                y & 0 & 0 & 0 & 0 & 0 &  0 & 0 & -w & \eps_{34}\\
                \hline
                0 & 0 & yz & 0 & xw & 0 &  0 & 0 & 0 & y_1\eps_{\varnothing}\\
                0 & 0 & 0 & 0 & 0  & 0 & yz & 0 & xw & y_2\eps_{\varnothing},\\
                \CodeAfter
  \SubMatrix[{2-1}{7-9}]
        \end{NiceArray}
    \]

    \[
        d_{\Sh}^{\lift{M}}|_4 =
        \left[\begingroup
        \renewcommand{\arraystretch}{1.5}\begin{array}{c|c}
            \overline{S_1|_2} & 0 \\
            \hline
            \overline{\mathbf{1}_2 \otimes d^{\M}_2}  &  \overline{S_2|_0}\\
\end{array}\endgroup\right]\]
which equals
\begin{center}\small\setlength{\arraycolsep}{2.5pt} $\begin{NiceArray}{cccc|cccc|c|c|cc}
y_1\eps_{12} & y_1 \eps_{14} & y_1 \eps_{23} & y_1 \eps_{34} & y_2\eps_{12} & y_2 \eps_{14} & y_2 \eps_{23} & y_2 \eps_{34} & y_1^{(2)}\eps_{\varnothing} & y_1y_2 \eps_{\varnothing} &y_2^{(2)}\eps_{\varnothing} & \\
                0 & 0 & 0 & x & 0 & -w & 0 & 0 & 0 & 0 & 0 & \eps_{134}\\
                \hline
                z & w & 0 & 0 & 0 & 0 & 0 & 0 & 1 & 0 & 0 & y_1\eps_1 \\
                -x & 0 & w & 0 & 0 & 0 & 0 & 0 & 0 & 0 & 0 & y_1 \eps_2\\
                0 & 0 & -y & x & 0 & 0 & 0 & 0 & 0 & 1 & 0 & y_1 \eps_3\\
                0 & -y & 0 & -z & 0 & 0 & 0 & 0 & 0 & 0 & 0 & y_1 \eps_4\\
                \hline
                0 & 0 & 0 & 0 & z & w & 0 & 0 & 0 & 1 & 0 & y_2 \eps_1\\
                0 & 0 & 0 & 0 & -x & 0 & w & 0 & 0 & 0 & 0 & y_2 \eps_2\\
                0 & 0 & 0 & 0 & 0 & 0 & -y & x & 0 & 0 & 1 & y_2 \eps_3\\
                0 & 0 & 0 & 0 & 0 & -y & 0 & -z & 0 & 0 & 0 & y_2 \eps_4.\\
                \CodeAfter
  \SubMatrix[{2-1}{10-11}]
        \end{NiceArray}$\end{center}

    The unit entries are the $1$ from $y^{\boldsymbol{u}} \eps_\varnothing$ to $y^{\boldsymbol{u} - \boldsymbol{e}_1} \eps_1$ whenever $u_1 \ge 1$ and to $y^{\boldsymbol{u} - \boldsymbol{e}_2} \eps_3$ whenever $u_2 \ge 1$, the two unit edges ${y_1 \eps_\varnothing \xrightarrow{1} \eps_1}$ and ${y_2 \eps_\varnothing \xrightarrow{1} \eps_3}$ repeated over every $\boldsymbol{u}$.
    Matching each $y^{\boldsymbol{u}} \eps_\varnothing$ with $|\boldsymbol{u}| \ge 1$ to one of them leaves no unit entry, so
    \[
        \beta^{\Rbar}_n(R/I) = 1, 2, 4, 6, 8, 10, \dots, \quad \text{that is, } \beta^{\Rbar}_n(R/I) = 2n \text{ for } n \ge 1, \qquad \text{for } \Rbar = R/(xy, zw).
    \]
\end{example}

\section{Algebraic Morse theory on the Shamash construction is again a Shamash construction}\label{sec:commute}
Throughout this section we fix a Morse matching $\M$ on $\F$, and split its differential as $d = \mu + \delta$. Let $(i, p, h)$ and $(i^{\M}, p^{\M}, h^{\M})$ be the contractions of Theorem~\ref{thm:amt}. Here $i$ and $p$ are the inclusion of and the projection onto the critical cells, and $h$ is minus the inverse of $\mu$ on the lower cells. The contraction of $\F$ onto $\F^{\M}$ is then $(i^{\M}, p^{\M}, h^{\M}) = (\Phi i, p + p\delta\Phi h, \Phi h)$, where $\Phi = \sum_{n \ge 0} (h\delta)^n$.

We state the main theorem now to set up the hypotheses for the section, and postpone its proof until the end of the section.

\begin{theorem}\label{thm:commute}
    Let $\F$ be a free resolution of $R/I$ over $R$ with a finite basis of each $F_q$, let $a_1, \dots, a_c$ be an $R$-regular sequence in $I$, let $\sigma = \{\sigma^{(j)}\}_{j=1}^c$ be a truncated system of higher homotopies for $a_1, \dots, a_c$ on $\F$, and let $\M$ be a Morse matching on $\F$.
    \begin{enumerate}[label={\rm(\roman*)},leftmargin=*,itemsep=4pt]
        \item Under the identity map on the common basis $\{y^{\boldsymbol{u}} \otimes \bar\kappa\}$, with $\boldsymbol{u} \in \N^c$ and $\kappa$ critical for $\M$, $\Sh(\F, \sigma)^{\lift{M}}$ is the free complex $R\{y_j\}_{j=1}^{c} \otimes_R \F^{\M} \otimes_R \Rbar$ over $\Rbar$ with differential
        \[
             d_{\Sh}^{\lift{M}}=\sum_{\boldsymbol{u} \in \N^c} t^{\boldsymbol{u}} \otimes \sigma^{\M}_{\boldsymbol{u}} \otimes \id .
        \]
        \item The maps $\sigma^{\M}$ form a system of higher homotopies for $a_1, \dots, a_c$ on $\F^{\M}$ over $R$, that is, $\sigma^{\M}_0 = d^{\M}$ and
        \[
            \sum_{\boldsymbol{v} + \boldsymbol{v'} = \boldsymbol{u}} \sigma^{\M}_{\boldsymbol{v}} \sigma^{\M}_{\boldsymbol{v'}} =
            \begin{cases}
                a_j \id & \text{if } \boldsymbol{u} = \boldsymbol{e}_j,\\
                0 & \text{otherwise.}
            \end{cases}
        \]
        So the complex in (i) is the Shamash resolution $\Sh(\F^{\M}, \sigma^{\M})$ of Definition~\ref{def:shamash}.
    \end{enumerate}
\end{theorem}

The following corollary can be used to compute the ranks of the reduced Eisenbud--Shamash complex; in fact it was applied in Examples \ref{ex:running} and \ref{ex:runningtwo}.

\begin{corollary}\label{cor:ranks}
Let $\mathcal{K}_q$ be the set of cells of homological degree $q$ that are critical for $\M$, the basis of the degree $q$ part of $K_\bullet$ in Theorem~\ref{thm:amt}, with $\mathcal{K}_q = \emptyset$ for $q < 0$. Then
    \[
        \operatorname{rank} \Sh(\F^{\M}, \sigma^{\M})_n = \sum_{\boldsymbol{u} \in \N^c} \#\mathcal{K}_{n - 2|\boldsymbol{u}|},
    \]
    and this is a free $\Rbar$-resolution of $R/I$.
\end{corollary}

When $I$ and $a_1, \dots, a_c$ are homogeneous, every $a_j$ lies in $\mathfrak{m} I$ and $\F^{\M}$ is the minimal free resolution of $R/I$ over $R$, the resolution $\Sh(\F^{\M}, \sigma^{\M})$ is minimal, by Remark~\ref{rem:mincond} applied to the system $\sigma^{\M}$ of Theorem~\ref{thm:commute}(ii).

For an $R$-linear map $f$ between $\F$ and $\F^{\M}$ we write $\lift{f} = \id \otimes f$, so $\lift{f}(y^{\boldsymbol{u}} \otimes \alpha) = y^{\boldsymbol{u}} \otimes f(\alpha)$, with no sign, and $\lift{f} \, \lift{g} = \id \otimes fg$. These are maps over $R$, and the script letter marks the copy over each $\boldsymbol{u}$, as in the lifted matching $\lift{M}$. On $\Sh(\F, \sigma)$ and $\Sh(\F^{\M}, \sigma^{\M})$ we use their reductions $\overline{\lift{f}}$. For the Greek maps $\delta$ and $\mu$ we write $\id \otimes \delta$ and $\id \otimes \mu$. We use $\lift{d}$, $\id \otimes \delta$, $\lift{d}^{\M}$, $(\lift{i}, \lift{p}, \lift{h})$ and $(\lift{i}^{\M}, \lift{p}^{\M}, \lift{h}^{\M})$, and $\lift{M}$ is the lifted matching.

Because the input is truncated, the transfer $\sigma^{\M}_{\boldsymbol{u}} = \sum p^{\M} \sigma^{(j_1)} h^{\M} \cdots h^{\M} \sigma^{(j_k)} i^{\M}$ of Definition~\ref{def:transferred} involves only the $c$ maps $\sigma^{(j)}$, and in closed form it is the perturbation formula $d^{\M} = p \delta \Phi i$ of Theorem~\ref{thm:amt} with the connecting part $\dcon$ in place of $\delta$.

Here $\lift{d} = \partial_0$ and $\dcon = \partial_1 + \dots + \partial_c$ is the connecting part of Definition~\ref{def:delcon}, so $D = \lift{d} + \dcon$, $d_{\Sh} = \overline{\lift{d}} + \bar\partial_{\mathrm{con}}$, and $\dcon(y^{\boldsymbol{u}} \otimes \alpha) = \sum_{j : u_j \ge 1} y^{\boldsymbol{u} - \boldsymbol{e}_j} \otimes \sigma^{(j)}(\alpha)$.

\begin{proposition}\label{prop:closed}
     In the setting of Theorem~\ref{thm:commute}, as operators from $R\{y_j\}_{j=1}^{c} \otimes_R \F^{\M}$ to itself,
    \[
        \sum_{\boldsymbol{u} \in \N^c} t^{\boldsymbol{u}} \otimes \sigma^{\M}_{\boldsymbol{u}} = \lift{d}^{\M} + \lift{p}^{\M} \dcon \Phi \lift{i}^{\M} ,
    \]
    where $\Phi = \sum_{m \ge 0} (\lift{h}^{\M} \dcon)^m$ is a finite sum on each basis element.
\end{proposition}

Every map in Proposition~\ref{prop:closed} is an $R$-linear map.
The left side of the equality is the map $D$ of the system of higher homotopies $\sigma^{\M}$, and its reduction modulo $\aaa$ is the differential of $\Sh(\F^{\M}, \sigma^{\M})$. The identity is used over $\Rbar$ only after reduction, in the proof of Theorem~\ref{thm:commute}.

\begin{proof}
    The term of the series with $m = k - 1$ contributes
     \[
        \lift{p}^{\M} \dcon \bigl( \lift{h}^{\M} \dcon \bigr)^{k-1} \lift{i}^{\M} = \sum_{(j_1, \dots, j_k)} \bigl( t_{j_1} \cdots t_{j_k} \bigr) \otimes \bigl( p^{\M} \sigma^{(j_1)} h^{\M} \cdots h^{\M} \sigma^{(j_k)} i^{\M} \bigr),
    \]
    the sum being over all sequences of colors $(j_1, \dots, j_k)$ of length $k$.
    To
     observe this, substitute $\dcon = \sum_j t_j \otimes \sigma^{(j)}$ into each of the $k$ factors
    and expand, using $$(t \otimes f)(\id \otimes g) = t \otimes (fg), (\id \otimes g)(t \otimes f) = t \otimes (gf) \text{ and } (t \otimes f)(t' \otimes f') = (tt') \otimes (ff').$$ None of these carries a sign, because every $t_j$ has even homological degree and $\id$ has homological degree $0$.

    Collecting the sequences with $\boldsymbol{e}_{j_1} + \dots + \boldsymbol{e}_{j_k} = \boldsymbol{u}$ and using $t_{j_1} \cdots t_{j_k} = t^{\boldsymbol{u}}$ gives $\sum_{|\boldsymbol{u}| = k} t^{\boldsymbol{u}} \otimes \sigma^{\M}_{\boldsymbol{u}}$ by \eqref{eq:transfer} of Definition \ref{def:transferred}. Summing over $k \ge 1$ and adding the term $\lift{d}^{\M} = \id \otimes \sigma^{\M}_0$ gives the display.
    The series $\Phi$ is finite on $y^{\boldsymbol{u}} \otimes \kappa$ because $\lift{h}^{\M} \dcon$ lowers $|\boldsymbol{u}|$ by one.
\end{proof}

So $\sigma^{\M}$ comes from $\dcon$ exactly as $d^{\M} = p \delta \Phi i$ comes from $\delta$. Proposition~\ref{prop:closed} uses only the degrees of $i^{\M}$, $p^{\M}$, $h^{\M}$, so it holds for any such triple for which the series is finite, with $\sigma^{\M}_{\boldsymbol{u}}$ defined by Definition \ref{def:transferred} \eqref{eq:transfer} for that triple. For one color it is a single geometric series.

\begin{corollary}\label{cor:conecolor}
     In the setting of Theorem~\ref{thm:commute} and Proposition \ref{prop:closed}, with $c = 1$,
     identify $\N^1$ with $\N$, and write $a = a_1$, $\sigma = \sigma^{(1)}$ and $\sigma^{\M}_k = \sigma^{\M}_{k \boldsymbol{e}_1}$. Then
    \[
        \sigma^{\M}_k = p^{\M} \sigma \bigl( h^{\M} \sigma \bigr)^{k-1} i^{\M} \quad (k \ge 1), \qquad \sum_{k \ge 1} \sigma^{\M}_k = p^{\M} \sigma \Phi i^{\M}, \quad \Phi = \sum_{m \ge 0} (h^{\M} \sigma)^m .
    \]
     So $\sigma^{\M}$ is truncated if and only if $p^{\M} \sigma (h^{\M} \sigma)^{k} i^{\M} = 0$ for every $k \ge 1$, and $\sigma h^{\M} \sigma = 0$ suffices.
\end{corollary}

The next lemma says the reduced resolution keeps the identity $D^2 = \Theta$,
with $\Theta = \sum_j t_j \otimes a_j \id$, of Remark~\ref{rem:piecessquare}. Over $\Rbar$ this would follow from the reduction alone, since $\Theta$ vanishes there, but over $R$ the reduction must reproduce $\Theta$. The idea is that the inclusion, perturbed along with the differential, still intertwines $D$ and $D^{\M}$, and $\Theta$, which commutes with all the maps involved, passes through it.

\begin{lemma}\label{lem:curved}
     In the setting of Theorem~\ref{thm:commute}, as operators from $R\{y_j\}_{j=1}^{c} \otimes_R \F^{\M}$ to itself, before the tensor with $\Rbar$,
    \[
        \Bigl( \sum_{\boldsymbol{u} \in \N^c} t^{\boldsymbol{u}} \otimes \sigma^{\M}_{\boldsymbol{u}} \Bigr) \Bigl( \sum_{\boldsymbol{u} \in \N^c} t^{\boldsymbol{u}} \otimes \sigma^{\M}_{\boldsymbol{u}} \Bigr) = \sum_{j=1}^{c} t_j \otimes a_j \id .
    \]
\end{lemma}

\begin{proof}
     Recall that $D = \lift{d} + \dcon$, so that $d_{\Sh} = \overline{D}$, and $\Theta = \sum_j t_j \otimes a_j \id$. By Proposition~\ref{prop:closed} the corresponding map for $\sigma^{\M}$, whose reduction is the differential in Theorem~\ref{thm:commute}(i), is
    \[
        D^{\M} = \sum_{\boldsymbol{u} \in \N^c} t^{\boldsymbol{u}} \otimes \sigma^{\M}_{\boldsymbol{u}} = \lift{d}^{\M} + \lift{p}^{\M} \dcon \Phi \lift{i}^{\M}, \qquad \Phi = \sum_{m \ge 0} (\lift{h}^{\M} \dcon)^m ,
    \]
    and we must show $(D^{\M})^2 = \Theta$. We use three facts.
    \begin{enumerate}[label=(\arabic*),leftmargin=*,itemsep=1pt]
        \item $D^2 = \Theta$ by Remark~\ref{rem:piecessquare}. Since $\lift{d}^2 = 0$, this says $\lift{d} \dcon + \dcon \lift{d} + \dcon \dcon = \Theta$.
        \item $\Theta$ commutes with $\dcon$ and with every map $\lift{f} = \id \otimes f$, since the CI operators commute with each other and with every $\lift{f}$, and each $a_j$ lies in $R$.
        \item The contraction identities of $(i^{\M}, p^{\M}, h^{\M})$ hold
         for the maps $\lift{i}^{\M}$, $\lift{p}^{\M}$, $\lift{h}^{\M}$ between $R\{y_j\}_{j=1}^{c} \otimes_R \F$ and $R\{y_j\}_{j=1}^{c} \otimes_R \F^{\M}$ (these are $R$-linear maps,
        not attached to the lifted matching, which lives on $\Sh(\F, \sigma)$), since $\lift{f} \, \lift{g} = \id \otimes fg$:
        \[
            \lift{p}^{\M} \lift{i}^{\M} = \id, \qquad \lift{d} \lift{h}^{\M} = \lift{i}^{\M} \lift{p}^{\M} - \id - \lift{h}^{\M} \lift{d}, \qquad \lift{h}^{\M} \lift{h}^{\M} = \lift{h}^{\M} \lift{i}^{\M} = \lift{p}^{\M} \lift{h}^{\M} = 0, \qquad \lift{d} \lift{i}^{\M} = \lift{i}^{\M} \lift{d}^{\M} .
        \]
    \end{enumerate}
    The series $\Phi$ is finite on each cell, since $\dcon$ lowers $|\boldsymbol{u}|$ by one and $\lift{h}^{\M}$ preserves it, and $\Phi$ inverts $\id - \lift{h}^{\M} \dcon$.

    We perturb the inclusion along with the differential:
    set $\eta = \Phi \lift{i}^{\M}$.
    Expanding the series $\Phi$ and factoring $\lift{h}^{\M}\dcon$ on the left from all terms with $m \geq 1$, we see that $$\eta=\Phi\lift{i}^{\M}=\lift{i}^{\M} + \lift{h}^{\M}\dcon\Phi\lift{i}^{\M}=\lift{i}^{\M}+\lift{h}^{\M}\dcon \eta.$$

    Replacing $\Phi \lift{i}^{\M}$ with $\eta$ in $D ^{\M}$, we obtain $D^{\M} = \lift{d}^{\M} + \lift{p}^{\M} \dcon \eta$. We show
    \[
        \lift{p}^{\M} \eta = \id \qquad \text{and} \qquad D \eta = \eta D^{\M} .
    \]
    Granting these, $\eta (D^{\M})^2 = D \eta D^{\M} = D^2 \eta = \Theta \eta = \eta \Theta$ by (1) and (2), and applying $\lift{p}^{\M}$ gives $(D^{\M})^2 = \Theta$.

    $\lift{p}^{\M}$ is a left inverse of $\eta$: By (3), $\lift{p}^{\M} \eta = \lift{p}^{\M} \lift{i}^{\M} + \lift{p}^{\M} \lift{h}^{\M} \dcon \eta = \id$.

    $\eta$ intertwines $D$ and $D^{\M}$: Let $E = D \eta - \eta D^{\M}$.
    Replacing $\eta$ by $\lift{i}^{\M}+\lift{h}^{\M}\dcon\eta$
    and using (3) for $\lift{d} \lift{i}^{\M}$ and $\lift{d} \lift{h}^{\M}$,
    \[
        \lift{d} \eta = \lift{d} \lift{i}^{\M} + \lift{d} \lift{h}^{\M} \dcon \eta = \lift{i}^{\M} \lift{d}^{\M} + \lift{i}^{\M} \lift{p}^{\M} \dcon \eta - \dcon \eta - \lift{h}^{\M} \lift{d} \dcon \eta ,
    \]
    so $D \eta = \lift{d} \eta + \dcon \eta = \lift{i}^{\M} D^{\M} - \lift{h}^{\M} \lift{d} \dcon \eta$.

    Using the replacement $\eta=\lift{i}^{\M}+\lift{h}^{\M}\dcon\eta$, we obtain
    $\eta D^{\M} = \lift{i}^{\M} D^{\M} + \lift{h}^{\M} \dcon \eta D^{\M}$. Subtracting,
    \[
        E = -\lift{h}^{\M} \bigl( \lift{d} \dcon \eta + \dcon \eta D^{\M} \bigr) = -\lift{h}^{\M} \bigl( \lift{d} \dcon \eta + \dcon D \eta \bigr) + \lift{h}^{\M} \dcon E ,
    \]
    using $\eta D^{\M} = D \eta - E$. By (1), $\lift{d} \dcon \eta + \dcon D \eta = (\lift{d} \dcon + \dcon \lift{d} + \dcon \dcon) \eta = \Theta \eta$, and by (2) and (3), $\lift{h}^{\M} \Theta \eta = \Theta \lift{h}^{\M} (\lift{i}^{\M} + \lift{h}^{\M} \dcon \eta) = 0$. So $(\id - \lift{h}^{\M} \dcon) E = 0$, and applying $\Phi$ gives $E = 0$.
\end{proof}

The following lemma is a version of \cite[Proposition~3.6]{Real2000} under different hypotheses. Real assumes that each of $\delta_1$, $\delta_2$ and $\delta_1 + \delta_2$ is a perturbation of $(N_\bullet, d)$ whose composite with $h$ is pointwise nilpotent, and concludes that the two perturbed contractions agree. We do not need $\delta_2$ to be a perturbation of $(N_\bullet, d)$ on its own, which need not hold in our application, and we only compare the perturbed differentials.
\begin{lemma}[Perturbation Additivity Lemma]\label{lem:additive}
    Let $(i, p, h)$ be a contraction of a complex $(N_\bullet, d)$ onto a complex $(N_\bullet', d')$ of $A$-modules, and let $\delta_1, \delta_2 \colon N_\bullet \to N_{\bullet - 1}$ be $A$-linear maps such that $\delta_1$ and $\delta_1 + \delta_2$ are perturbations of $(N_\bullet, d)$. Suppose that $$\Phi_1 = \sum_{n \ge 0} (h\delta_1)^n = (\id - h\delta_1)^{-1} \text{ and } \Phi_{12} = \sum_{n \ge 0} \bigl( h(\delta_1 + \delta_2) \bigr)^n = \bigl(\id - h(\delta_1 + \delta_2)\bigr)^{-1}$$ are finite sums on each element,
    let
    \[
        (i_1, p_1, h_1) = \bigl( \Phi_1 i, \; p + p \delta_1 \Phi_1 h, \; \Phi_1 h \bigr)
    \]
    be the contraction of $(N_\bullet, d + \delta_1)$ onto $(N_\bullet', d' + p \delta_1 \Phi_1 i)$ of Lemma~\ref{lem:perturbation}, and suppose that $$\Phi_2 = \sum_{n \ge 0} (h_1\delta_2)^n = (\id - h_1\delta_2)^{-1}$$ is a finite sum on each element. Then perturbing by $\delta_1$ and then by $\delta_2$ gives the same differential on $N_\bullet'$ as perturbing by $\delta_1 + \delta_2$ at once:
    \[
        d' + p (\delta_1 + \delta_2) \Phi_{12} i
        = d' + p \delta_1 \Phi_1 i + p_1 \delta_2 \Phi_2 i_1 .
    \]
\end{lemma}

\begin{figure}[htbp]
    \centering
        \begin{tikzpicture}[>=Stealth, x=1cm, y=1cm, font=\small]
            \node (L1) at (0, 6.3) {$(N_\bullet, d)$};
            \node (R1) at (6.4, 6.3) {$(N'_\bullet, d')$};
            \node (L2) at (0, 4.2) {$(N_\bullet, d + \delta_1)$};
            \node (R2) at (6.4, 4.2) {$\bigl(N_\bullet', d' + p \delta_1 \Phi_1 i\bigr)$};
            \node (L3) at (0, 2.1) {$(N_\bullet, d + \delta_1 + \delta_2)$};
            \node[align=center] (R3) at (6.4, 2.1) {$\bigl(N'_\bullet, d' + p \delta_1 \Phi_1 i$\\[1pt] $+ p_1 \delta_2 \Phi_2 i_1\bigr)$};
            \node (L4) at (0, 0) {$(N_\bullet, d + \delta_1 + \delta_2)$};
            \node (R4) at (6.4, 0) {$\bigl(N'_\bullet, d' + p (\delta_1 + \delta_2) \Phi_{12} i\bigr)$};
            \foreach \r/\p/\i in {1/{p}/{i}, 2/{p_1}/{i_1}, 3/{p_2}/{i_2}, 4/{p_{12}}/{i_{12}}}{
            \draw[->] ([yshift=3.5pt]L\r.east) -- node[above, font=\footnotesize]{$\p$} ([yshift=3.5pt]R\r.west);
            \draw[<-] ([yshift=-3.5pt]L\r.east) -- node[below, font=\footnotesize]{$\i$} ([yshift=-3.5pt]R\r.west);}
            \foreach \r/\h in {1/{h}, 2/{h_1}, 3/{h_2}, 4/{h_{12}}}{
            \draw[->] ([yshift=3pt]L\r.west) .. controls +(-0.9, 0.4) and +(-0.9, -0.4) .. node[left, font=\footnotesize]{$\h$} ([yshift=-3pt]L\r.west);}
            \draw[->, thick] (L1) -- node[right, font=\footnotesize]{$+\delta_1$} (L2);
            \draw[->, thick] (L2) -- node[right, font=\footnotesize]{$+\delta_2$} (L3);
            \draw[->, thick] (-2.6, 6.1) .. controls (-4.0, 5.0) and (-4.0, 1.1) .. node[left, font=\footnotesize, align=right]{$+(\delta_1 + \delta_2)$\\ at once} (-2.6, 0.2);
            \foreach \s in {-1.2pt, 1.2pt}{
            \draw ([xshift=\s, yshift=-1pt]L3.south) -- ([xshift=\s, yshift=1pt]L4.north);
            \draw ([xshift=\s, yshift=-1pt]R3.south) -- ([xshift=\s, yshift=1pt]R4.north);}
            \node[right, font=\footnotesize] at ($(R3.south)!0.5!(R4.north) + (0.1, 0)$) {Lemma~\ref{lem:additive}};
        \end{tikzpicture}
    \caption{Perturbing by $\delta_1$ and then by $\delta_2$, or by $\delta_1 + \delta_2$ at once.}
    \label{fig:additive}
\end{figure}

\begin{proof}
    Figure~\ref{fig:additive} shows the two ways of perturbing by $\delta_1 + \delta_2$. In Figure~\ref{fig:additive} $(i_2, p_2, h_2)$ is the contraction of Lemma~\ref{lem:perturbation} for $(i_1, p_1, h_1)$ and $\delta_2$, which is a perturbation of $(N_\bullet, d + \delta_1)$ since $(d + \delta_1 + \delta_2)^2 = 0$, and $(i_{12}, p_{12}, h_{12})$ the one for $(i, p, h)$ and $\delta_1 + \delta_2$. So every row is a contraction, and it remains to compare the two perturbation terms.

    Since $(\id - h\delta_1)\Phi_1 = \id$ and $h_1 = \Phi_1 h$,
    \[
        \id - h(\delta_1 + \delta_2) = (\id - h\delta_1)(\id - \Phi_1 h \delta_2) = (\id - h\delta_1)(\id - h_1 \delta_2),
    \]
    so $\Phi_{12} = \Phi_2 \Phi_1$. The left side and right side respectively are
    $$p \delta_1 \Phi_2 \Phi_1 i + p \delta_2 \Phi_2 \Phi_1 i \text{ and }
    p \delta_1 \Phi_1 i + p \delta_2 \Phi_2 \Phi_1 i + p \delta_1 \Phi_1 h \delta_2 \Phi_2 \Phi_1 i.$$ Their difference is
    \[
        p \delta_1 (\Phi_2 - \id) \Phi_1 i - p \delta_1 \Phi_1 h \delta_2 \Phi_2 \Phi_1 i ,
    \]
    which is zero because $\Phi_2 - \id = h_1 \delta_2 \Phi_2 = \Phi_1 h \delta_2 \Phi_2$.
\end{proof}

The Perturbation Additivity Lemma is the main tool for the proof of the main theorem, and we are now ready to begin.

\begin{proof}[Proof of Theorem~\ref{thm:commute}]
    \textbf{Part (i):} The edges of $\lift{M}$ lie over those of $\M$, so the matched part of $d_{\Sh}$ is $\id \otimes \bar\mu$ and the unmatched part is $\id \otimes \bar\delta + \bar\partial_{\mathrm{con}}$. Theorem~\ref{thm:amt} attaches to $\lift{M}$ the contraction $(\overline{\lift{i}}, \overline{\lift{p}}, \overline{\lift{h}})$, and since $\lift{M}$ is a Morse matching (Proposition~\ref{prop:lift}), the theorem applies and gives
    \[
        d_{\Sh}^{\lift{M}} = \overline{\lift{p}} \, (\id \otimes \bar\delta + \bar\partial_{\mathrm{con}}) \, \overline{\Phi}_{12} \, \overline{\lift{i}}, \qquad \Phi_{12} = \sum_{n \ge 0} \bigl( \lift{h} (\id \otimes \delta + \dcon) \bigr)^n ,
    \]
    where we use $\bar f \, \bar g = \overline{fg}$ for $R$-linear maps $f$ and $g$.

    In other words, $\Sh(\F, \sigma)^{\lift{M}}$ comes from perturbing $\id \otimes \bar\mu$ by $\id \otimes \bar\delta + \bar\partial_{\mathrm{con}}$ all at once. We show that perturbing first by $\id \otimes \bar\delta$ and then by $\bar\partial_{\mathrm{con}}$ gives the same result, as in Figure~\ref{fig:M2ways}. This is exactly the situation of Lemma~\ref{lem:additive}: Figure~\ref{fig:M2ways} is Figure~\ref{fig:additive} for the contraction $(\overline{\lift{i}}, \overline{\lift{p}}, \overline{\lift{h}})$, with $d = \id \otimes \bar\mu$, $d' = 0$, $\delta_1 = \id \otimes \bar\delta$ and $\delta_2 = \bar\partial_{\mathrm{con}}$. In particular $(i_2, p_2, h_2)$ is the contraction obtained by perturbing by $\delta_1$ and then by $\delta_2$.

    \begin{figure}[htbp]
        \centering
        \resizebox{\linewidth}{!}{%
        \begin{tikzpicture}[>=Stealth, x=1cm, y=1cm, font=\small]
            \node (L1) at (0, 6.9) {$\bigl( R\{y_j\}_{j=1}^{c} \otimes_R \F \otimes_R \Rbar, \id \otimes \bar\mu \bigr)$};
            \node (R1) at (7.2, 6.9) {$\bigl( R\{y_j\}_{j=1}^{c} \otimes_R \F^{\M} \otimes_R \Rbar, 0 \bigr)$};
            \node (L2) at (0, 4.6) {$\bigl( R\{y_j\}_{j=1}^{c} \otimes_R \F \otimes_R \Rbar, \overline{\lift{d}} \bigr)$};
            \node (R2) at (7.2, 4.6) {$\bigl( R\{y_j\}_{j=1}^{c} \otimes_R \F^{\M} \otimes_R \Rbar, \overline{\lift{d}}^{\M} \bigr)$};
            \node (L3) at (0, 2.3) {$\bigl( R\{y_j\}_{j=1}^{c} \otimes_R \F \otimes_R \Rbar, d_{\Sh} \bigr)$};
            \node[align=center] (R3) at (7.2, 2.3) {$\bigl( R\{y_j\}_{j=1}^{c} \otimes_R \F^{\M} \otimes_R \Rbar, \sum_{\boldsymbol{u}} t^{\boldsymbol{u}} \otimes \bar\sigma^{\M}_{\boldsymbol{u}} \bigr)$\\[1pt] $= \Sh(\F^{\M}, \sigma^{\M})$};
            \node (L4) at (0, 0) {$\bigl( R\{y_j\}_{j=1}^{c} \otimes_R \F \otimes_R \Rbar, d_{\Sh} \bigr)$};
            \node[align=center] (R4) at (7.2, 0) {$\bigl( R\{y_j\}_{j=1}^{c} \otimes_R \F^{\M} \otimes_R \Rbar, d^{\lift{M}} \bigr)$\\[1pt] $= \Sh(\F, \sigma)^{\lift{M}}$};
            \foreach \r/\p/\i in {1/{\overline{\lift{p}}}/{\overline{\lift{i}}}, 2/{\overline{\lift{p}}^{\M}}/{\overline{\lift{i}}^{\M}}, 3/{p_2}/{i_2}, 4/{p^{\lift{M}}}/{i^{\lift{M}}}}{
            \draw[->] ([yshift=3.5pt]L\r.east) -- node[above, font=\footnotesize]{$\p$} ([yshift=3.5pt]R\r.west);
            \draw[<-] ([yshift=-3.5pt]L\r.east) -- node[below, font=\footnotesize]{$\i$} ([yshift=-3.5pt]R\r.west);}
            \foreach \r/\h in {1/{\overline{\lift{h}}}, 2/{\overline{\lift{h}}^{\M}}, 3/{h_2}, 4/{h^{\lift{M}}}}{
            \draw[->] ([yshift=3pt]L\r.west) .. controls +(-0.9, 0.4) and +(-0.9, -0.4) .. node[left, font=\footnotesize]{$\h$} ([yshift=-3pt]L\r.west);}
            \draw[->, thick] (L1) -- node[right, font=\footnotesize]{$+\id \otimes \bar\delta$} (L2);
            \draw[->, thick] (L2) -- node[right, font=\footnotesize]{$+\bar\partial_{\mathrm{con}}$} (L3);
            \draw[->, thick] (-4.0, 6.7) .. controls (-5.2, 5.6) and (-5.2, 1.3) .. node[left, font=\footnotesize, align=right]{$+(\id \otimes \bar\delta + \bar\partial_{\mathrm{con}})$\\ at once} (-4.0, 0.2);
            \foreach \s in {-1.2pt, 1.2pt}{
            \draw ([xshift=\s, yshift=-1pt]L3.south) -- ([xshift=\s, yshift=1pt]L4.north);
            \draw ([xshift=\s, yshift=-1pt]R3.south) -- ([xshift=\s, yshift=1pt]R4.north);}
            \node[right, font=\footnotesize] at ($(R3.south)!0.5!(R4.north) + (0.1, 0)$) {Theorem~\ref{thm:commute}};
        \end{tikzpicture}}
        \caption{Matching done in two ways.}
        \label{fig:M2ways}
    \end{figure}

    The maps in the top two rows of Figure~\ref{fig:M2ways} are of the form $\overline{\lift{f}}$ for maps $f$ on $\F$, since $\lift{M}$ and $\id \otimes \bar\delta$ are copies of $\M$ and $\bar\delta$ over each $\boldsymbol{u}$. This fails in the bottom two rows, because $\dcon$ lowers $|\boldsymbol{u}|$ and is not of the form $\lift{f}$. The proof of Lemma~\ref{lem:additive} also gives $(i_2, p_2, h_2) = (i^{\lift{M}}, p^{\lift{M}}, h^{\lift{M}})$, although we do not need this.

    Lemma~\ref{lem:additive} says that perturbing by $\delta_1$ and then by $\delta_2$ gives the same differential, $$p \delta_1 \Phi_1 i + p_1 \delta_2 \Phi_2 i_1,$$ as perturbing by $\delta_1 + \delta_2$ at once. We apply it over $\Rbar$ to the contraction $(\overline{\lift{i}}, \overline{\lift{p}}, \overline{\lift{h}})$ of $$\bigl(R\{y_j\}_{j=1}^{c} \otimes_R \F \otimes_R \Rbar, \id \otimes \bar\mu\bigr) \text{ onto } \bigl(R\{y_j\}_{j=1}^{c} \otimes_R K_\bullet \otimes_R \Rbar, 0\bigr),$$ with $\delta_1 = \id \otimes \bar\delta$ and $\delta_2 = \bar\partial_{\mathrm{con}}$. Both $\delta_1$ and $\delta_1 + \delta_2$ are perturbations, since
    \[
        (\id \otimes \bar\mu + \id \otimes \bar\delta)^2 = \overline{\lift{d}}^{\,2} = 0 \quad \text{and} \quad (\id \otimes \bar\mu + \id \otimes \bar\delta + \bar\partial_{\mathrm{con}})^2 = \overline{D^2} = d_{\Sh}^2 = 0 ,
    \]
    the second because every $a_j$ is zero in $\Rbar$. The series in the lemma are $\overline{\Phi}_{12}$, $\overline{\Phi}_1$ and $\overline{\Phi}_2$, with $\Phi_{12}$ as above, $\Phi_1 = \sum_{n \ge 0} \bigl(\lift{h} (\id \otimes \delta)\bigr)^n$ and $\Phi_2 = \sum_{n \ge 0} (\lift{h}^{\M} \dcon)^n$. Each is finite on each cell: $\Phi_{12}$ as above, $\Phi_1$ by Theorem~\ref{thm:amt} for $\M$, and $\Phi_2$
    by Proposition~\ref{prop:closed}, since $\dcon$ lowers $|\boldsymbol{u}|$. Perturbing by $\delta_1$ alone is just Theorem~\ref{thm:amt} for $\M$, applied in each copy $y^{\boldsymbol{u}} \otimes \F$. That theorem gives $d^{\M} = p \delta \Phi i$ and $(i^{\M}, p^{\M}, h^{\M}) = (\Phi i, p + p \delta \Phi h, \Phi h)$ with $\Phi = \sum_{n \ge 0}(h\delta)^n$, so lifting and then reducing gives
    \[
        \lift{p} \, (\id \otimes \delta) \, \Phi_1 \, \lift{i} = \lift{d}^{\M}, \qquad
        (i_1, p_1, h_1) = (\overline{\lift{i}}^{\M}, \overline{\lift{p}}^{\M}, \overline{\lift{h}}^{\M}) .
    \]
    Now Lemma~\ref{lem:additive} and Proposition~\ref{prop:closed}, tensored with $\Rbar$,
    give
    \[
        d_{\Sh}^{\lift{M}} = \overline{\lift{d}}^{\M} + \overline{\lift{p}}^{\M} \, \bar\partial_{\mathrm{con}} \, \overline{\Phi}_2 \, \overline{\lift{i}}^{\M}
        = \sum_{\boldsymbol{u} \in \N^c} t^{\boldsymbol{u}} \otimes \bar\sigma^{\M}_{\boldsymbol{u}} ,
    \]
    which is the Shamash differential of the system $\sigma^{\M}$ on the cells $y^{\boldsymbol{u}} \otimes \bar\kappa$ critical for $\lift{M}$.

    \textbf{Part (ii):} Expanding the square in Lemma~\ref{lem:curved} gives
    \[
        \sum_{\boldsymbol{u}} t^{\boldsymbol{u}} \otimes \Bigl( \sum_{\boldsymbol{u_1} + \boldsymbol{u_2} = \boldsymbol{u}} \sigma^{\M}_{\boldsymbol{u_1}} \sigma^{\M}_{\boldsymbol{u_2}} \Bigr) = \sum_{j} t_j \otimes a_j \id .
    \]
    On $y^{\boldsymbol{v}} \otimes \kappa$ the terms for distinct $\boldsymbol{u} \le \boldsymbol{v}$ lie in distinct summands $y^{\boldsymbol{v}-\boldsymbol{u}} \otimes \F^{\M}$, so we can compare the coefficients of each $t^{\boldsymbol{u}}$, and $$\sum\limits_{\boldsymbol{u_1} + \boldsymbol{u_2} = \boldsymbol{u}} \sigma^{\M}_{\boldsymbol{u_1}} \sigma^{\M}_{\boldsymbol{u_2}}=\begin{cases} a_j \id \text{ when } \boldsymbol{u} = \boldsymbol{e}_j\\ 0 \text{ otherwise.} \\
    \end{cases}$$
\end{proof}

\begin{remark}\label{rem:unitconnecting}
    By Theorem~\ref{thm:commute}, $\Sh(\F^{\M}, \sigma^{\M})$ is a free $\Rbar$-resolution of $R/I$. Its level edges have the weights of the edges of $G(\F^{\M})$, and its color-$j$ connecting edges have the weights of the entries of $\overline{(\sigma^{\M})^{(j)}}$. So when a connecting edge has unit weight, $\Sh(\F^{\M}, \sigma^{\M})$ is not minimal, even when $\F^{\M}$ is, and reducing along a lifted matching does not produce the minimal free resolution of $R/I$ over $\Rbar$. Such edges can still be reduced by Proposition~\ref{thm:A}.
\end{remark}

\section{Additional Examples}\label{sec:moreexamples}
In this section we introduce some other interesting examples that describe behaviors that we could not observe from the running example. 

\begin{remark}\label{rem:defect}
    The transferred maps need not square to zero: by Theorem~\ref{thm:commute}(ii),
    \[
        \sigma^{\M}_{\boldsymbol{e}_j} \sigma^{\M}_{\boldsymbol{e}_j} = p^{\M} \sigma^{(j)} (d h^{\M} + h^{\M} d) \sigma^{(j)} i^{\M} = -(d^{\M} \sigma^{\M}_{2\boldsymbol{e}_j} + \sigma^{\M}_{2\boldsymbol{e}_j} d^{\M}).
    \]
    This is the same thing we saw for the $c = 2$ lattice of Figure~\ref{fig:Sham2_g}, where a general system has $\sigma^{(1)}\sigma^{(2)} + \sigma^{(2)}\sigma^{(1)} = -(d\sigma_{(1,1)} + \sigma_{(1,1)}d)$ and the diagonal arrow of $\sigma_{(1,1)}$ accounts for the difference. Here the arrow of $\sigma^{\M}_{2\boldsymbol{e}_j}$, which jumps two levels, plays that role.
    The transferred system is truncated if
    \begin{equation}\label{eq:trunccond}
        \sigma^{(j)} h^{\M} \sigma^{(k)} = 0 \qquad \text{for all colors } j, k,
    \end{equation}
    since every summand of $\sigma^{\M}_{\boldsymbol{u}}$ with $|\boldsymbol{u}| \ge 2$ contains a factor $\sigma^{(j)} h^{\M} \sigma^{(k)}$.

    If no entry of $\delta$ joins an upper cell to a lower cell, then $h \delta h = 0$, so $h^{\M} = h$ and
    \[
        (\sigma^{(j)} h \sigma^{(k)})(\gamma) = -\sum_{(\alpha, \beta) \in \M} d_{\alpha,\beta}^{-1} \sigma^{(k)}_{\gamma,\beta} \sigma^{(j)}(\alpha),
    \]
    with $\alpha$ upper and $\beta$ lower. So \eqref{eq:trunccond} holds as soon as every pair $(\alpha, \beta) \in \M$ satisfies one of the following.
    \begin{enumerate}[label=(\arabic*),leftmargin=*,itemsep=1pt]
        \item No $\sigma^{(k)}$ has a nonzero entry into $\beta$.
        \item Every $\sigma^{(j)}$ vanishes on $\alpha$.
    \end{enumerate}
    That is, no gradient path in $\Sh(\F, \sigma)$ leaves a cell along a connecting edge, crosses one reversed matched edge, and leaves along another connecting edge. Without $h^{\M} = h$ the path may cross several reversed matched edges, as in the next example.
\end{remark}

\begin{example}\label{ex:defect}
    This example shows that the transferred system need not be truncated. In particular, when we prune too few edges from the free resolution, it is often not truncated. Let $I = (x^3, x^2y, xy^2, y^3) \subset \kk[x, y]$, and let $\F$ be its Taylor resolution on the generators $m_1, \dots, m_4$ in that order. We take $c = 1$ and $a = x^3 = m_1$. Sobieska's homotopy (Section~\ref{sec:sobieska}) for the one term expression $a = m_1$ is the truncated system
    \[
        \sigma(\eps_S) = (-1)^{\#\{s \in S : s > 1\}} \gcd(m_1, m_S) \eps_{S \cup 1} \quad (1 \notin S), \qquad \sigma(\eps_S) = 0 \quad (1 \in S).
    \]

    We match $\eps_{124}$ with $\eps_{14}$ and $\eps_{234}$ with $\eps_{24}$, two edges of weight $-1$. Since $\Grep(\F, \M)$ has no directed cycle, $\M$ is a Morse matching. We leave $(\eps_{123}, \eps_{13})$ and $(\eps_{1234}, \eps_{134})$ unmatched on purpose. Matching them as well would reach the minimal resolution, which has length two, and there the transferred system is truncated by Remark~\ref{rem:length}.

    The critical cells are $\eps_\varnothing$, then $\eps_1, \eps_2, \eps_3, \eps_4$, then $\eps_{12}, \eps_{13}, \eps_{23}, \eps_{34}$, then $\eps_{123}, \eps_{134}$, and $\eps_{1234}$, so $\F^{\M}$ has ranks $1, 4, 4, 2, 1$. In these bases its differential is
    \[
        d^{\M}_1 = \begin{bmatrix} x^3 & x^2y & xy^2 & y^3 \end{bmatrix}, \qquad
        d^{\M}_2 = \begin{bmatrix} y & y^2 & 0 & 0\\ -x & 0 & y & 0\\ 0 & -x^2 & -x & y\\ 0 & 0 & 0 & -x \end{bmatrix},
    \]
    \[
        d^{\M}_3 = \begin{bmatrix} y & -y^2\\ -1 & y\\ x & -xy\\ 0 & 0 \end{bmatrix}, \qquad
        d^{\M}_4 = \begin{bmatrix} y\\ 1 \end{bmatrix}
    \]
    and the transferred system is
    \[
        \sigma^{\M}_1|_0 = \begin{bmatrix} 1\\ 0\\ 0\\ 0 \end{bmatrix}, \quad
        \sigma^{\M}_1|_1 = \begin{bmatrix} 0 & -x^2 & 0 & -y^2\\ 0 & 0 & -x & 0\\ 0 & 0 & 0 & -xy\\ 0 & 0 & 0 & -x^2 \end{bmatrix}, \quad
        \sigma^{\M}_1|_2 = \begin{bmatrix} 0 & 0 & x^2 & 0\\ 0 & 0 & 0 & x \end{bmatrix}, \quad
        \sigma^{\M}_1|_3 = \begin{bmatrix} 0 & x^3 \end{bmatrix},
    \]

    \[
        \sigma^{\M}_2|_1 = \begin{bmatrix} 0 & 0 & 0 & x^3 \end{bmatrix} \colon F^{\M}_1 \to F^{\M}_4, \qquad \sigma^{\M}_2|_q = 0 \quad (q \neq 1).
    \]

    Only one path contributes to $\sigma^{\M}_2 = p^{\M} \sigma h^{\M} \sigma i^{\M}$. Starting from $\eps_4$, the map $\sigma$ goes up to the lower cell $-\eps_{14}$, then $h^{\M}$ goes up to $-\eps_{124} - x \eps_{234}$, and $\sigma$ goes up again, so
    \[
        \sigma^{\M}_2(\eps_4) = \sigma\bigl(-\eps_{124} - x\eps_{234}\bigr) = x^3 \eps_{1234} .
    \]
    This path crosses two reversed matched edges between its two connecting edges. Each matched pair satisfies one of the conditions of Remark~\ref{rem:defect}, (2) at $(\eps_{124}, \eps_{14})$ and (1) at $(\eps_{234}, \eps_{24})$. Those conditions need $h^{\M} = h$, however, and here the entry $x$ of $\delta$ from $\eps_{124}$ to $\eps_{24}$ makes $h^{\M} \neq h$.

    Finally, Theorem~\ref{thm:commute}(ii) gives $(\sigma^{\M}_1)^2 = -(d^{\M} \sigma^{\M}_2 + \sigma^{\M}_2 d^{\M})$, and here both sides are nonzero. They vanish except on $F^{\M}_1$ and $F^{\M}_2$, where
    \[
        (\sigma^{\M}_1)^2|_1 = \begin{bmatrix} 0 & 0 & 0 & -x^3y\\ 0 & 0 & 0 & -x^3 \end{bmatrix}, \qquad
        \bigl(d^{\M} \sigma^{\M}_2 + \sigma^{\M}_2 d^{\M}\bigr)\big|_1 = \begin{bmatrix} 0 & 0 & 0 & x^3y\\ 0 & 0 & 0 & x^3 \end{bmatrix},
    \]
    \[
        (\sigma^{\M}_1)^2|_2 = \begin{bmatrix} 0 & 0 & 0 & x^4 \end{bmatrix}, \qquad
        \bigl(d^{\M} \sigma^{\M}_2 + \sigma^{\M}_2 d^{\M}\bigr)\big|_2 = \begin{bmatrix} 0 & 0 & 0 & -x^4 \end{bmatrix}.
    \]
    Over $\Rbar = R/(x^3)$ all of these entries are zero, so the failure of truncation can only be seen over $R$.
\end{example}

In Examples~\ref{ex:running} and~\ref{ex:runningtwo} the transferred system is truncated, while in Example~\ref{ex:defect} it is not, even though there $\pd_R R/I = 2$. So whether $\sigma^{\M}$ is truncated depends on the matching, and not only on the ideal.

The next two examples are interesting because the reduced free resolution is minimal but does not have a DG algebra structure. The ideal is due to Avramov \cite{Avramov1981}, see also Katth\"an \cite[Example~3.8]{Katthan}. We found a matching in which the transferred system of higher homotopies is truncated despite the lack of DG algebra structure, and one where it is not.
\begin{example}\label{ex:avramovtrunc}
    This example shows that the transferred system can be truncated on a minimal free resolution with no DG algebra structure. Let $I = (xy, yz, zu, uv, vw) \subset \kk[x, y, z, u, v, w]$, and let $\F$ be its Taylor resolution on the generators $m_1, \dots, m_5$ in that order. We take $c = 1$ and $a = z^2u + v^2w = z m_3 + v m_5$. Sobieska's homotopy (Section~\ref{sec:sobieska}) for this two term expression is the truncated system
    \[
        \sigma(\eps_S) = (-1)^{\#\{s \in S : s > 3\}} \, z \gcd(m_3, m_S) \, \eps_{S \cup 3} + v \gcd(m_5, m_S) \, \eps_{S \cup 5},
    \]
    where the first term is omitted if $3 \in S$ and the second if $5 \in S$.

    We use the matching produced by the pruning algorithm \cite[Algorithm~3.1]{AlvarezFernandezGimenez},
    \[
        \M = \{(\eps_{123}, \eps_{13}), (\eps_{234}, \eps_{24}), (\eps_{345}, \eps_{35}), (\eps_{1234}, \eps_{134}), (\eps_{2345}, \eps_{245}), (\eps_{1235}, \eps_{135}), (\eps_{12345}, \eps_{1345})\}.
    \]
    Since $\Grep(\F, \M)$ has no directed cycle, $\M$ is a Morse matching.

    The critical cells are $\eps_\varnothing$, then $\eps_1, \dots, \eps_5$, then $\eps_{12}, \eps_{14}, \eps_{15}, \eps_{23}, \eps_{25}, \eps_{34}, \eps_{45}$, then $\eps_{124}, \eps_{125}, \eps_{145}, \eps_{235}$, and $\eps_{1245}$, so $\F^{\M}$ has ranks $1, 5, 7, 4, 1$. It is the minimal free resolution of $R/I$ over $R$. In these bases its differential is
    \[
        d^{\M}_1 = \begin{bmatrix} xy & yz & zu & uv & vw \end{bmatrix}, \qquad
        d^{\M}_2 = \begin{bmatrix} z & uv & vw & 0 & 0 & 0 & 0\\ -x & 0 & 0 & u & vw & 0 & 0\\ 0 & 0 & 0 & -y & 0 & v & 0\\ 0 & -xy & 0 & 0 & 0 & -z & w\\ 0 & 0 & -xy & 0 & -yz & 0 & -u \end{bmatrix},
    \]
    \[
        d^{\M}_3 = \begin{bmatrix} uv & vw & 0 & 0\\ -z & 0 & w & 0\\ 0 & -z & -u & 0\\ xv & 0 & 0 & vw\\ 0 & x & 0 & -u\\ xy & 0 & 0 & yw\\ 0 & 0 & xy & yz \end{bmatrix}, \qquad
        d^{\M}_4 = \begin{bmatrix} w\\ -u\\ z\\ -x \end{bmatrix}
    \]
    and the transferred system is
    \[
        \sigma^{\M}_1|_0 = \begin{bmatrix} 0\\ 0\\ z\\ 0\\ v \end{bmatrix}, \quad
        \sigma^{\M}_1|_1 = \begin{bmatrix} zu & 0 & 0 & 0 & 0\\ 0 & 0 & 0 & 0 & 0\\ v & 0 & 0 & 0 & 0\\ xz & z^2 & 0 & 0 & 0\\ 0 & v & 0 & 0 & 0\\ 0 & 0 & vw & -zu & -zw\\ 0 & 0 & zv & v^2 & -z^2 \end{bmatrix},
    \]
    \[
        \sigma^{\M}_1|_2 = \begin{bmatrix} 0 & -zu & 0 & 0 & 0 & 0 & 0\\ v & 0 & -zu & 0 & 0 & 0 & 0\\ 0 & v^2 & 0 & 0 & 0 & 0 & 0\\ 0 & 0 & -xz & v & -z^2 & 0 & 0 \end{bmatrix}, \quad
        \sigma^{\M}_1|_3 = \begin{bmatrix} v^2 & 0 & zu & 0 \end{bmatrix},
    \]
    with $\sigma^{\M}_k = 0$ for $k \ge 2$. So the transferred system is truncated, and $(\sigma^{\M}_1)^2 = 0$.

    Since $a \in \mathfrak{m} I$, Remark~\ref{rem:mincond} says that $\Sh(\F^{\M}, \sigma^{\M})$ is the minimal free resolution of $R/I$ over $\Rbar = R/(z^2u + v^2w)$. By Corollary~\ref{cor:ranks} its ranks are $1, 5, 8, 9, 9, 9, \dots$.
\end{example}

\begin{example}\label{ex:avramovnontrunc}
    This example shows that on the same minimal free resolution another matching gives a transferred system that is not truncated. We keep the ideal $I$, the element $a = z^2u + v^2w$ and the homotopy $\sigma$ of Example~\ref{ex:avramovtrunc}, and we change the matching to
    \[
        \M' = \{(\eps_{123}, \eps_{13}), (\eps_{234}, \eps_{24}), (\eps_{345}, \eps_{35}), (\eps_{1234}, \eps_{134}), (\eps_{2345}, \eps_{245}), (\eps_{1345}, \eps_{135}), (\eps_{12345}, \eps_{1235})\}.
    \]
    It differs from $\M$ only in how the four cells $\eps_{135}$, $\eps_{1235}$, $\eps_{1345}$, $\eps_{12345}$ are paired. Since $\Grep(\F, \M')$ has no directed cycle, $\M'$ is a Morse matching.

    The critical cells are those of Example~\ref{ex:avramovtrunc}, and the Morse complex is the same, $\F^{\M'} = \F^{\M}$, with the differential given there. The maps $\sigma^{\M'}_1|_0$ and $\sigma^{\M'}_1|_1$ are also the same as there, but
    \[
        \sigma^{\M'}_1|_2 = \begin{bmatrix} 0 & -zu & -zw & 0 & 0 & 0 & 0\\ v & 0 & 0 & 0 & 0 & 0 & 0\\ 0 & v^2 & -z^2 & 0 & 0 & 0 & 0\\ 0 & 0 & 0 & v & -z^2 & 0 & 0 \end{bmatrix}, \quad
        \sigma^{\M'}_1|_3 = \begin{bmatrix} v^2 & -z^2 & 0 & 0 \end{bmatrix},
    \]

    \[
        \sigma^{\M'}_2|_1 = \begin{bmatrix} zv & 0 & 0 & 0 & 0 \end{bmatrix} \colon F^{\M'}_1 \to F^{\M'}_4, \qquad \sigma^{\M'}_2|_q = 0 \quad (q \neq 1),
    \]
    and $\sigma^{\M'}_k = 0$ for $k \ge 3$.

    Only one path contributes to $\sigma^{\M'}_2 = p^{\M'} \sigma h^{\M'} \sigma i^{\M'}$. Starting from $\eps_1$, the map $\sigma$ goes up to $z \eps_{13} + v \eps_{15}$. The cell $\eps_{15}$ is critical, and $\eps_{13}$ is a lower cell, so $h^{\M'}$ goes up to $z \eps_{123}$, and $\sigma$ goes up again to $zv \, \eps_{1235}$. In $\M'$ the cell $\eps_{1235}$ is a lower cell, matched to $\eps_{12345}$, and $p^{\M'}(\eps_{1235}) = \eps_{1245}$, so
    \[
        \sigma^{\M'}_2(\eps_1) = zv \, \eps_{1245} .
    \]
    In the matching $\M$ of Example~\ref{ex:avramovtrunc} the same path reaches $\eps_{1235}$ as an upper cell, and $p^{\M}$ sends it to $0$.

    Finally, Theorem~\ref{thm:commute}(ii) gives $(\sigma^{\M'}_1)^2 = -(d^{\M'} \sigma^{\M'}_2 + \sigma^{\M'}_2 d^{\M'})$, and here both sides are nonzero. They vanish except on $F^{\M'}_1$ and $F^{\M'}_2$, where
    \[
        (\sigma^{\M'}_1)^2|_1 = \begin{bmatrix} -zvw & 0 & 0 & 0 & 0\\ zuv & 0 & 0 & 0 & 0\\ -z^2v & 0 & 0 & 0 & 0\\ xzv & 0 & 0 & 0 & 0 \end{bmatrix}, \qquad
        (\sigma^{\M'}_1)^2|_2 = \begin{bmatrix} -z^2v & -zuv^2 & -zv^2w & 0 & 0 & 0 & 0 \end{bmatrix}.
    \]

    As in Example~\ref{ex:avramovtrunc}, $\Sh(\F^{\M'}, \sigma^{\M'})$ is the minimal free resolution of $R/I$ over $\Rbar = R/(z^2u + v^2w)$, by Remark~\ref{rem:mincond}. So the same minimal free resolution over $R$ carries a truncated system and one that is not truncated, depending on the matching.
\end{example}

\section*{Declaration of generative AI and AI-assisted technologies in the writing process}
During the preparation of this work the authors used Claude (Anthropic) to make the diagrams and to help spot mathematical errors throughout the paper. After using this tool, the authors reviewed and edited the content as needed and take full responsibility for the content of the publication.

\bibliographystyle{plain}
\bibliography{ref}
\end{document}